\documentclass[a4paper,twoside,10pt]{amsart}
\usepackage{amsmath,amsfonts,amssymb,amsthm}
\newtheorem{theorem}{Theorem}[section]
\usepackage{mathrsfs}
\usepackage{color,graphicx}
\usepackage[a4paper]{geometry}
\numberwithin{equation}{section}

\author[G. Nemes]{Gerg\H{o} Nemes}
\address{Central European University, Department of Mathematics and its Applications, H-1051 Budapest, N\'ador utca 9, Hungary}
\email{nemesgery@gmail.com}


\begin{document}

\title[Resurgence of Hankel and Bessel functions]{The resurgence properties of\\ the Hankel and Bessel functions of\\ nearly equal order and argument}

\begin{abstract} The aim of this paper is to derive new representations for the Hankel functions, the Bessel functions and their derivatives, exploiting the reformulation of the method of steepest descents by M. V. Berry and C. J. Howls (Berry and Howls, Proc. R. Soc. Lond. A \textbf{434} (1991) 657--675). Using these representations, we obtain a number of properties of the asymptotic expansions of the Hankel and Bessel functions and their derivatives of nearly equal order and argument, including explicit and numerically computable error bounds, asymptotics for the late coefficients, exponentially improved asymptotic expansions, and the smooth transition of the Stokes discontinuities.
\end{abstract}
\maketitle

\section{Introduction and main results}\label{section1}

At about the same time, Nicholson \cite{Nicholson} and Debye \cite{Debye} investigated the asymptotic behaviour of the Hankel functions $H_\nu ^{\left( 1 \right)} \left( {z } \right)$, $H_\nu ^{\left( 2 \right)} \left( {z } \right)$, and the Bessel functions $J_\nu  \left( {z} \right)$, $Y_\nu  \left( {z} \right)$ for nearly equal order and argument, both of which were assumed to be positive. In modern notation, the asymptotic series of Nicholson and Debye can be written as
\begin{equation}\label{eq17}
H_\nu ^{\left( 1 \right)} \left( {z } \right) \sim - \frac{2}{{3\pi }}\sum\limits_{n = 0}^\infty  {6^{\frac{{n + 1}}{3}} B_n \left( \kappa  \right)e^{\frac{{2\left( {n + 1} \right)\pi i}}{3}} \sin \left( {\frac{{\left( {n + 1} \right)\pi }}{3}} \right)\frac{{\Gamma \left( {\frac{{n + 1}}{3}} \right)}}{{z^{\frac{{n + 1}}{3}} }}} ,
\end{equation}
\begin{equation}\label{eq18}
H_\nu ^{\left( 2 \right)} \left( {z } \right) \sim - \frac{2}{{3\pi }}\sum\limits_{n = 0}^\infty  {6^{\frac{{n + 1}}{3}} B_n \left( \kappa  \right)e^{-\frac{{2\left( {n + 1} \right)\pi i}}{3}} \sin \left( {\frac{{\left( {n + 1} \right)\pi }}{3}} \right)\frac{{\Gamma \left( {\frac{{n + 1}}{3}} \right)}}{{z^{\frac{{n + 1}}{3}} }}} ,
\end{equation}
\begin{equation}\label{eq19}
J_\nu  \left( {z} \right) \sim \frac{1}{{3\pi }}\sum\limits_{n = 0}^\infty  {6^{\frac{{n + 1}}{3}} B_n \left( \kappa  \right)\sin \left( {\frac{{\left( {n + 1} \right)\pi }}{3}} \right)\frac{{\Gamma \left( {\frac{{n + 1}}{3}} \right)}}{{z^{\frac{{n + 1}}{3}} }}}
\end{equation}
and
\begin{equation}\label{eq20}
Y_\nu  \left( {z } \right) \sim - \frac{2}{{3\pi }}\sum\limits_{n = 0}^\infty  {\left( { - 1} \right)^n 6^{\frac{{n + 1}}{3}} B_n \left( \kappa  \right)\sin ^2 \left( {\frac{{\left( {n + 1} \right)\pi }}{3}} \right)\frac{{\Gamma \left( {\frac{{n + 1}}{3}} \right)}}{{z^{\frac{{n + 1}}{3}} }}}, 
\end{equation}
as $z\to +\infty$ with $\kappa=z-\nu$. For these series, to be asymptotic series in Poincar\'{e}'s sense, it has to be assumed that $\kappa = o\left(z^{1/3}\right)$. The coefficients $B_n \left( \kappa  \right)$ are polynomials in $\kappa$ of degree $n$, the first few being
\[
B_0 \left( \kappa  \right) = 1,\; B_1 \left( \kappa  \right) = \kappa ,\; B_2 \left( \kappa  \right) = \frac{{\kappa ^2 }}{2} - \frac{1}{{20}},\; B_3 \left( \kappa  \right) = \frac{{\kappa ^3 }}{6} - \frac{\kappa }{{15}},\; B_4 \left( \kappa  \right) = \frac{{\kappa ^4 }}{{24}} - \frac{{\kappa ^2 }}{{24}} + \frac{1}{{280}}.
\]
For higher coefficients, see Airey \cite{Airey} and Sch\"{o}be \cite{Schobe}. Watson \cite[pp. 245--248]{Watson} extended these expansions to the sector $\left| {\arg z} \right| \le \pi - \delta  < \pi$, with $0<\delta \leq \pi$ being fixed. In Appendix \ref{appendixb}, we will show that the asymptotic expansions \eqref{eq17} and \eqref{eq18} are valid in the even larger regions  $ - \pi  + \delta \leq \arg z \leq 2\pi  - \delta$ and $ - 2\pi  + \delta  \leq \arg z \leq \pi  - \delta$, respectively. The Stokes phenomenon associated with these series shows that these are the largest possible regions of validity. Another type of asymptotic series, valid when $\kappa = \mathcal{O}\left(z^{1/3}\right)$, were obtained by Olver \cite{Olver} (see also Sch\"{o}be \cite{Schobe}).

The case when $\kappa =0$, i.e., $z=\nu$ was studied in detail in the previous paper \cite{Nemes} of the author.

The aim of this paper is to establish new resurgence-type integral representations for the remainders of the asymptotic expansions \eqref{eq17}--\eqref{eq20}. Our derivation is based on the reformulation of the method of steepest descents by Berry and Howls \cite{Berry} (see also Boyd \cite{Boyd} and Paris \cite[pp. 94--99]{Paris}). Here resurgence has to be understood in the sense of Berry and Howls, meaning that the function reappears in the remainder of its own asymptotic series. We also consider the corresponding asymptotic series for the derivatives. Using these representations, we obtain several new properties of Nicholson's and Debye's classical expansions, including asymptotics for the late coefficients, exponentially improved asymptotic expansions, and the smooth transition of the Stokes discontinuities. As a consequence of the late term formula, we solve a problem of Watson regarding the approximation of the polynomials $B_n \left( \kappa  \right)$.

We also obtain explicit and numerically computable error bounds for these asymptotic series when $\kappa=0$. These bounds are extensions and improvements of the results in the paper \cite{Nemes}. A brief discussion about the general $\kappa$ case is given in Section \ref{section6}.

Our first theorem describes the resurgence properties of the asymptotic expansions of $H_\nu^{\left(1\right)}\left(z\right)$, $H_\nu^{\left(2\right)}\left(z\right)$, $J_\nu \left(z\right)$ and $Y_\nu \left(z\right)$. Throughout this paper, empty sums are taken to be zero.

\begin{theorem}\label{thm1} Let $N$ be a non-negative integer, and $\kappa$ be a complex number such that $\nu+\kappa=z$ and $\left|\Re \left( \kappa  \right)\right| < \frac{N + 1}{3}$. Then we have
\begin{equation}\label{eq12}
H_\nu ^{\left( 1 \right)} \left( {z } \right) =  - \frac{2}{{3\pi }}\sum\limits_{n = 0}^{N - 1} {6^{\frac{{n + 1}}{3}} B_n \left( \kappa  \right)e^{\frac{{2\left( {n + 1} \right)\pi i}}{3}} \sin \left( {\frac{{\left( {n + 1} \right)\pi }}{3}} \right)\frac{{\Gamma \left( {\frac{{n + 1}}{3}} \right)}}{{z^{\frac{{n + 1}}{3}} }}}  + R_N^{\left( H \right)} \left( {z ,\kappa } \right)
\end{equation}
for $-\frac{\pi}{2} < \arg z < \frac{3\pi}{2}$;
\begin{equation}\label{eq7}
H_\nu ^{\left( 2 \right)} \left( {z} \right) =  - \frac{2}{{3\pi }}\sum\limits_{n = 0}^{N - 1} {6^{\frac{{n + 1}}{3}} B_n \left( \kappa  \right)e^{-\frac{{2\left( {n + 1} \right)\pi i}}{3}} \sin \left( {\frac{{\left( {n + 1} \right)\pi }}{3}} \right)\frac{{\Gamma \left( {\frac{{n + 1}}{3}} \right)}}{{z^{\frac{{n + 1}}{3}} }}}  - R_N^{\left( H \right)} \left( {ze^{\pi i} , - \kappa } \right)
\end{equation}
for $-\frac{3\pi}{2} < \arg z < \frac{\pi}{2}$;
\begin{equation}\label{eq13}
J_\nu  \left( z \right) = \frac{1}{{3\pi }}\sum\limits_{n = 0}^{N - 1} {6^{\frac{{n + 1}}{3}} B_n \left( \kappa  \right)\sin \left( {\frac{{\left( {n + 1} \right)\pi }}{3}} \right)\frac{{\Gamma \left( {\frac{{n + 1}}{3}} \right)}}{{z^{\frac{{n + 1}}{3}} }}}  + R_N^{\left( J \right)} \left( {z,\kappa } \right),
\end{equation}
\begin{equation}\label{eq14}
Y_\nu  \left( z \right) =  - \frac{2}{{3\pi }}\sum\limits_{n = 0}^{N - 1} {\left( { - 1} \right)^n 6^{\frac{{n + 1}}{3}} B_n \left( \kappa  \right)\sin ^2 \left( {\frac{{\left( {n + 1} \right)\pi }}{3}} \right)\frac{{\Gamma \left( {\frac{{n + 1}}{3}} \right)}}{{z^{\frac{{n + 1}}{3}} }}}  + R_N^{\left( Y \right)} \left( {z,\kappa } \right)
\end{equation}
for $\left|\arg \nu\right|<\frac{\pi}{2}$. The coefficients $B_n \left( \kappa  \right)$ are given by
\begin{gather}\label{eq8}
\begin{split}
B_n \left( \kappa  \right) = \; & \frac{1}{{n!}}\left[ {\frac{{d^n }}{{dt^n }}\left( {e^{\kappa t} \left( {\frac{1}{6}\frac{{t^3 }}{{\sinh t - t}}} \right)^{\frac{{n + 1}}{3}} } \right)} \right]_{t = 0} \\ = \; & \frac{{6^{ - \frac{{n + 1}}{3}} }}{{2\Gamma \left( {\frac{{n + 1}}{3}} \right)}}\int_0^{ + \infty } {t^{\frac{{n - 2}}{3}} e^{ - 2\pi t} i\left( {e^{\left( {2\pi \kappa  - \frac{\pi }{2}n} \right)i} H_{it + \kappa }^{\left( 1 \right)} \left( {it} \right) + e^{ - \left( {2\pi \kappa  - \frac{\pi }{2}n} \right)i} H_{it - \kappa }^{\left( 1 \right)} \left( {it} \right)} \right)dt} ,
\end{split}
\end{gather}
with the second representation being true if $\left|\Re \left( \kappa  \right)\right| < \frac{n + 1}{3}$. The remainder terms can be expressed as
\begin{gather}\label{eq9}
\begin{split}
R_N^{\left( H \right)} \left( {z ,\kappa } \right) = \; & \frac{e^{\left( {2\pi \kappa  - \frac{\pi }{2}N} \right)i}}{{6\pi z^{\frac{{N + 1}}{3}} }}\int_0^{ + \infty } {t^{\frac{{N - 2}}{3}}e^{ - 2\pi t} \left( {\frac{{e^{\frac{{\left( {N + 1} \right)\pi i}}{3}} }}{{1 + i\left( {t/z} \right)^{\frac{1}{3}} e^{\frac{\pi }{3}i} }} - \frac{{e^{\left( {N + 1} \right)\pi i} }}{{1 - i\left( {t/z} \right)^{\frac{1}{3}} }}} \right) H_{it + \kappa }^{\left( 1 \right)} \left( {it} \right)dt} 
\\ & + \frac{e^{-\left( {2\pi \kappa  - \frac{\pi }{2}N} \right)i}}{{6\pi z^{\frac{{N + 1}}{3}} }}\int_0^{ + \infty } { t^{\frac{{N - 2}}{3}}e^{ - 2\pi t} \left( {\frac{{e^{\frac{{\left( {N + 1} \right)\pi i}}{3}} }}{{1 - i\left( {t/z} \right)^{\frac{1}{3}} e^{\frac{\pi }{3}i} }} - \frac{{e^{\left( {N + 1} \right)\pi i} }}{{1 + i\left( {t/z} \right)^{\frac{1}{3}} }}} \right) H_{it - \kappa }^{\left( 1 \right)} \left( {it} \right)dt} ,
\end{split}
\end{gather}
\begin{gather}\label{eq15}
\begin{split}
R_N^{\left( J \right)} \left( {z,\kappa } \right) = \; & \frac{{e^{\left( {2\pi \kappa  - \frac{\pi }{2}N} \right)i} }}{{12\pi z^{\frac{{N + 1}}{3}} }}\int_0^{ + \infty } {t^{\frac{{N - 2}}{3}} e^{ - 2\pi t} \left( {\frac{{e^{\frac{{\left( {N + 1} \right)\pi i}}{3}} }}{{1 + i\left( {t/z} \right)^{\frac{1}{3}} e^{\frac{\pi }{3}i} }} - \frac{{e^{ - \frac{{\left( {N + 1} \right)\pi i}}{3}} }}{{1 + i\left( {t/z} \right)^{\frac{1}{3}} e^{ - \frac{\pi }{3}i} }}} \right)H_{it + \kappa }^{\left( 1 \right)} \left( {it} \right)dt} \\ & + \frac{{e^{ - \left( {2\pi \kappa  - \frac{\pi }{2}N} \right)i} }}{{12\pi z^{\frac{{N + 1}}{3}} }}\int_0^{ + \infty } {t^{\frac{{N - 2}}{3}} e^{ - 2\pi t} \left( {\frac{{e^{\frac{{\left( {N + 1} \right)\pi i}}{3}} }}{{1 - i\left( {t/z} \right)^{\frac{1}{3}} e^{\frac{\pi }{3}i} }} - \frac{{e^{ - \frac{{\left( {N + 1} \right)\pi i}}{3}} }}{{1 - i\left( {t/z} \right)^{\frac{1}{3}} e^{ - \frac{\pi }{3}i} }}} \right)H_{it - \kappa }^{\left( 1 \right)} \left( {it} \right)dt}
\end{split}
\end{gather}
and
\begin{gather}\label{eq16}
\begin{split}
& R_N^{\left( Y \right)} \left( {z,\kappa } \right) =\frac{{e^{\left( {2\pi \kappa  - \frac{\pi }{2}N} \right)i} }}{{12\pi iz^{\frac{{N + 1}}{3}} }}\int_0^{ + \infty } {t^{\frac{{N - 2}}{3}} e^{ - 2\pi t} \left( {\frac{{e^{\frac{{\left( {N + 1} \right)\pi i}}{3}} }}{{1 + i\left( {t/z} \right)^{\frac{1}{3}} e^{\frac{\pi }{3}i} }} + \frac{{e^{ - \frac{{\left( {N + 1} \right)\pi i}}{3}} }}{{1 + i\left( {t/z} \right)^{\frac{1}{3}} e^{ - \frac{\pi }{3}i} }} - \frac{{2e^{\left( {N + 1} \right)\pi i} }}{{1 - i\left( {t/z} \right)^{\frac{1}{3}} }}} \right)H_{it + \kappa }^{\left( 1 \right)} \left( {it} \right)dt} 
\\ & + \frac{{e^{ - \left( {2\pi \kappa  - \frac{\pi }{2}N} \right)i} }}{{12\pi iz^{\frac{{N + 1}}{3}} }}\int_0^{ + \infty } {t^{\frac{{N - 2}}{3}} e^{ - 2\pi t} \left( {\frac{{e^{\frac{{\left( {N + 1} \right)\pi i}}{3}} }}{{1 - i\left( {t/z} \right)^{\frac{1}{3}} e^{\frac{\pi }{3}i} }} + \frac{{e^{ - \frac{{\left( {N + 1} \right)\pi i}}{3}} }}{{1 - i\left( {t/z} \right)^{\frac{1}{3}} e^{ - \frac{\pi }{3}i} }} - \frac{{2e^{\left( {N + 1} \right)\pi i} }}{{1 + i\left( {t/z} \right)^{\frac{1}{3}} }}} \right)H_{it - \kappa }^{\left( 1 \right)} \left( {it} \right)dt} .
\end{split}
\end{gather}
The cube roots are defined to be positive on the positive real line and are defined by analytic continuation elsewhere.
\end{theorem}

We remark that an equivalent form of the first representation in \eqref{eq8} was also given by Watson \cite[p. 246]{Watson}. Some other formulas for the coefficients
$B_n \left( \kappa  \right)$ can be found in Appendix \ref{appendixa}.

Using the continuation formulas (see, e.g., \cite[\S 10.11]{NIST})
\begin{gather}\label{eq75}
\begin{split}
\sin \left( {\pi \nu } \right)H_{\nu e^{2\pi im} }^{\left( 1 \right)} \left( {z e^{2\pi im} } \right) & = \sin \left( {\pi \nu } \right)H_\nu ^{\left( 1 \right)} \left( {z e^{2\pi im} } \right)\\ & =  - \sin \left( {\left( {2m - 1} \right)\pi \nu } \right)H_\nu ^{\left( 1 \right)} \left( z \right) - e^{ - \pi i\nu } \sin \left( {2\pi m\nu } \right)H_\nu ^{\left( 2 \right)} \left( z \right),
\end{split}
\end{gather}
\begin{align*}
\sin \left( {\pi \nu } \right)H_{\nu e^{2\pi im} }^{\left( 2 \right)} \left( {z e^{2\pi im} } \right) & = \sin \left( {\pi \nu } \right)H_\nu ^{\left( 2 \right)} \left( {z e^{2\pi im} } \right)\\ & = \sin \left( {\left( {2m + 1} \right)\pi \nu } \right)H_\nu ^{\left( 2 \right)} \left( z \right) + e^{\pi i\nu } \sin \left( {2\pi m\nu } \right)H_\nu ^{\left( 1 \right)} \left( z \right),
\end{align*}
and the resurgence formulas \eqref{eq12}, \eqref{eq7} and \eqref{eq9}, we can derive analogous representations in sectors of the form
\[
\left( {2m - \frac{1}{2}} \right)\pi  < \arg z  < \left( {2m + \frac{3}{2}} \right)\pi \; \text{ and } \; \left( {2m - \frac{3}{2}} \right)\pi  < \arg z < \left( {2m + \frac{1}{2}} \right)\pi, \; m \in \mathbb{Z},
\]
respectively. The lines $\arg z  = \left( {2m - \frac{1}{2}} \right)\pi$ are the Stokes lines for the function $H_\nu ^{\left( 1 \right)} \left(z\right)$, and the lines $\arg z = \left( {2m + \frac{1}{2}} \right)\pi$ are the Stokes lines for the function $H_\nu ^{\left( 2 \right)} \left( z \right)$.

Similarly, applying the continuation formulas
\begin{align*}
& J_{\nu e^{\left( {2m + 1} \right)\pi i} } \left( {z e^{\left( {2m + 1} \right)\pi i} } \right) = J_{ - \nu } \left( {z e^{\left( {2m + 1} \right)\pi i} } \right)\\
& =  e^{2\pi im\nu } J_\nu  \left( z \right) - i\sin \left( {2\pi m\nu } \right)H_\nu ^{\left( 1 \right)} \left( z \right)  - ie^{ - \pi i\nu } \sin \left( {\left( {2m + 1} \right)\nu \pi } \right)H_\nu ^{\left( 2 \right)} \left( z \right),
\end{align*}
\begin{align*}
& Y_{\nu e^{\left( {2m + 1} \right)\pi i} } \left( {z e^{\left( {2m + 1} \right)\pi i} } \right) = Y_{ - \nu } \left( {z e^{\left( {2m + 1} \right)\pi i}} \right) \\
& = e^{ - 2\left( {m + 1} \right)\pi i\nu } Y_\nu  \left( z \right) + 2ie^{ - \pi i\nu } \sin \left( {\left( {2m + 1} \right)\pi \nu } \right)\cot \left( {\pi \nu } \right)J_\nu  \left( z \right) \\
 & \;\; - \sin \left( {2\pi m\nu } \right)H_\nu ^{\left( 1 \right)} \left( z \right) - e^{ - \pi i\nu } \sin \left( {\left( {2m + 1} \right)\pi \nu } \right)H_\nu ^{\left( 2 \right)} \left( z \right)
\end{align*}
and the representations \eqref{eq12}--\eqref{eq16}, we can obtain analogous formulas in any sector of the form
\[
\left( {2m + \frac{1}{2}} \right)\pi  < \arg z  < \left( {2m + \frac{3}{2}} \right)\pi , \; m \in \mathbb{Z}.
\]
Finally, from
\[
J_{\nu e^{2\pi im} } \left( {z e^{2\pi im}} \right) = J_\nu  \left( {z e^{2\pi im} } \right) = e^{2\pi im\nu } J_\nu  \left( z \right),
\]
\[
Y_{\nu e^{2\pi im} } \left( {z e^{2\pi im} } \right) = Y_\nu  \left( {z e^{2\pi im}} \right) = e^{ - 2\pi im\nu } Y_\nu  \left( z \right) + 2i\sin \left( {2\pi m\nu } \right)\cot \left( {\pi \nu } \right)J_\nu  \left( z \right)
\]
and the resurgence formulas \eqref{eq13}, \eqref{eq14}, \eqref{eq15} and \eqref{eq16}, we can derive the corresponding representations in sectors of the form
\[
\left( {2m - \frac{1}{2}} \right)\pi  < \arg z  < \left( {2m + \frac{1}{2}} \right)\pi , \; m \in \mathbb{Z}.
\]
The lines $\arg z  = \left( {2m \pm \frac{1}{2}} \right)\pi$ are the Stokes lines for the functions $J_\nu\left( z \right)$ and $Y_\nu \left( z \right)$.

When $\nu$ is an integer, the limiting values have to be taken in these continuation formulas.

The second theorem provides resurgence formulas for the $z$-derivatives of the Hankel and Bessel functions. The results are direct consequences of Theorem \ref{thm1} and the functional relation $2\mathscr{C}'_\nu  \left( z \right) = \mathscr{C}_{\nu  - 1} \left( z \right) - \mathscr{C}_{\nu  + 1} \left( z \right)$, where $\mathscr{C}_\nu  \left( z \right)$ denotes any of the functions $H_\nu^{\left(1\right)}\left(z\right)$, $H_\nu^{\left(2\right)}\left(z\right)$, $J_\nu \left(z\right)$ or $Y_\nu \left(z\right)$.

\begin{theorem}\label{thm2} Let $N>2$ be an integer, and $\kappa$ be a complex number such that $\nu+\kappa=z$ and $\left|\Re \left( \kappa  \right)\right| +1 < \frac{N + 1}{3}$. Then we have
\[
H_\nu ^{\left( 1 \right) \prime}  \left( z \right) =  - \frac{2}{{3\pi }}\sum\limits_{n = 1}^{N - 1} {6^{\frac{{n + 1}}{3}} D_n \left( \kappa  \right)e^{\frac{{2\left( {n + 1} \right)\pi i}}{3}} \sin \left( {\frac{{\left( {n + 1} \right)\pi }}{3}} \right)\frac{{\Gamma \left( {\frac{{n + 1}}{3}} \right)}}{{z^{\frac{{n + 1}}{3}} }}}  + R_N^{\left( {H'} \right)} \left( {z,\kappa } \right)
\]
for $-\frac{\pi}{2} < \arg z < \frac{3\pi}{2}$;
\[
H_\nu ^{\left( 2 \right) \prime} \left( z \right) =  - \frac{2}{{3\pi }}\sum\limits_{n = 1}^{N - 1} {6^{\frac{{n + 1}}{3}} D_n \left( \kappa  \right)e^{ - \frac{{2\left( {n + 1} \right)\pi i}}{3}} \sin \left( {\frac{{\left( {n + 1} \right)\pi }}{3}} \right)\frac{{\Gamma \left( {\frac{{n + 1}}{3}} \right)}}{{z^{\frac{{n + 1}}{3}} }}}  + R_N^{\left( {H'} \right)} \left( {ze^{\pi i} , - \kappa } \right)
\]
for $-\frac{3\pi}{2} < \arg z < \frac{\pi}{2}$;
\[
J'_\nu  \left( z \right) = \frac{1}{{3\pi }}\sum\limits_{n = 1}^{N - 1} {6^{\frac{{n + 1}}{3}} D_n \left( \kappa  \right)\sin \left( {\frac{{\left( {n + 1} \right)\pi }}{3}} \right)\frac{{\Gamma \left( {\frac{{n + 1}}{3}} \right)}}{{z^{\frac{{n + 1}}{3}} }}}  + R_N^{\left( {J'} \right)} \left( {z,\kappa } \right),
\]
\[
Y'_\nu  \left( z \right) =  - \frac{2}{{3\pi }}\sum\limits_{n = 1}^{N - 1} {\left( { - 1} \right)^n 6^{\frac{{n + 1}}{3}} D_n \left( \kappa  \right)\sin ^2 \left( {\frac{{\left( {n + 1} \right)\pi }}{3}} \right)\frac{{\Gamma \left( {\frac{{n + 1}}{3}} \right)}}{{z^{\frac{{n + 1}}{3}} }}}  + R_N^{\left( {Y'} \right)} \left( {z,\kappa } \right)
\]
for $\left|\arg \nu\right|<\frac{\pi}{2}$. The coefficients $D_n \left( \kappa  \right)$ are given by
\begin{align*}
D_n \left( \kappa  \right) & = \frac{{B_n \left( {\kappa  + 1} \right) - B_n \left( {\kappa  - 1} \right)}}{2}
\\ & = \frac{1}{{n!}}\left[ {\frac{{d^n }}{{dt^n }}\left( {e^{\kappa t} \sinh t\left( {\frac{1}{6}\frac{{t^3 }}{{\sinh t - t}}} \right)^{\frac{{n + 1}}{3}} } \right)} \right]_{t = 0} 
\\ & = - \frac{{6^{ - \frac{{n + 1}}{3}} }}{{2\Gamma \left( {\frac{{n + 1}}{3}} \right)}}\int_0^{ + \infty } {t^{\frac{{n - 2}}{3}} e^{ - 2\pi t} i\left( {e^{\left( {2\pi \kappa  - \frac{\pi }{2}n} \right)i} H_{it + \kappa }^{\left( 1 \right) \prime} \left( {it} \right) - e^{ - \left( {2\pi \kappa  - \frac{\pi }{2}n} \right)i} H_{it - \kappa }^{\left( 1 \right) \prime}  \left( {it} \right)} \right)dt} ,
\end{align*}
with the second representation being true if $\left|\Re \left( \kappa  \right)\right|+1 < \frac{n + 1}{3}$. The remainder terms can be expressed as
\begin{gather}\label{eq32}
\begin{split}
R_N^{\left( {H'} \right)} \left( {z,\kappa } \right) = & - \frac{{e^{\left( {2\pi \kappa  - \frac{\pi }{2}N} \right)i} }}{{6\pi z^{\frac{{N + 1}}{3}} }}\int_0^{ + \infty } {t^{\frac{{N - 2}}{3}} e^{ - 2\pi t} \left( {\frac{{e^{\frac{{\left( {N + 1} \right)\pi i}}{3}} }}{{1 + i\left( {t/z} \right)^{\frac{1}{3}} e^{\frac{\pi }{3}i} }} - \frac{{e^{\left( {N + 1} \right)\pi i} }}{{1 - i\left( {t/z} \right)^{\frac{1}{3}} }}} \right)H_{it + \kappa }^{\left( 1 \right) \prime} \left( {it} \right)dt} \\ & + \frac{{e^{ - \left( {2\pi \kappa  - \frac{\pi }{2}N} \right)i} }}{{6\pi z^{\frac{{N + 1}}{3}} }}\int_0^{ + \infty } {t^{\frac{{N - 2}}{3}} e^{ - 2\pi t} \left( {\frac{{e^{\frac{{\left( {N + 1} \right)\pi i}}{3}} }}{{1 - i\left( {t/z} \right)^{\frac{1}{3}} e^{\frac{\pi }{3}i} }} - \frac{{e^{\left( {N + 1} \right)\pi i} }}{{1 + i\left( {t/z} \right)^{\frac{1}{3}} }}} \right)H_{it - \kappa }^{\left( 1 \right) \prime}  \left( {it} \right)dt} ,
\end{split}
\end{gather}
\begin{gather}\label{eq65}
\begin{split}
R_N^{\left( {J'} \right)} \left( {z,\kappa } \right) =  & - \frac{{e^{\left( {2\pi \kappa  - \frac{\pi }{2}N} \right)i} }}{{12\pi z^{\frac{{N + 1}}{3}} }}\int_0^{ + \infty } {t^{\frac{{N - 2}}{3}} e^{ - 2\pi t} \left( {\frac{{e^{\frac{{\left( {N + 1} \right)\pi i}}{3}} }}{{1 + i\left( {t/z} \right)^{\frac{1}{3}} e^{\frac{\pi }{3}i} }} - \frac{{e^{ - \frac{{\left( {N + 1} \right)\pi i}}{3}} }}{{1 + i\left( {t/z} \right)^{\frac{1}{3}} e^{ - \frac{\pi }{3}i} }}} \right)H_{it + \kappa }^{\left( {1} \right)\prime} \left( {it} \right)dt} \\ & + \frac{{e^{ - \left( {2\pi \kappa  - \frac{\pi }{2}N} \right)i} }}{{12\pi z^{\frac{{N + 1}}{3}} }}\int_0^{ + \infty } {t^{\frac{{N - 2}}{3}} e^{ - 2\pi t} \left( {\frac{{e^{\frac{{\left( {N + 1} \right)\pi i}}{3}} }}{{1 - i\left( {t/z} \right)^{\frac{1}{3}} e^{\frac{\pi }{3}i} }} - \frac{{e^{ - \frac{{\left( {N + 1} \right)\pi i}}{3}} }}{{1 - i\left( {t/z} \right)^{\frac{1}{3}} e^{ - \frac{\pi }{3}i} }}} \right)H_{it - \kappa }^{\left( {1} \right)\prime} \left( {it} \right)dt} 
\end{split}
\end{gather}
and
\begin{gather}\label{eq67}
\begin{split}
& R_N^{\left( Y' \right)} \left( {z,\kappa } \right) = - \frac{{e^{\left( {2\pi \kappa  - \frac{\pi }{2}N} \right)i} }}{{12\pi iz^{\frac{{N + 1}}{3}} }}\int_0^{ + \infty } {t^{\frac{{N - 2}}{3}} e^{ - 2\pi t} \left( {\frac{{e^{\frac{{\left( {N + 1} \right)\pi i}}{3}} }}{{1 + i\left( {t/z} \right)^{\frac{1}{3}} e^{\frac{\pi }{3}i} }} + \frac{{e^{ - \frac{{\left( {N + 1} \right)\pi i}}{3}} }}{{1 + i\left( {t/z} \right)^{\frac{1}{3}} e^{ - \frac{\pi }{3}i} }} - \frac{{2e^{\left( {N + 1} \right)\pi i} }}{{1 - i\left( {t/z} \right)^{\frac{1}{3}} }}} \right)H_{it + \kappa }^{\left( 1 \right) \prime} \left( {it} \right)dt} 
\\ & + \frac{{e^{ - \left( {2\pi \kappa  - \frac{\pi }{2}N} \right)i} }}{{12\pi iz^{\frac{{N + 1}}{3}} }}\int_0^{ + \infty } {t^{\frac{{N - 2}}{3}} e^{ - 2\pi t} \left( {\frac{{e^{\frac{{\left( {N + 1} \right)\pi i}}{3}} }}{{1 - i\left( {t/z} \right)^{\frac{1}{3}} e^{\frac{\pi }{3}i} }} + \frac{{e^{ - \frac{{\left( {N + 1} \right)\pi i}}{3}} }}{{1 - i\left( {t/z} \right)^{\frac{1}{3}} e^{ - \frac{\pi }{3}i} }} - \frac{{2e^{\left( {N + 1} \right)\pi i} }}{{1 + i\left( {t/z} \right)^{\frac{1}{3}} }}} \right)H_{it - \kappa }^{\left( 1 \right) \prime} \left( {it} \right)dt} .
\end{split}
\end{gather}
The cube roots are defined to be positive on the positive real line and are defined by analytic continuation elsewhere.
\end{theorem}

Again, these resurgence formulas can be extended to other sectors of the complex plane using the various connection formulas between the Hankel and Bessel functions.

If we neglect the remainder terms and extend the sums to $N = \infty$ in Theorem \ref{thm2}, we obtain the asymptotic expansions of the derivatives of the Hankel and Bessel functions of nearly equal order and argument.

In the following theorem, we give exponentially improved asymptotic expansions for the functions $H_\nu ^{\left( 1 \right)} \left( z \right)$ and $H_\nu ^{\left( 1 \right)\prime} \left( z \right)$. The related expansions for the functions $H_\nu^{\left(2\right)}\left(z\right)$, $J_\nu \left(z\right)$, $Y_\nu \left(z\right)$ and their derivatives may be derived from the corresponding connection formulas. In this theorem we truncate the asymptotic series of $H_\nu ^{\left( 1 \right)} \left( z \right)$ and $H_\nu ^{\left( 1 \right)\prime} \left( z \right)$ at about their least terms and re-expand the remainders into new asymptotic expansions. The resulting exponentially improved asymptotic series are valid in larger regions than the original Nicholson--Debye expansions. The terms in these new series involve the Terminant function $\widehat T_p\left(w\right)$, which allows the smooth transition through the Stokes lines $\arg z = -\frac{\pi}{2}$ and $\arg z = \frac{3\pi}{2}$. For the definition and basic properties of the Terminant function, see Section \ref{section5}. Throughout this paper, we use subscripts in the $\mathcal{O}$ notations to indicate the dependence of the implied constant on certain parameters.

\begin{theorem}\label{thm3} Let $\kappa$ be a fixed complex number such that $\nu+\kappa=z$. Suppose that $ - 2\pi  <  - 2\pi  + \delta  \le \arg z \le 3\pi  - \delta  < 3\pi$ with a small positive $\delta$. Define $R_{N,M}^{\left( {H} \right)} \left( {z,\kappa } \right)$ and $R_{N,M}^{\left( {H'} \right)} \left( {z,\kappa } \right)$ by
\begin{multline*}
H_\nu ^{\left( 1 \right)} \left( z \right) = \frac{{e^{ - \frac{\pi }{3}i} }}{{\sqrt 3 \pi z^{\frac{1}{3}} }}\sum\limits_{n = 0}^{N - 1} {\left( { - 1} \right)^n 6^{n + \frac{1}{3}} B_{3n} \left( \kappa  \right)\frac{{\Gamma \left( {n + \frac{1}{3}} \right)}}{{z^n }}} \\ + \frac{{e^{\frac{\pi }{3}i} }}{{\sqrt 3 \pi z^{\frac{2}{3}} }}\sum\limits_{m = 0}^{M - 1} {\left( { - 1} \right)^m 6^{m + \frac{2}{3}} B_{3m + 1} \left( \kappa  \right)\frac{{\Gamma \left( {m + \frac{2}{3}} \right)}}{{z^m }}}  + R_{N,M}^{\left( {H} \right)} \left( {z,\kappa } \right),
\end{multline*}
and
\begin{multline*}
H_\nu ^{\left( 1 \right)\prime} \left( z \right) = \frac{{e^{ - \frac{\pi }{3}i} }}{{\sqrt 3 \pi z^{\frac{1}{3}} }}\sum\limits_{n = 0}^{N - 1} {\left( { - 1} \right)^n 6^{n + \frac{1}{3}} D_{3n} \left( \kappa  \right)\frac{{\Gamma \left( {n + \frac{1}{3}} \right)}}{{z^n }}} \\ + \frac{{e^{\frac{\pi }{3}i} }}{{\sqrt 3 \pi z^{\frac{2}{3}} }}\sum\limits_{m = 0}^{M - 1} {\left( { - 1} \right)^m 6^{m + \frac{2}{3}} D_{3m + 1} \left( \kappa  \right)\frac{{\Gamma \left( {m + \frac{2}{3}} \right)}}{{z^m }}}  + R_{N,M}^{\left( {H'} \right)} \left( {z,\kappa } \right),
\end{multline*}
where
\[
N = 2\pi \left|z\right| +\rho \; \text{ and } \; M = 2\pi \left|z\right| +\sigma,
\]
$\left|z\right|$ being large, $\rho$ and $\sigma$ being bounded quantities. Then
\begin{align*}
R_{N,M}^{\left( {H} \right)} \left( {z,\kappa } \right) = \; & ie^{ - \frac{\pi }{3}i} \frac{{e^{ - 2\pi i\nu } }}{{\sqrt 3 }}\frac{2}{{3\pi }}\sum\limits_{k = 0}^{K - 1} {6^{\frac{{k + 1}}{3}} B_k \left( \kappa  \right)\sin \left( {\frac{{\left( {k + 1} \right)\pi }}{3}} \right)\frac{{\Gamma \left( {\frac{{k + 1}}{3}} \right)}}{{z^{\frac{{k + 1}}{3}} }}\widehat T_{N - \frac{k}{3}} \left( { - 2\pi iz} \right)} \\
& - i\frac{{e^{2\pi i\nu } }}{{\sqrt 3 }}\frac{2}{{3\pi }}\sum\limits_{k = 0}^{K - 1} {6^{\frac{{k + 1}}{3}} B_k \left( \kappa  \right)e^{\frac{{2\left( {k + 1} \right)\pi i}}{3}} \sin \left( {\frac{{\left( {k + 1} \right)\pi }}{3}} \right)\frac{{\Gamma \left( {\frac{{k + 1}}{3}} \right)}}{{z^{\frac{{k + 1}}{3}} }}\widehat T_{N - \frac{k}{3}} \left( {2\pi iz} \right)} \\
& - ie^{\frac{\pi }{3}i} \frac{{e^{ - 2\pi i\nu } }}{{\sqrt 3 }}\frac{2}{{3\pi }}\sum\limits_{\ell  = 0}^{L - 1} {6^{\frac{{\ell  + 1}}{3}} B_\ell  \left( \kappa  \right)\sin \left( {\frac{{\left( {\ell  + 1} \right)\pi }}{3}} \right)\frac{{\Gamma \left( {\frac{{\ell  + 1}}{3}} \right)}}{{z^{\frac{{\ell  + 1}}{3}} }}\widehat T_{M - \frac{{\ell  - 1}}{3}} \left( { - 2\pi iz} \right)} \\
& + i\frac{{e^{2\pi i\nu } }}{{\sqrt 3 }}\frac{2}{{3\pi }}\sum\limits_{\ell  = 0}^{L - 1} {6^{\frac{{\ell  + 1}}{3}} B_\ell  \left( \kappa  \right)e^{\frac{{2\left( {\ell  + 1} \right)\pi i}}{3}} \sin \left( {\frac{{\left( {\ell  + 1} \right)\pi }}{3}} \right)\frac{{\Gamma \left( {\frac{{\ell  + 1}}{3}} \right)}}{{z^{\frac{{\ell  + 1}}{3}} }}\widehat T_{M - \frac{{\ell  - 1}}{3}} \left( {2\pi iz} \right)} \\
& + R_{N,M,K,L}^{\left( {H} \right)} \left( {z,\kappa } \right),
\end{align*}
and
\begin{align*}
R_{N,M}^{\left( {H'} \right)} \left( {z,\kappa } \right) = \; & ie^{ - \frac{\pi }{3}i} \frac{{e^{ - 2\pi i\nu } }}{{\sqrt 3 }}\frac{2}{{3\pi }}\sum\limits_{k = 0}^{K - 1} {6^{\frac{{k + 1}}{3}} D_k \left( \kappa  \right)\sin \left( {\frac{{\left( {k + 1} \right)\pi }}{3}} \right)\frac{{\Gamma \left( {\frac{{k + 1}}{3}} \right)}}{{z^{\frac{{k + 1}}{3}} }}\widehat T_{N - \frac{k}{3}} \left( { - 2\pi iz} \right)} \\
& - i\frac{{e^{2\pi i\nu } }}{{\sqrt 3 }}\frac{2}{{3\pi }}\sum\limits_{k = 0}^{K - 1} {6^{\frac{{k + 1}}{3}} D_k \left( \kappa  \right)e^{\frac{{2\left( {k + 1} \right)\pi i}}{3}} \sin \left( {\frac{{\left( {k + 1} \right)\pi }}{3}} \right)\frac{{\Gamma \left( {\frac{{k + 1}}{3}} \right)}}{{z^{\frac{{k + 1}}{3}} }}\widehat T_{N - \frac{k}{3}} \left( {2\pi iz} \right)} \\
& - ie^{\frac{\pi }{3}i} \frac{{e^{ - 2\pi i\nu } }}{{\sqrt 3 }}\frac{2}{{3\pi }}\sum\limits_{\ell  = 0}^{L - 1} {6^{\frac{{\ell  + 1}}{3}} D_\ell  \left( \kappa  \right)\sin \left( {\frac{{\left( {\ell  + 1} \right)\pi }}{3}} \right)\frac{{\Gamma \left( {\frac{{\ell  + 1}}{3}} \right)}}{{z^{\frac{{\ell  + 1}}{3}} }}\widehat T_{M - \frac{{\ell  - 1}}{3}} \left( { - 2\pi iz} \right)} \\
& + i\frac{{e^{2\pi i\nu } }}{{\sqrt 3 }}\frac{2}{{3\pi }}\sum\limits_{\ell  = 0}^{L - 1} {6^{\frac{{\ell  + 1}}{3}} D_\ell  \left( \kappa  \right)e^{\frac{{2\left( {\ell  + 1} \right)\pi i}}{3}} \sin \left( {\frac{{\left( {\ell  + 1} \right)\pi }}{3}} \right)\frac{{\Gamma \left( {\frac{{\ell  + 1}}{3}} \right)}}{{z^{\frac{{\ell  + 1}}{3}} }}\widehat T_{M - \frac{{\ell  - 1}}{3}} \left( {2\pi iz} \right)} \\
& + R_{N,M,K,L}^{\left( {H'} \right)} \left( {z,\kappa } \right),
\end{align*}
where $K$ and $L$ are arbitrary fixed non-negative integers and
\begin{equation}\label{eq60}
R_{N,M,K,L}^{\left( {H} \right)} \left( {z,\kappa } \right), R_{N,M,K,L}^{\left( {H'} \right)} \left( {z,\kappa } \right) = \mathcal{O}_{\kappa ,K,\rho } \left( {\frac{{e^{ - 2\pi \left| z \right|} }}{{\left| z \right|^{\frac{{K + 1}}{3}} }}} \right) + \mathcal{O}_{\kappa ,L,\sigma } \left( {\frac{{e^{ - 2\pi \left| z \right|} }}{{\left| z \right|^{\frac{{L + 1}}{3}} }}} \right)
\end{equation}
for $\left|\arg z\right|\leq \frac{\pi}{2}$;
\begin{equation}\label{eq76}
R_{N,M,K,L}^{\left( H \right)} \left( {z,\kappa } \right), R_{N,M,K,L}^{\left( {H'} \right)} \left( {z,\kappa } \right) = \mathcal{O}_{\kappa ,K,\rho } \left( {\frac{{e^{\mp 2\pi \Im \left( z \right)} }}{{\left| z \right|^{\frac{{K + 1}}{3}} }}} \right) + \mathcal{O}_{\kappa ,L,\sigma } \left( {\frac{{e^{\mp 2\pi \Im \left( z \right)} }}{{\left| z \right|^{\frac{{L + 1}}{3}} }}} \right)
\end{equation}
for $\frac{\pi}{2} \leq \pm \arg z \leq \frac{3\pi}{2}$;
\[
R_{N,M,K,L}^{\left( H \right)} \left( {z,\kappa } \right), R_{N,M,K,L}^{\left( {H'} \right)} \left( {z,\kappa } \right) = \mathcal{O}_{\kappa ,K,\rho ,\delta } \left( {\frac{{\cosh \left( {2\pi \Im \left( z \right)} \right)}}{{\left| z \right|^{\frac{{K + 1}}{3}} }}} \right) + \mathcal{O}_{\kappa ,L,\sigma ,\delta } \left( {\frac{{\cosh \left( {2\pi \Im \left( z \right)} \right)}}{{\left| z \right|^{\frac{{L + 1}}{3}} }}} \right)
\]
for $ - 2\pi  <  - 2\pi  + \delta  \le \arg z \le  - \frac{{3\pi }}{2}$ and $\frac{{3\pi }}{2} \le \arg z \le 3\pi  - \delta  < 3\pi$. Moreover, if $K=L$ then the estimates \eqref{eq60} remain valid in the larger sector $-\frac{\pi}{2} \leq  \arg z \leq \frac{3\pi}{2}$; and \eqref{eq76} holds in the range $ - 2\pi  <  - 2\pi  + \delta  \le \arg z \le  - \frac{{3\pi }}{2}$ with the ``$-$" sign, and in the range $\frac{{3\pi }}{2} \le \arg z \le 3\pi  - \delta  < 3\pi$ with the ``$+$" sign and with an implied constant that also depends on $\delta$.
\end{theorem}

The rest of the paper is organised as follows. In Section \ref{section2}, we prove the resurgence formulas stated in Theorem \ref{thm1}. In Section \ref{section3}, we give explicit and numerically computable error bounds for the Nicholson--Debye expansions when $\kappa = 0$ using the results of Section \ref{section2}. In Section \ref{section4}, asymptotic approximations for the coefficients $B_n\left(\kappa\right)$ are given. In Section \ref{section5}, we prove the exponentially improved expansions presented in Theorem \ref{thm3}, and provide a detailed discussion of the Stokes phenomenon related to the expansion of $H_\nu ^{\left( 1 \right)} \left( z \right)$. The paper concludes with a discussion in Section \ref{section6}.

\section{Proofs of the resurgence formulas}\label{section2}

In this section we prove the resurgence formulas given in Theorem \ref{thm1}. In the first subsection, we prove the results related to the Hankel functions. In the second subsection we show how the corresponding formulas for the Bessel functions can be derived using the results for the Hankel functions.

\subsection{The Hankel functions $H_\nu^{\left(1\right)}\left(z\right)$ and $H_\nu^{\left(2\right)}\left(z\right)$} Our analysis is based on the Schl\"afli--Sommerfeld integral representation
\[
H_\nu ^{\left( 1 \right)} \left( z \right) = \frac{1}{{\pi i}}\int_{ - \infty }^{\infty  + \pi i} {e^{z\sinh t - \nu t} dt} \quad \left|\arg z\right| <\frac{\pi}{2}
\]
and the connection formula $H_\nu ^{\left( 2 \right)} \left( z \right) = \overline {H_{ \bar{\nu} }^{\left( 1 \right)} \left( {\bar z} \right)}$ \cite[eqs. 10.9.18 and 10.11.9]{NIST}. If $z = \nu + \kappa$, where $\kappa$ is a complex constant, then
\begin{equation}\label{eq1}
H_\nu ^{\left( 1 \right)} \left( {z} \right) = \frac{1}{{\pi i}}\int_{ - \infty }^{\infty  + \pi i} {e^{z\left( {\sinh t - t} \right)} e^{\kappa t} dt} \quad \left| {\arg z} \right| < \frac{\pi }{2}.
\end{equation}
The saddle points of the integrand are the roots of the equation $\cosh t = 1$. Hence, the saddle points are given by $t^{\left( k \right)} = 2\pi ik$. We denote by $\mathscr{C}^{\left( k \right)}\left(\theta\right)$ the portion of the steepest descent paths that pass through the saddle point $t^{\left( k \right)}$. Here, and subsequently, we write $\theta = \arg z$.

\begin{figure}[!t]
\def\svgwidth{0.8\columnwidth}
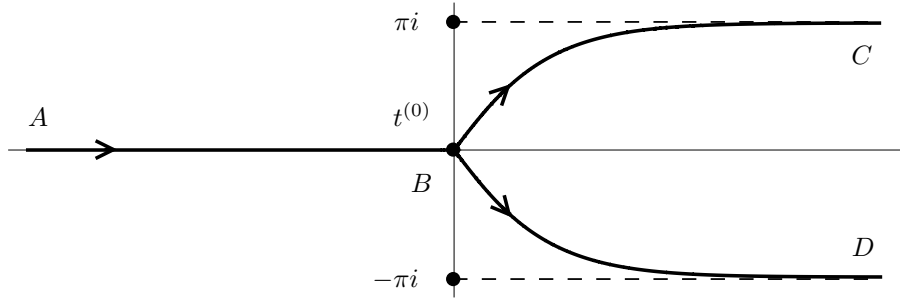
\caption{The steepest descent paths for $H_\nu ^{\left( 1 \right)} \left(z \right)$ through the saddle point $t^{\left( 0 \right)} = 0$ when $\arg z =0$.}
\label{fig1}
\end{figure}

For the integral \eqref{eq1}, we need to consider only the saddle point $t^{\left( 0 \right)} = 0$. Let $f\left(t\right) = \sinh t-t$. The steepest descent paths in this case are given by
\[
\arg \left[ {e^{i\theta } \left( {f\left( t^{\left( 0 \right)} \right) - f\left( t \right)} \right)} \right] = 0 .
\]
The path $\mathscr{C}^{\left( 0 \right)}\left(\theta\right)$ for $\theta=0$, with an appropriate orientation, is shown in Figure \ref{fig1}. We assume that $\theta = 0$ and later we shall use an analytic continuation argument to extend the results to complex values of $z$. Denote by $\mathscr{P}^{\left( 0 \right)}\left(0\right)$ and $\mathscr{L}^{\left( 0 \right)}\left(0\right)$ the curves $ABC$ and $ABD$, respectively (see Figure \ref{fig1}). By Cauchy's theorem, we can deform the path of integration in \eqref{eq1} to $\mathscr{P}^{\left( 0 \right)}\left(0\right)$. (The path $\mathscr{L}^{\left( 0 \right)}\left(0\right)$ can be taken as the integration path in the corresponding integral representation for $H_\nu^{\left(2\right)}\left(z\right)$.) Let
\[
\tau = f(t^{\left( 0 \right)} ) - f\left( t \right)=  - f\left( t \right).
\]
By the definition of $\mathscr{C}^{\left( 0 \right)} \left( 0  \right)$, for $t \in \mathscr{P}^{\left( 0 \right)} \left( 0  \right)$, $\tau$ is real and positive. Then, corresponding to each positive $\tau$, there are two values of $t$, called $t_1$ and $t_2$, such that $t_1$ is a complex number with a positive real part and $t_2$ is a negative real number. In terms of $\tau$, \eqref{eq1} becomes
\begin{align*}
H_\nu ^{\left( 1 \right)} \left( {z} \right) & = \frac{1}{{\pi i}}\int_0^{ + \infty } {e^{ - z\tau } \left( {e^{\kappa t_1 } \frac{{dt_1 }}{{d\tau }} - e^{\kappa t_2 } \frac{{dt_2 }}{{d\tau }}} \right)d\tau } 
\\ & = \frac{1}{{\pi i}}\int_0^{ + \infty } {e^{ - z\tau } \left( {\frac{{e^{\kappa t_1 \left( \tau  \right)} }}{{1 - \cosh t_1 \left( \tau  \right)}} - \frac{{e^{\kappa t_2 \left( \tau  \right)} }}{{1 - \cosh t_2 \left( \tau  \right)}}} \right)d\tau } .
\end{align*}
Following Berry and Howls, the function in the large parentheses can be represented as a contour integral using the residue theorem, to yield
\[
H_\nu ^{\left( 1 \right)} \left( {z} \right) = \frac{1}{{3\pi i}}\int_0^{ + \infty } {\tau ^{ - \frac{2}{3}} e^{ - z\tau } \frac{1}{{2\pi i}}\oint_{\Gamma ^{\left( 0 \right)} } {\left( {\frac{{e^{\frac{\pi }{3}i} }}{{1 - \tau ^{\frac{1}{3}} f^{ - \frac{1}{3}} \left( u \right)e^{\frac{\pi }{3}i} }} + \frac{1}{{1 + \tau ^{\frac{1}{3}} f^{ - \frac{1}{3}} \left( u \right)}}} \right)\frac{{e^{\kappa u} }}{{f^{\frac{1}{3}} \left( u \right)}}du} d\tau } ,
\]
where the contour $\Gamma ^{\left( 0 \right)}$ encloses $\mathscr{P}^{\left( 0 \right)} \left( 0  \right)$ in an anti-clockwise loop, and does not enclose any of the saddle points $t ^{\left( k \right)} \neq t^{\left( 0 \right)}$. The cube root is defined so that $\left(-f\left( t \right) \right)^{\frac{1}{3}}$ is positive on the portion of $\mathscr{P}^{\left( 0 \right)}\left(0\right)$ that starts at $t=-\infty$ and ends at $t= t^{\left( 0 \right)}$. The contour integral converges as long as $\left|\Re \left( \kappa  \right)\right| < \frac{1}{3}$. Next we apply the expression
\begin{equation}\label{eq48}
\frac{1}{1 - x} = \sum\limits_{n = 0}^{N-1} {x^n}  + \frac{x^N}{1 - x},\; x \neq 1,
\end{equation}
to expand the function under the contour integral in powers of $\tau ^{\frac{1}{3}} f^{ - \frac{1}{3}} \left( u \right)$. The result is
\[
H_\nu ^{\left( 1 \right)} \left( {z} \right) = -\frac{1}{{3\pi i}}\sum\limits_{n = 0}^{N - 1} {\left( {e^{ \left( {n + 1} \right)\pi i} - e^{\frac{ \left( {n + 1} \right)\pi i}{3}} } \right)\int_0^{ + \infty } {\tau ^{\frac{{n - 2}}{3}} e^{ - z\tau } \frac{1}{{2\pi i}}\oint_{\Gamma ^{\left( 0 \right)} } {\frac{{e^{\kappa u} }}{{f^{\frac{{n + 1}}{3}} \left( u \right)}}du} d\tau } }  + R_N^{\left( H \right)} \left( {z ,\kappa } \right),
\]
where
\begin{equation}\label{eq2}
R_N^{\left( H \right)} \left( {z ,\kappa } \right) = \frac{1}{{3\pi i}}\int_0^{ + \infty } {\tau ^{\frac{{N - 2}}{3}} e^{ - z\tau } \frac{1}{{2\pi i}}\oint_{\Gamma ^{\left( 0 \right)} } {\left( {\frac{{e^{\frac{\left( {N + 1} \right)\pi i}{3}} }}{{1 - \tau ^{\frac{1}{3}} f^{ - \frac{1}{3}} \left( u \right)e^{\frac{\pi }{3}i} }} - \frac{e^{\left( {N + 1} \right)\pi i}}{{1 + \tau ^{\frac{1}{3}} f^{ - \frac{1}{3}} \left( u \right)}}} \right)\frac{{e^{\kappa u} }}{{f^{\frac{{N + 1}}{3}} \left( u \right)}}du} d\tau } .
\end{equation}
By analytic continuation, this expression is valid as long as $\left|\Re \left( \kappa  \right)\right| < \frac{N + 1}{3}$. The path $\Gamma ^{\left( 0 \right)}$ in the sum can be shrunk into a small circle around $t^{\left( 0 \right)} = 0$, and we arrive at
\begin{equation}\label{eq4}
H_\nu ^{\left( 1 \right)} \left( {z} \right) =  - \frac{2}{{3\pi }}\sum\limits_{n = 0}^{N - 1} {6^{\frac{{n + 1}}{3}} B_n \left( \kappa  \right)e^{\frac{{2\left( {n + 1} \right)\pi i}}{3}} \sin \left( {\frac{{\left( {n + 1} \right)\pi }}{3}} \right)\frac{{\Gamma \left( {\frac{{n + 1}}{3}} \right)}}{{z^{\frac{{n + 1}}{3}} }}}  + R_N^{\left( H \right)} \left( {z ,\kappa } \right),
\end{equation}
with
\[
B_n \left( \kappa  \right) = \frac{{6^{ - \frac{{n + 1}}{3}} }}{{2\pi i}}\oint_{\left( {0^ +  } \right)} {\frac{{e^{\kappa u} }}{{f^{\frac{{n + 1}}{3}} \left( u \right)}}du}  = \frac{1}{{n!}}\left[ {\frac{{d^n }}{{dt^n }}\left( {e^{\kappa t} \left( {\frac{1}{6}\frac{{t^3 }}{{\sinh t - t}}} \right)^{\frac{{n + 1}}{3}} } \right)} \right]_{t = 0} .
\]

\begin{figure}[!t]
\def\svgwidth{\columnwidth}
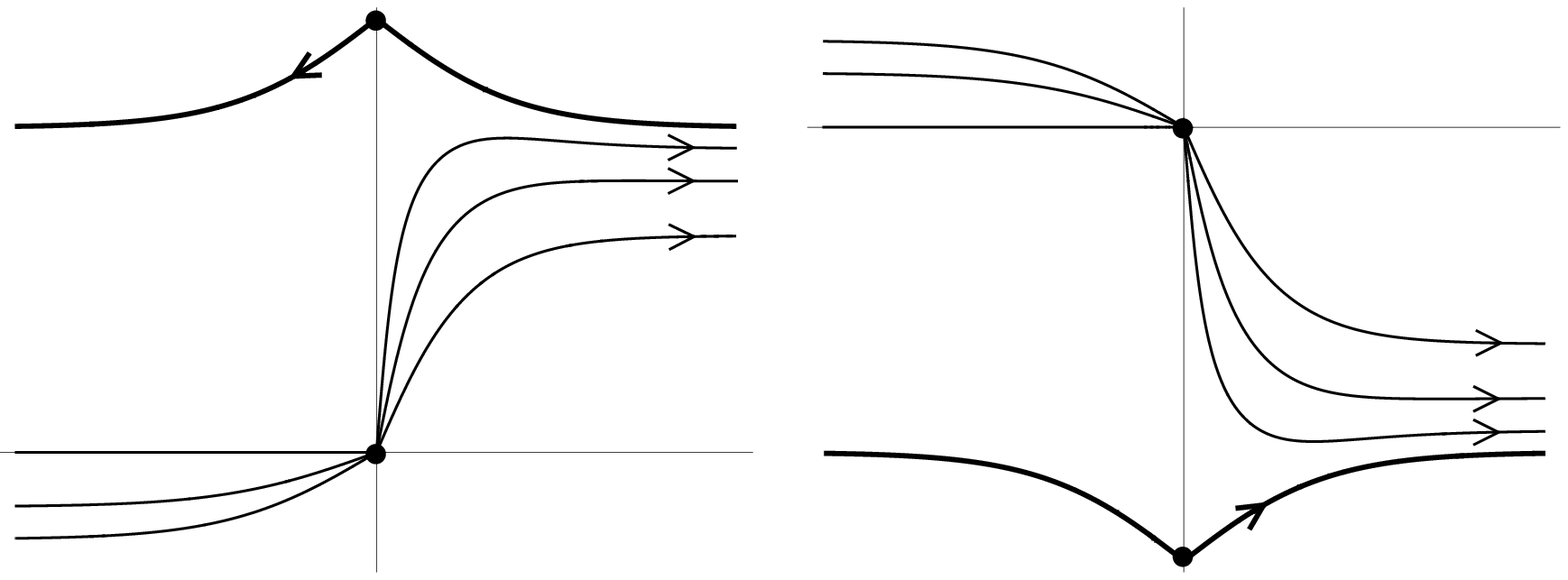
\caption{(a) The steepest descent path for $H_\nu ^{\left( 1 \right)} \left(z\right)$ through the saddle point $t^{\left( 0 \right)}$ when (i) $\arg z =0$, (ii) $\arg z =-\frac{\pi}{4}$, (iii) $\arg z =-\frac{2\pi}{5}$. The path $\mathscr{L}^{\left( 1 \right)} \left( { - \frac{\pi }{2}} \right)$ is an adjacent contour to $t^{\left( 0 \right)}$. (b) The steepest descent path $\mathscr{L}^{\left( 0 \right)} \left( \theta \right)$ through the saddle point $t^{\left( 0 \right)}$ when (i) $\arg z =0$, (ii) $\arg z = \frac{\pi}{4}$, (iii) $\arg z = \frac{2\pi}{5}$. The path $\mathscr{P}^{\left(-1\right)}\left({\frac{\pi}{2}}\right)$ is an adjacent contour to $t^{\left( 0 \right)}$.}
\label{fig2}
\end{figure}

Applying the change of variable $z\tau  = s$ in \eqref{eq2} gives
\begin{gather}\label{eq3}
\begin{split}
R_N^{\left( H \right)} \left( {z ,\kappa } \right) = \; &\frac{1}{{3\pi iz^{\frac{{N + 1}}{3}} }}\int_0^{ + \infty } s^{\frac{{N - 2}}{3}} e^{ - s} \\ & \times \frac{1}{{2\pi i}}\oint_{\Gamma ^{\left( 0 \right)} } {\left( {\frac{{e^{\frac{{\left( {N + 1} \right)\pi i}}{3}} }}{{1 - \left( {s/zf\left( u \right)} \right)^{\frac{1}{3}} e^{\frac{\pi }{3}i} }} - \frac{{e^{\left( {N + 1} \right)\pi i} }}{{1 + \left( {s/zf\left( u \right)} \right)^{\frac{1}{3}} }}} \right)\frac{{e^{\kappa u} }}{{f^{\frac{{N + 1}}{3}} \left( u \right)}}du} ds.
\end{split}
\end{gather}
As $\theta$ varies, the path $\mathscr{C}^{\left(0\right)}\left(\theta\right)$ varies smoothly, so we can define the paths $\mathscr{P}^{\left(0\right)}\left(\theta\right)$ and $\mathscr{L}^{\left(0\right)}\left(\theta\right)$. The representation \eqref{eq3} of $R_N^{\left( H \right)} \left( {z,\kappa } \right)$ and the formula \eqref{eq4} can be continued analytically if we choose $\Gamma^{\left( 0 \right)} = \Gamma^{\left( 0 \right)}\left(\theta\right)$ to be an infinite contour that surrounds the path $\mathscr{P}^{\left(0\right)}\left(\theta\right)$ in the anti-clockwise direction and that does not encircle any of the saddle points $t^{\left( k \right)} \neq t^{\left( 0 \right)}$. This continuation argument works until the path $\mathscr{P}^{\left(0\right)}\left(\theta\right)$ runs into an other saddle point. In the terminology of Berry and Howls, such saddle points are called adjacent to $t^{\left(0\right)}$. In our case, the relevant path is $\mathscr{P}^{\left(0\right)}\left(\theta\right)$, for the adjacency problem, however, we also need to consider $\mathscr{L}^{\left(0\right)}\left(\theta\right)$. When $\theta = -\frac{\pi}{2}$ or $\frac{3\pi}{2}$, the path $\mathscr{P}^{\left(0\right)}\left(\theta\right)$ connects to the saddle point $t^{\left( 1 \right)} = 2\pi i$. Similarly, when $\theta = \frac{\pi}{2}$ or $-\frac{3\pi}{2}$, the path $\mathscr{L}^{\left(0\right)}\left(\theta\right)$ connects to the saddle point $t^{\left( -1 \right)} = -2\pi i$. Therefore, the adjacent saddles are $t^{\left( \pm 1 \right)}$ (see Figure \ref{fig2}). The set
\[
\Delta ^{\left( 0 \right)}  = \left\{ {u \in \mathscr{P}^{\left( 0 \right)} \left( \theta  \right) : - \frac{\pi }{2} < \theta  < \frac{{3\pi }}{2}} \right\} \cup \left\{ {u \in \mathscr{L}^{\left( 0 \right)} \left( \theta  \right) : - \frac{3\pi }{2} < \theta  < \frac{{\pi }}{2}} \right\}
\]
forms a domain in the complex plane whose boundaries are the steepest descent paths $\mathscr{L}^{\left( 1 \right)} \left( { - \frac{\pi }{2}} \right)$ and $\mathscr{P}^{\left( -1 \right)} \left( { \frac{\pi }{2}} \right)$, the adjacent contours to $t^{\left( 0 \right)}$ (they are defined analogously to $\mathscr{P}^{\left(0\right)}$ and $\mathscr{L}^{\left(0\right)}$). For $N\geq 0$, $\left|\Re \left( \kappa  \right)\right| < \frac{N + 1}{3}$ and fixed $z$, the function under the contour integral in \eqref{eq3} is an analytic function of $u$ in the domain $\Delta ^{\left( 0 \right)}$ excluding $\mathscr{P}^{\left(0\right)}\left(\theta\right)$ and $\mathscr{L}^{\left(0\right)}\left(\theta\right)$, therefore we can deform $\Gamma^{\left( 0 \right)}$ over the adjacent contours. We thus find that for $-\frac{\pi}{2} < \theta < \frac{3\pi}{2}$ and $N \geq 0$, \eqref{eq3} can be written
\begin{gather}\label{eq5}
\begin{split}
R_N^{\left( H \right)} \left( {z ,\kappa } \right) =\; & \frac{1}{{3\pi iz^{\frac{{N + 1}}{3}} }}\int_0^{ + \infty } s^{\frac{{N - 2}}{3}} e^{ - s} \\ & \times \frac{1}{{2\pi i}}\int_{\mathscr{L}^{\left( 1 \right)} \left( { - \frac{\pi }{2}} \right)} {\left( {\frac{{e^{\frac{{\left( {N + 1} \right)\pi i}}{3}} }}{{1 - \left( {s/zf\left( u \right)} \right)^{\frac{1}{3}} e^{\frac{\pi }{3}i} }} - \frac{{e^{\left( {N + 1} \right)\pi i} }}{{1 + \left( {s/zf\left( u \right)} \right)^{\frac{1}{3}} }}} \right)\frac{{e^{\kappa u} }}{{f^{\frac{{N + 1}}{3}} \left( u \right)}}du} ds
\\ & + \frac{1}{{3\pi iz^{\frac{{N + 1}}{3}} }}\int_0^{ + \infty } s^{\frac{{N - 2}}{3}} e^{ - s} \\ & \times \frac{1}{{2\pi i}}\int_{\mathscr{P}^{\left( { - 1} \right)} \left( {\frac{\pi }{2}} \right)} {\left( {\frac{{e^{\frac{{\left( {N + 1} \right)\pi i}}{3}} }}{{1 - \left( {s/zf\left( u \right)} \right)^{\frac{1}{3}} e^{\frac{\pi }{3}i} }} - \frac{{e^{\left( {N + 1} \right)\pi i} }}{{1 + \left( {s/zf\left( u \right)} \right)^{\frac{1}{3}} }}} \right)\frac{{e^{\kappa u} }}{{f^{\frac{{N + 1}}{3}} \left( u \right)}}du} ds,
\end{split}
\end{gather}
provided that $\left|\Re \left( \kappa  \right)\right| < \frac{N + 1}{3}$. Now we perform the changes of variable
\[
s = t\frac{{\left| {f\left( {2\pi i} \right) - f\left( 0 \right)} \right|}}{{f\left( {2\pi i} \right) - f\left( 0 \right)}}\left( {f\left( u \right) - f\left( 0 \right)} \right) = itf\left( u \right)
\]
in the first, and
\[
s = t\frac{{\left| {f\left( { - 2\pi i} \right) - f\left( 0 \right)} \right|}}{{f\left( { - 2\pi i} \right) - f\left( 0 \right)}}\left( {f\left( u \right) - f\left( 0 \right)} \right) =  - itf\left( u \right)
\]
in the second double integral. The quantities $f\left( { \pm 2\pi i} \right) - f\left( 0 \right)=  \mp 2\pi i$ were essentially called the ``singulants" by Dingle \cite[p. 147]{Dingle}. When using these changes of variable, we should take $i^{\frac{1}{3}}  = -i$ in the first, and $\left( { - i} \right)^{\frac{1}{3}}  = i$ in the second double integral. With these changes of variable, the representation \eqref{eq5} for $R_N^{\left( H \right)} \left( {z ,\kappa } \right)$ becomes
\begin{gather}\label{eq6}
\begin{split}
R_N^{\left( H \right)} \left( {z ,\kappa } \right) = \; & \frac{\left( { - i} \right)^N}{{3\pi z^{\frac{{N + 1}}{3}} }}\int_0^{ + \infty }  t^{\frac{{N - 2}}{3}} \left( {\frac{{e^{\frac{{\left( {N + 1} \right)\pi i}}{3}} }}{{1 + i\left( {t/z} \right)^{\frac{1}{3}} e^{\frac{\pi }{3}i} }} - \frac{{e^{\left( {N + 1} \right)\pi i} }}{{1 - i\left( {t/z} \right)^{\frac{1}{3}} }}} \right) \\ & \times \frac{{ - 1}}{{2\pi i}}\int_{\mathscr{L}^{\left( 1 \right)} \left( { - \frac{\pi }{2}} \right)} {e^{ - itf\left( u \right)} e^{\kappa u} du} dt 
\\ & + \frac{i^N}{{3\pi z^{\frac{{N + 1}}{3}} }}\int_0^{ + \infty }  t^{\frac{{N - 2}}{3}} \left( {\frac{{e^{\frac{{\left( {N + 1} \right)\pi i}}{3}} }}{{1 - i\left( {t/z} \right)^{\frac{1}{3}} e^{\frac{\pi }{3}i} }} - \frac{{e^{\left( {N + 1} \right)\pi i} }}{{1 + i\left( {t/z} \right)^{\frac{1}{3}} }}} \right) \\ & \times \frac{1}{{2\pi i}}\int_{\mathscr{P}^{\left( { - 1} \right)} \left( {\frac{\pi }{2}} \right)} {e^{itf\left( u \right)} e^{\kappa u} du} dt ,
\end{split}
\end{gather}
for $-\frac{\pi}{2} < \theta < \frac{3\pi}{2}$, $N \geq 0$ and $\left|\Re \left( \kappa  \right)\right| < \frac{N + 1}{3}$. Finally, the contour integrals can themselves be represented in terms of the Hankel functions since
\begin{multline*}
\frac{{ - 1}}{{2\pi i}}\int_{\mathscr{L}^{\left( 1 \right)} \left( { - \frac{\pi }{2}} \right)} {e^{ - itf\left( u \right)} e^{\kappa u} du}  = \frac{1}{{2\pi i}}\int_{\mathscr{L}^{\left( 0 \right)} \left( { - \frac{\pi }{2}} \right)} {e^{ - itf\left( {u + 2\pi i} \right)} e^{\kappa \left( {u + 2\pi i} \right)} du} 
\\ = \frac{{e^{2\pi i\left( {it + \kappa } \right)} }}{{2\pi i}}\int_{\mathscr{L}^{\left( 0 \right)} \left( { - \frac{\pi }{2}} \right)} {e^{ - itf\left( u \right)} e^{\kappa u} du}  =  - \frac{{e^{2\pi i\left( {it + \kappa } \right)} }}{2}H_{ - it - \kappa }^{\left( 2 \right)} \left( { - it} \right) = \frac{{e^{2\pi i\left( {it + \kappa } \right)} }}{2}H_{it + \kappa }^{\left( 1 \right)} \left( {it} \right),
\end{multline*}
and
\begin{multline*}
\frac{1}{{2\pi i}}\int_{\mathscr{P}^{\left( { - 1} \right)} \left( {\frac{\pi }{2}} \right)} {e^{itf\left( u \right)} e^{\kappa u} du}  = \frac{1}{{2\pi i}}\int_{\mathscr{P}^{\left( 0 \right)} \left( {\frac{\pi }{2}} \right)} {e^{itf\left( {u - 2\pi i} \right)} e^{\kappa \left( {u - 2\pi i} \right)} du} \\  = \frac{{e^{2\pi i\left( {it - \kappa } \right)} }}{{2\pi i}}\int_{\mathscr{P}^{\left( 0 \right)} \left( {\frac{\pi }{2}} \right)} {e^{itf\left( u \right)} e^{\kappa u} du}  = \frac{{e^{2\pi i\left( {it - \kappa } \right)} }}{2}H_{it - \kappa }^{\left( 1 \right)} \left( {it} \right).
\end{multline*}
Substituting these into \eqref{eq6} gives \eqref{eq9}. The expansion \eqref{eq7} follows from the connection formula $H_\nu ^{\left( 2 \right)} \left( z \right) = \overline {H_{ \bar{\nu} }^{\left( 1 \right)} \left( {\bar{z } } \right)}$. Finally, we prove the second representation in \eqref{eq8}. We note that
\begin{multline*}
\Gamma \left( {\frac{{n + 1}}{3}} \right)6^{\frac{{n + 1}}{3}} B_n \left( \kappa  \right) = \int_0^{ + \infty } {s^{\frac{{n - 2}}{3}} e^{ - s} \frac{1}{{2\pi i}}\oint_{\left( {0^ +  } \right)} {\frac{{e^{\kappa u} }}{{f^{\frac{{n + 1}}{3}} \left( u \right)}}du} ds} 
\\ = \int_0^{ + \infty } {s^{\frac{{n - 2}}{3}} e^{ - s} \frac{1}{{2\pi i}}\int_{\mathscr{L}^{\left( 1 \right)} \left( { - \frac{\pi }{2}} \right)} {\frac{{e^{\kappa u} }}{{f^{\frac{{n + 1}}{3}} \left( u \right)}}du} ds}  + \int_0^{ + \infty } {s^{\frac{{n - 2}}{3}} e^{ - s} \frac{1}{{2\pi i}}\int_{\mathscr{P}^{\left( { - 1} \right)} \left( {\frac{\pi }{2}} \right)} {\frac{{e^{\kappa u} }}{{f^{\frac{{n + 1}}{3}} \left( u \right)}}du} ds} ,
\end{multline*}
as long as $\left|\Re \left( \kappa  \right)\right| < \frac{n + 1}{3}$. For these integrals we can apply the same changes of variable as above, then using the same argument yields the second representation in \eqref{eq8}.

\subsection{The Bessel functions $J_\nu\left(z\right)$ and $Y_\nu\left(z\right)$} To prove the resurgence representations for the Bessel functions, we apply the well-known connection formulas
\begin{equation}\label{eq11}
J_\nu  \left( z \right) = \frac{1}{2}\left( {H_\nu ^{\left( 1 \right)} \left( z \right) + H_\nu ^{\left( 2 \right)} \left( z \right)} \right),
\end{equation}
\begin{equation}\label{eq21}
Y_\nu  \left( z \right) = \frac{1}{{2i}}\left( {H_\nu ^{\left( 1 \right)} \left( z \right) - H_\nu ^{\left( 2 \right)} \left( z \right)} \right).
\end{equation}
To show \eqref{eq13} and \eqref{eq15}, we substitute \eqref{eq12} and \eqref{eq7} into \eqref{eq11} and note that
\[
\frac{{e^{\frac{{2\left( {n + 1} \right)\pi i}}{3}}  + e^{ - \frac{{2\left( {n + 1} \right)\pi i}}{3}} }}{2}\sin \left( {\frac{{\left( {n + 1} \right)\pi }}{3}} \right) = \cos \left( {\frac{{2\left( {n + 1} \right)\pi }}{3}} \right)\sin \left( {\frac{{\left( {n + 1} \right)\pi }}{3}} \right) =  - \frac{1}{2}\sin \left( {\frac{{\left( {n + 1} \right)\pi }}{3}} \right).
\]
To prove \eqref{eq14} and \eqref{eq16}, we substitute \eqref{eq12} and \eqref{eq7} into \eqref{eq21} and use the fact that
\begin{align*}
& \frac{{e^{\frac{{2\left( {n + 1} \right)\pi i}}{3}}  - e^{ - \frac{{2\left( {n + 1} \right)\pi i}}{3}} }}{{2i}}\sin \left( {\frac{{\left( {n + 1} \right)\pi }}{3}} \right) = \sin \left( {\frac{{2\left( {n + 1} \right)\pi }}{3}} \right)\sin \left( {\frac{{\left( {n + 1} \right)\pi }}{3}} \right) \\ & =  2\cos \left( {\frac{{\left( {n + 1} \right)\pi }}{3}} \right)\sin ^2 \left( {\frac{{\left( {n + 1} \right)\pi }}{3}} \right) = \left( { - 1} \right)^n \sin ^2 \left( {\frac{{\left( {n + 1} \right)\pi }}{3}} \right).
\end{align*}

\section{Error bounds}\label{section3} In this section we derive error bounds for the asymptotic expansions \eqref{eq17}--\eqref{eq20} and the analogous asymptotic series for the $z$-derivatives, in the special case $\kappa=0$, i.e., when $z =\nu$. Since $B_n\left(0\right) =0$ for odd $n$, we find it convenient to introduce the following simplified notations. For any $N\geq 0$, we write
\begin{equation}\label{eq24}
H_\nu ^{\left( 1 \right)} \left( \nu  \right) =  - \frac{2}{{3\pi }}\sum\limits_{n = 0}^{N - 1} {6^{\frac{{2n + 1}}{3}} B_{2n} \left( 0 \right)e^{\frac{{2\left( {2n + 1} \right)\pi i}}{3}} \sin \left( {\frac{{\left( {2n + 1} \right)\pi }}{3}} \right)\frac{{\Gamma \left( {\frac{{2n + 1}}{3}} \right)}}{{\nu ^{\frac{{2n + 1}}{3}} }}}  + R_N^{\left( H \right)} \left( \nu  \right),
\end{equation}
\[
H_\nu ^{\left( 2 \right)} \left( \nu  \right) =  - \frac{2}{{3\pi }}\sum\limits_{n = 0}^{N - 1} {6^{\frac{{2n + 1}}{3}} B_{2n} \left( 0 \right)e^{ - \frac{{2\left( {2n + 1} \right)\pi i}}{3}} \sin \left( {\frac{{\left( {2n + 1} \right)\pi }}{3}} \right)\frac{{\Gamma \left( {\frac{{2n + 1}}{3}} \right)}}{{\nu ^{\frac{{2n + 1}}{3}} }}}  - R_N^{\left( H \right)} \left( {\nu e^{\pi i} } \right),
\]
\begin{equation}\label{eq41}
J_\nu  \left( \nu  \right) = \frac{1}{{3\pi }}\sum\limits_{n = 0}^{N - 1} {6^{\frac{{2n + 1}}{3}} B_{2n} \left( 0 \right)\sin \left( {\frac{{\left( {2n + 1} \right)\pi }}{3}} \right)\frac{{\Gamma \left( {\frac{{2n + 1}}{3}} \right)}}{{\nu ^{\frac{{2n + 1}}{3}} }}}  + R_N^{\left( J \right)} \left( \nu  \right)
\end{equation}
and
\begin{equation}\label{eq42}
Y_\nu  \left( \nu  \right) =  - \frac{2}{{3\pi }}\sum\limits_{n = 0}^{N - 1} {6^{\frac{{2n + 1}}{3}} B_{2n} \left( 0 \right)\sin ^2 \left( {\frac{{\left( {2n + 1} \right)\pi }}{3}} \right)\frac{{\Gamma \left( {\frac{{2n + 1}}{3}} \right)}}{{\nu ^{\frac{{2n + 1}}{3}} }}}  + R_N^{\left( Y \right)} \left( \nu  \right),
\end{equation}
with the notations $R_N^{\left( H \right)} \left( \nu  \right) = R_{2N}^{\left( {H} \right)} \left( {\nu ,0} \right)$, $R_N^{\left( J \right)} \left( \nu  \right) = R_{2N}^{\left( J \right)} \left( {\nu ,0} \right)$ and $R_N^{\left( Y \right)} \left( \nu  \right) = R_{2N}^{\left( Y \right)} \left( {\nu ,0} \right)$.

Similarly, $D_n\left(0\right) =0$ for even $n$, whence we write for $N\geq 1$,
\begin{equation}\label{eq33}
H_\nu ^{\left( 1 \right)\prime  } \left( \nu  \right) =  - \frac{2}{{3\pi }}\sum\limits_{n = 0}^{N - 1} {6^{\frac{{2n + 2}}{3}} D_{2n + 1} \left( 0 \right)e^{\frac{{2\left( {2n + 2} \right)\pi i}}{3}} \sin \left( {\frac{{\left( {2n + 2} \right)\pi }}{3}} \right)\frac{{\Gamma \left( {\frac{{2n + 2}}{3}} \right)}}{{\nu ^{\frac{{2n + 2}}{3}} }}}  + R_{N}^{\left( {H'} \right)} \left( {\nu} \right),
\end{equation}
\[
H_\nu ^{\left( 2 \right)\prime  } \left( \nu  \right) =  - \frac{2}{{3\pi }}\sum\limits_{n = 0}^{N - 1} {6^{\frac{{2n + 2}}{3}} D_{2n + 1} \left( 0 \right)e^{ - \frac{{2\left( {2n + 2} \right)\pi i}}{3}} \sin \left( {\frac{{\left( {2n + 2} \right)\pi }}{3}} \right)\frac{{\Gamma \left( {\frac{{2n + 2}}{3}} \right)}}{{\nu ^{\frac{{2n + 2}}{3}} }}}  + R_N^{\left( {H'} \right)} \left( {\nu e^{\pi i} } \right),
\]
\begin{equation}\label{eq64}
J'_\nu  \left( \nu  \right) = \frac{1}{{3\pi }}\sum\limits_{n = 0}^{N - 1} {6^{\frac{{2n + 2}}{3}} D_{2n + 1} \left( 0 \right)\sin \left( {\frac{{\left( {2n + 2} \right)\pi }}{3}} \right)\frac{{\Gamma \left( {\frac{{2n + 2}}{3}} \right)}}{{\nu ^{\frac{{2n + 2}}{3}} }}}  + R_N^{\left( {J'} \right)} \left( \nu  \right)
\end{equation}
and
\begin{equation}\label{eq68}
Y'_\nu  \left( \nu  \right) = \frac{2}{{3\pi }}\sum\limits_{n = 0}^{N - 1} {6^{\frac{{2n + 2}}{3}} D_{2n + 1} \left( 0 \right)\sin ^2 \left( {\frac{{\left( {2n + 2} \right)\pi }}{3}} \right)\frac{{\Gamma \left( {\frac{{2n + 2}}{3}} \right)}}{{\nu ^{\frac{{2n + 2}}{3}} }}}  + R_N^{\left( {Y'} \right)} \left( \nu  \right),
\end{equation}
with the notations $R_N^{\left( {H'} \right)} \left( \nu  \right) = R_{2N + 1}^{\left( {H'} \right)} \left( {\nu ,0} \right)$, $R_N^{\left( {J'} \right)} \left( \nu  \right) = R_{2N + 1}^{\left( {J'} \right)} \left( {\nu ,0} \right)$ and $R_N^{\left( {Y'} \right)} \left( \nu  \right) = R_{2N + 1}^{\left( {Y'} \right)} \left( {\nu ,0} \right)$.

Some of the results for $R_N^{\left( H \right)} \left( \nu  \right)$, $R_N^{\left( J \right)} \left( \nu  \right)$ and $R_N^{\left( Y \right)} \left( \nu  \right)$ are improvements over the recently obtained error bounds \cite{Nemes} due to the present author. As far as we know, numerically computable bounds for these remainders were first obtained by Gatteschi \cite{Gatteschi}. The advantage of our bounds is that they are directly related to the absolute value of the first non-vanishing omitted term(s) of the asymptotic series.

\subsection{Error bounds for the asymptotic series of $H_\nu ^{\left( 1 \right)} \left( \nu  \right)$ and $H_\nu ^{\left( 2 \right)} \left( \nu  \right)$} Since the remainder in the asymptotic series of $H_\nu ^{\left( 2 \right)} \left( \nu  \right)$ is $- R_N^{\left( H \right)} \left( {\nu e^{\pi i} } \right)$, it is enough to estimate $R_N^{\left( H \right)} \left( \nu  \right)$. We shall estimate $R_N^{\left( H \right)} \left( \nu  \right)$ in the large sector $-\pi <\arg \nu <2\pi$. The representation \eqref{eq9} is valid only for $-\frac{\pi}{2} <\arg \nu <\frac{3\pi}{2}$, for the larger sector we define $R_N^{\left( H \right)} \left( \nu  \right)$ via \eqref{eq24}. Algebraic manipulation of \eqref{eq9} shows that
\begin{align}
R_N^{\left( H \right)} \left( \nu  \right) = \; & \frac{{\left( { - 1} \right)^{N + 1} }}{{\sqrt 3 \pi \nu ^{\frac{{2N + 1}}{3}} }}\int_0^{ + \infty } { \frac{t^{\frac{{2N - 2}}{3}} e^{ - 2\pi t} e^{\frac{{2\left( {2N + 1} \right)\pi i}}{3}}}{{\left( {1 + \left( {t/\nu } \right)^{\frac{2}{3}} e^{\frac{{2\pi i}}{3}} } \right)\left( {1 + \left( {t/\nu } \right)^{\frac{2}{3}} } \right)}}iH_{it}^{\left( 1 \right)} \left( {it} \right)dt}  \label{eq29}
\\ & + \frac{{\left( { - 1} \right)^{N + 1} }}{{\sqrt 3 \pi \nu ^{\frac{{2N + 3}}{3}} }}\int_0^{ + \infty } { \frac{{t^{\frac{{2N}}{3}} e^{ - 2\pi t} e^{\frac{{2\left( {2N + 1} \right)\pi i}}{3}} e^{\frac{\pi}{3}i} }}{{\left( {1 + \left( {t/\nu } \right)^{\frac{2}{3}} e^{\frac{{2\pi i}}{3}} } \right)\left( {1 + \left( {t/\nu } \right)^{\frac{2}{3}} } \right)}}iH_{it}^{\left( 1 \right)} \left( {it} \right)dt} ,
\nonumber \\ R_N^{\left( H \right)} \left( \nu  \right) = \; & \frac{{\left( { - 1} \right)^N }}{{\sqrt 3 \pi \nu ^{\frac{{2N + 3}}{3}} }}\int_0^{ + \infty } { \frac{{t^{\frac{{2N}}{3}} e^{ - 2\pi t} e^{\frac{{2\left( {2N + 1} \right)\pi i}}{3}} e^{\frac{\pi}{3}i} }}{{\left( {1 + \left( {t/\nu } \right)^{\frac{2}{3}} e^{\frac{{2\pi i}}{3}} } \right)\left( {1 + \left( {t/\nu } \right)^{\frac{2}{3}} } \right)}}i H_{it}^{\left( 1 \right)} \left( {it} \right)dt},
\label{eq27} \\ R_N^{\left( H \right)} \left( \nu  \right) = \; & \frac{{\left( { - 1} \right)^N }}{{\sqrt 3 \pi \nu ^{\frac{{2N + 1}}{3}} }}\int_0^{ + \infty } { \frac{t^{\frac{{2N - 2}}{3}} e^{ - 2\pi t} e^{\frac{{2\left( {2N + 1} \right)\pi i}}{3}}}{{\left( {1 + \left( {t/\nu } \right)^{\frac{2}{3}} e^{\frac{{2\pi i}}{3}} } \right)\left( {1 + \left( {t/\nu } \right)^{\frac{2}{3}} } \right)}}i H_{it}^{\left( 1 \right)} \left( {it} \right)dt} \label{eq28}
\end{align}
according to whether $N\equiv 0 \mod 3$, $N\equiv 1 \mod 3$ or $N\equiv 2 \mod 3$, respectively. Bounds for the sector $-\frac{\pi}{2}<\arg \nu <\frac{3\pi}{2}$ were derived in the earlier paper \cite{Nemes}. With the present notations, these bounds are
\begin{align*}
\left| {R_N^{\left( H \right)} \left( \nu  \right)} \right| \le \; & \left(\frac{2}{{3\pi }}6^{\frac{{2N + 1}}{3}} \left| {B_{2N} \left( 0 \right)} \right|\frac{{\sqrt 3 }}{2}\frac{{\Gamma \left( {\frac{{2N + 1}}{3}} \right)}}{{\left| \nu  \right|^{\frac{{2N + 1}}{3}} }} \right. \\ & \left.+ \frac{2}{{3\pi }}6^{\frac{{2N + 3}}{3}} \left| {B_{2N + 2} \left( 0 \right)} \right|\frac{{\sqrt 3 }}{2}\frac{{\Gamma \left( {\frac{{2N + 3}}{3}} \right)}}{{\left| \nu  \right|^{\frac{{2N + 3}}{3}} }}\right) \begin{cases} \left|\sec \theta \right| & \; \text{ if } \; { - \frac{\pi }{2} < \theta  < 0 \; \text{ or } \; \pi  < \theta  < \frac{{3\pi }}{2}}, \\ 1 & \; \text{ if } \; {0 \le \theta  \le \pi }, \end{cases}
\end{align*}
if $N\equiv 0 \mod 3$;
\begin{equation}\label{eq31}
\left| {R_N^{\left( H \right)} \left( \nu  \right)} \right| \le \frac{2}{{3\pi }}6^{\frac{{2N + 3}}{3}} \left| {B_{2N + 2} \left( 0 \right)} \right|\frac{{\sqrt 3 }}{2}\frac{{\Gamma \left( {\frac{{2N + 3}}{3}} \right)}}{{\left| \nu  \right|^{\frac{{2N + 3}}{3}} }} \begin{cases} \left|\sec \theta \right| & \; \text{ if } \; { - \frac{\pi }{2} < \theta  < 0 \; \text{ or } \; \pi  < \theta  < \frac{{3\pi }}{2}}, \\ 1 & \; \text{ if } \; {0 \le \theta  \le \pi }, \end{cases}
\end{equation}
when $N\equiv 1 \mod 3$; and
\[
\left| {R_N^{\left( H \right)} \left( \nu  \right)} \right| \le \frac{2}{{3\pi }}6^{\frac{{2N + 1}}{3}} \left| {B_{2N} \left( 0 \right)} \right|\frac{{\sqrt 3 }}{2}\frac{{\Gamma \left( {\frac{{2N + 1}}{3}} \right)}}{{\left| \nu  \right|^{\frac{{2N + 1}}{3}} }} \begin{cases} \left|\sec \theta \right| & \; \text{ if } \; { - \frac{\pi }{2} < \theta  < 0 \; \text{ or } \; \pi  < \theta  < \frac{{3\pi }}{2}}, \\ 1 & \; \text{ if } \; {0 \le \theta  \le \pi }, \end{cases}
\]
if $N\equiv 2 \mod 3$. Here, as before, $\theta = \arg \nu$. Now we focus on the sectors $-\pi <\arg \nu<0$ and $\pi <\arg \nu <2\pi$. We consider the case $N\equiv 1 \mod 3$. We introduce the integral representation of $e^{ - 2\pi t} i H_{it}^{\left( 1 \right)} \left( {it} \right)$ and perform the change of variable $is/f\left( u \right) = t$ in \eqref{eq27}, to obtain
\begin{gather}\label{eq22}
\begin{split}
& R_N^{\left( H \right)} \left( \nu  \right) = \frac{{\left( { - 1} \right)^N }}{{\sqrt 3 \pi \nu ^{\frac{{2N + 3}}{3}} }}\int_0^{ + \infty } {\frac{{t^{\frac{{2N}}{3}} e^{\frac{{2\left( {2N + 1} \right)\pi i}}{3}}  e^{\frac{\pi }{3}i} }}{{\left( {1 + \left( {t/\nu } \right)^{\frac{2}{3}} e^{\frac{{2\pi i}}{3}} } \right)\left( {1 + \left( {t/\nu } \right)^{\frac{2}{3}} } \right)}}\frac{1}{\pi }\int_{\mathscr{P}^{\left( { - 1} \right)} \left( {\frac{\pi }{2}} \right)} {e^{itf\left( u \right)} du} dt} 
\\ & = \frac{{\left( { - 1} \right)^N }}{{\sqrt 3 \pi \nu ^{\frac{{2N + 3}}{3}} }}\int_0^{ + \infty } {s^{\frac{{2N}}{3}} e^{ - s} e^{\frac{{2\left( {2N + 1} \right)\pi i}}{3}} e^{\frac{\pi }{3}i} \frac{1}{\pi }\int_{\mathscr{P}^{\left( { - 1} \right)} \left( {\frac{\pi }{2}} \right)} {\frac{{\left( {i/f\left( u \right)} \right)^{\frac{{2N + 3}}{3}} }}{{\left( {1 + \left( {is/\nu f\left( u \right) } \right)^{\frac{2}{3}} e^{\frac{{2\pi i}}{3}} } \right)\left( {1 + \left( {is/\nu f\left( u \right) } \right)^{\frac{2}{3}} } \right)}}du} ds} .
\end{split}
\end{gather}
Note that along the path $\mathscr{P}^{\left( { - 1} \right)} \left( {\frac{\pi }{2}} \right)$, we have
\[
0 \le i\left( {f\left( { - 2\pi i} \right) - f\left( u \right)} \right) =  - 2\pi  - if\left( u \right) <  - if\left( u \right),
\]
thus, $s$ is indeed positive. Set $u=x+iy$. Along the path $\mathscr{P}^{\left( { - 1} \right)} \left( {\frac{\pi }{2}} \right)$, we also have
\[
0 <  - if\left( u \right) =  - i\left( {\sinh \left( {x + iy} \right) - \left( {x + iy} \right)} \right) = \Re \left( { - i\left( {\sinh \left( {x + iy} \right) - \left( {x + iy} \right)} \right)} \right) = \cosh x\sin y - y,
\]
and therefore, $- if\left( {x + iy} \right) =  - if\left( { - x + iy} \right)$. Denote by $\mathscr{P}_\pm ^{\left( { - 1} \right)}\left( {\frac{\pi }{2}} \right)$ the portions of the steepest path $\mathscr{P}^{\left( { - 1} \right)}\left( {\frac{\pi }{2}} \right)$ that lie in the left- and in the right-half plane, respectively. Using the previous observation and the fact that $\mathscr{P}^{\left( { - 1} \right)}\left( {\frac{\pi }{2}} \right)$ is symmetric with respect to the imaginary axis (see Figure \ref{fig2}), we find
\begin{align*}
& \int_{\mathscr{P}^{\left( { - 1} \right)} \left( {\frac{\pi }{2}} \right)} {\frac{{\left( {i/f\left( u \right)} \right)^{\frac{{2N + 3}}{3}} }}{{\left( {1 + \left( {is/f\left( u \right)\nu } \right)^{\frac{2}{3}} e^{\frac{{2\pi i}}{3}} } \right)\left( {1 + \left( {is/f\left( u \right)\nu } \right)^{\frac{2}{3}} } \right)}}du} 
\\ = \; & \int_{\mathscr{P}_ - ^{\left( { - 1} \right)} \left( {\frac{\pi }{2}} \right)} {\frac{{\left( {i/f\left( {x + iy} \right)} \right)^{\frac{{2N + 3}}{3}} }}{{\left( {1 + \left( {is/f\left( {x + iy} \right)\nu } \right)^{\frac{2}{3}} e^{\frac{{2\pi i}}{3}} } \right)\left( {1 + \left( {is/f\left( {x + iy} \right)\nu } \right)^{\frac{2}{3}} } \right)}}d\left( {x + iy} \right)} 
\\ & + \int_{\mathscr{P}_ + ^{\left( { - 1} \right)} \left( {\frac{\pi }{2}} \right)} {\frac{{\left( {i/f\left( x+iy \right)} \right)^{\frac{{2N + 3}}{3}} }}{{\left( {1 + \left( {is/f\left( {x + iy} \right)\nu } \right)^{\frac{2}{3}} e^{\frac{{2\pi i}}{3}} } \right)\left( {1 + \left( {is/f\left( {x + iy} \right)\nu } \right)^{\frac{2}{3}} } \right)}}d\left( {x + iy} \right)}
\\ = \; & \int_{\mathscr{P}_ + ^{\left( { - 1} \right)} \left( {\frac{\pi }{2}} \right)} {\frac{{\left( {i/f\left( { - x + iy} \right)} \right)^{\frac{{2N + 3}}{3}} }}{{\left( {1 + \left( {is/f\left( { - x + iy} \right)\nu } \right)^{\frac{2}{3}} e^{\frac{{2\pi i}}{3}} } \right)\left( {1 + \left( {is/f\left( { - x + iy} \right)\nu } \right)^{\frac{2}{3}} } \right)}}d\left( {x - iy} \right)} 
\\ & + \int_{\mathscr{P}_ + ^{\left( { - 1} \right)} \left( {\frac{\pi }{2}} \right)} {\frac{{\left( {i/f\left( {x + iy} \right)} \right)^{\frac{{2N + 3}}{3}} }}{{\left( {1 + \left( {is/f\left( {x + iy} \right)\nu } \right)^{\frac{2}{3}} e^{\frac{{2\pi i}}{3}} } \right)\left( {1 + \left( {is/f\left( {x + iy} \right)\nu } \right)^{\frac{2}{3}} } \right)}}d\left( {x + iy} \right)} 
\\ = \; & 2\int_{\mathscr{P}_ + ^{\left( { - 1} \right)} \left( {\frac{\pi }{2}} \right)} {\frac{{\left( {i/f\left( u \right)} \right)^{\frac{{2N + 3}}{3}} }}{{\left( {1 + \left( {is/f\left( u \right)\nu } \right)^{\frac{2}{3}} e^{\frac{{2\pi i}}{3}} } \right)\left( {1 + \left( {is/f\left( u \right)\nu } \right)^{\frac{2}{3}} } \right)}}dx} .
\end{align*}
Let $0 < \varphi  < \frac{\pi }{2}$ be an acute angle that may depend on $N$ and suppose that $\pi + \varphi < \theta < \frac{3\pi }{2} + \varphi$. We rotate the path of integration in \eqref{eq22} by $\varphi$, to obtain the analytic continuation of the representation \eqref{eq9} to the sector $\pi + \varphi < \theta < \frac{3\pi }{2} + \varphi$:
\begin{gather}\label{eq23}
\begin{split}
R_N^{\left( H \right)} \left( \nu  \right) = \; & \frac{{\left( { - 1} \right)^N }}{{\sqrt 3 \pi \nu ^{\frac{{2N + 3}}{3}} }}\int_0^{ + \infty e^{i\varphi } } s^{\frac{{2N}}{3}} e^{ - s} e^{\frac{{2\left( {2N + 1} \right)\pi i}}{3}} e^{\frac{\pi }{3}i} \\ & \times \frac{2}{\pi }\int_{\mathscr{P}_ + ^{\left( { - 1} \right)} \left( {\frac{\pi }{2}} \right)} {\frac{{\left( {i/f\left( u \right)} \right)^{\frac{{2N + 3}}{3}} }}{{\left( {1 + \left( {is/f\left( u \right)\nu } \right)^{\frac{2}{3}} e^{\frac{{2\pi i}}{3}} } \right)\left( {1 + \left( {is/f\left( u \right)\nu } \right)^{\frac{2}{3}} } \right)}}dx} ds .
\end{split}
\end{gather}
Note that $x$ is positive and monotonically increasing along $\mathscr{P}_ + ^{\left( { - 1} \right)} \left( {\frac{\pi }{2}} \right)$. It was shown in the previous paper \cite[Appendix B]{Nemes}, that
\begin{equation}\label{eq36}
\frac{1}{{\left| {1 + re^{ - \frac{2}{3}i\vartheta } e^{\frac{{2\pi i}}{3}} } \right|\left| {1 + re^{ - \frac{2}{3}i\vartheta } } \right|}} \le \begin{cases} \left|\sec \vartheta \right| & \; \text{ if } \; { - \frac{\pi }{2} < \vartheta  < 0 \; \text{ or } \; \pi  < \vartheta < \frac{{3\pi }}{2}}, \\ 1 & \; \text{ if } \; {0 \le \vartheta  \le \pi }, \end{cases}
\end{equation}
for any $r>0$. Employing this inequality in \eqref{eq23} yields the bound
\[
\left| {R_N^{\left( H \right)} \left( \nu  \right)} \right| \le  - \frac{{\sec \left( {\theta  - \varphi } \right)}}{{\cos ^{\frac{{2N + 3}}{3}} \varphi }}\frac{1}{{\sqrt 3 \pi \left|\nu \right|^{\frac{{2N + 3}}{3}} }}\int_0^{ + \infty } {s^{\frac{{2N}}{3}} e^{ - s} \frac{2}{\pi }\int_{\mathscr{P}_ + ^{\left( { - 1} \right)} \left( {\frac{\pi }{2}} \right)} {\left( {\frac{i}{{f\left( u \right)}}} \right)^{\frac{{2N + 3}}{3}} dx} ds} .
\]
The double integral can be evaluated in terms of $B_{2N + 2} \left( 0 \right)$ since
\begin{align*}
\Gamma \left( {\frac{{2N + 3}}{3}} \right)6^{\frac{{2N + 3}}{3}} B_{2N + 2} \left( 0 \right) & = \left( { - 1} \right)^{N + 1} \int_0^{ + \infty } {t^{\frac{{2N}}{3}} e^{ - 2\pi t} iH_{it}^{\left( 1 \right)} \left( {it} \right)dt} 
\\ & = \left( { - 1} \right)^{N + 1} \int_0^{ + \infty } {t^{\frac{{2N}}{3}} \frac{1}{\pi }\int_{\mathscr{P}^{\left( { - 1} \right)} \left( {\frac{\pi }{2}} \right)} {e^{itf\left( u \right)} du} dt} 
\\ & = \left( { - 1} \right)^{N + 1} \int_0^{ + \infty } {s^{\frac{{2N}}{3}} e^{ - s} \frac{2}{\pi }\int_{\mathscr{P}_ + ^{\left( { - 1} \right)} \left( {\frac{\pi }{2}} \right)} {\left( {\frac{i}{{f\left( u \right)}}} \right)^{\frac{{2N + 3}}{3}} dx} ds} .
\end{align*}
Noting that $\left( { - 1} \right)^{N + 1} B_{2N + 2} \left( 0 \right) = \left| {B_{2N + 2} \left( 0 \right)} \right|$, we arrive at the estimate
\begin{equation}\label{eq25}
\left| {R_N^{\left( H \right)} \left( \nu  \right)} \right| \le  - \frac{{\sec \left( {\theta  - \varphi } \right)}}{{\cos ^{\frac{{2N + 3}}{3}} \varphi }}\frac{2}{{3\pi }}6^{\frac{{2N + 3}}{3}} \left| {B_{2N + 2} \left( 0 \right)} \right|\frac{{\sqrt 3 }}{2}\frac{{\Gamma \left( {\frac{{2N + 3}}{3}} \right)}}{{\left|\nu\right| ^{\frac{{2N + 3}}{3}} }} .
\end{equation}
The minimisation of the factor $ - \sec \left( {\theta  - \varphi } \right)\cos ^{ - \frac{{2N + 3}}{3}} \varphi$ as a function of $\varphi$ can be done using a lemma of Meijer \cite[pp. 953--954]{Meijer}. In our case, Meijer's lemma gives that the minimising value $\varphi = \varphi^\ast$ in \eqref{eq25}, is the unique solution of the equation
\begin{equation}\label{eq26}
\sin \left( {\theta  - 2\varphi ^\ast} \right) = \frac{N}{{N + 3}}\sin \theta ,
\end{equation}
that satisfies $ - \frac{{3\pi }}{2} + \theta  < \varphi ^\ast < \frac{\pi }{2} $ if $ \frac{{3\pi }}{2} \le \theta  < 2\pi$, and $0 < \varphi ^\ast  <  - \pi  + \theta$ if $\pi  < \theta  < \frac{{3\pi }}{2}$. To obtain the bound for the sector $-\pi <\theta < 0$, we rotate the path of integration in \eqref{eq22} by $-\frac{\pi}{2}<\varphi<0$, and an argument similar to the above one shows that
\[
\left| {R_N^{\left( H \right)} \left( \nu  \right)} \right| \le  \frac{{\sec \left( {\theta  - \varphi^\ast } \right)}}{{\cos ^{\frac{{2N + 3}}{3}} \varphi^\ast }}\frac{2}{{3\pi }}6^{\frac{{2N + 3}}{3}} \left| {B_{2N + 2} \left( 0 \right)} \right|\frac{{\sqrt 3 }}{2}\frac{{\Gamma \left( {\frac{{2N + 3}}{3}} \right)}}{{\left|\nu\right| ^{\frac{{2N + 3}}{3}} }} ,
\]
where $\varphi^\ast$ is the unique solution of the equation \eqref{eq26}, that satisfies $ - \frac{\pi}{2} < \varphi ^\ast < \frac{\pi}{2} +\theta$ if $-\pi < \theta  \le -\frac{\pi}{2}$, and $\theta < \varphi ^\ast  < 0$ if $-\frac{\pi}{2} < \theta  < 0$.

A very similar argument applied for the integral \eqref{eq28}, gives that
\[
\left| {R_N^{\left( H \right)} \left( \nu  \right)} \right| \le \frac{{\left| {\sec \left( {\theta  - \varphi ^{ \ast\ast }} \right)} \right|}}{{\cos ^{\frac{{2N + 1}}{3}} \varphi ^{ \ast\ast } }}\frac{2}{{3\pi }}6^{\frac{{2N + 1}}{3}} \left| {B_{2N} \left( 0 \right)} \right|\frac{{\sqrt 3 }}{2}\frac{{\Gamma \left( {\frac{{2N + 1}}{3}} \right)}}{{\left|\nu\right| ^{\frac{{2N + 1}}{3}} }}
\]
if $N\equiv 2 \mod 3$, where $\varphi^{ \ast\ast }$ is the unique solution of the equation
\begin{equation}\label{eq30}
\sin \left( {\theta  - 2\varphi ^{ \ast\ast }} \right) = \frac{N-1}{N + 2}\sin \theta ,
\end{equation}
that satisfies $ - \frac{{3\pi }}{2} + \theta  < \varphi ^{ \ast\ast } < \frac{\pi }{2} $ if $ \frac{{3\pi }}{2} \le \theta  < 2\pi$; $0 < \varphi ^{ \ast\ast }  <  - \pi  + \theta$ if $\pi  < \theta  < \frac{{3\pi }}{2}$; $ - \frac{\pi}{2} < \varphi ^{ \ast\ast } < \frac{\pi}{2} +\theta$ if $-\pi < \theta  \le -\frac{\pi}{2}$; and $\theta < \varphi ^{ \ast\ast }  < 0$ if $-\frac{\pi}{2} < \theta  < 0$.

Finally, from \eqref{eq29}, one finds that
\[
\left| {R_N^{\left( H \right)} \left( \nu  \right)} \right| \le \frac{{\left| {\sec \left( {\theta  - \varphi ^ \ast } \right)} \right|}}{{\cos ^{\frac{{2N + 1}}{3}} \varphi ^\ast }}\frac{2}{{3\pi }}6^{\frac{{2N + 1}}{3}} \left| {B_{2N} \left( 0 \right)} \right|\frac{{\sqrt 3 }}{2}\frac{{\Gamma \left( {\frac{{2N + 1}}{3}} \right)}}{{\left| \nu  \right|^{\frac{{2N + 1}}{3}} }} + \frac{{\left| {\sec \left( {\theta  - \varphi ^{\ast\ast} } \right)} \right|}}{{\cos ^{\frac{{2N + 3}}{3}} \varphi ^{ \ast\ast } }}\frac{2}{{3\pi }}6^{\frac{{2N + 3}}{3}} \left| {B_{2N + 2} \left( 0 \right)} \right|\frac{{\sqrt 3 }}{2}\frac{{\Gamma \left( {\frac{{2N + 3}}{3}} \right)}}{{\left| \nu  \right|^{\frac{{2N + 3}}{3}} }}
\]
if $N\equiv 0 \mod 3$, where $\varphi^\ast$ and $\varphi^{\ast \ast}$ are the unique solutions of the equations \eqref{eq26} and \eqref{eq30}, that satisfy $ - \frac{{3\pi }}{2} + \theta  < \varphi ^ \ast,\varphi ^{ \ast\ast } < \frac{\pi }{2} $ if $ \frac{{3\pi }}{2} \le \theta  < 2\pi$; $0 < \varphi ^ \ast,\varphi ^{ \ast\ast }  <  - \pi  + \theta$ if $\pi  < \theta  < \frac{{3\pi }}{2}$; $ - \frac{\pi}{2} < \varphi ^ \ast,\varphi ^{ \ast\ast } < \frac{\pi}{2} +\theta$ if $-\pi < \theta  \le -\frac{\pi}{2}$; and $\theta < \varphi ^ \ast,\varphi ^{ \ast\ast }  < 0$ if $-\frac{\pi}{2} < \theta  < 0$.

We can make our bounds simpler if $\arg \nu$ is close to $\frac{3\pi}{2}$ or $-\frac{\pi}{2}$ (i.e., close to the Stokes rays) as follows. Consider first the case $N\equiv 1 \mod 3$. When $\arg \nu = \frac{3\pi}{2}$, the minimising value $\varphi^\ast$ is given explicitly by
\[
\varphi^\ast = \mathop{\text{arccot}} \left( {\sqrt {\frac{{2N + 3}}{3}} } \right),
\]
and therefore we have
\[
 - \frac{{\sec \left( {\theta  - \varphi^\ast } \right)}}{{\cos ^{\frac{{2N + 3}}{3}} \varphi^\ast }} \le  - \frac{{\sec \left( {\frac{{3\pi }}{2} - \varphi^\ast } \right)}}{{\cos ^{\frac{{2N + 3}}{3}} \varphi^\ast }} = \frac{1}{{\sqrt 3 }}\left( {1 + \frac{3}{{2N + 3}}} \right)^{\frac{{N + 3}}{3}} \sqrt {2N + 3}  \le \sqrt {\frac{e}{3}\left( {2N + \frac{9}{2}} \right)} ,
\]
as long as $\pi  + \varphi^\ast  < \theta  \le \frac{{3\pi }}{2}$. We also have
\[
\sqrt {\frac{e}{3}\left( {2N + \frac{9}{2}} \right)}  \ge \sqrt {\frac{e}{3}\left( {2 + \frac{9}{2}} \right)}  \ge \left| {\sec \theta } \right|,
\]
for $\pi  < \theta  \le \pi  + \varphi^\ast \le \pi  + \mathop{\text{arccot}} \left( {\sqrt {\frac{{2 + 3}}{3}} } \right)$, therefore by \eqref{eq31} and \eqref{eq25} (with $\varphi = \varphi^\ast$), we deduce
\[
\left| {R_N^{\left( H \right)} \left( \nu  \right)} \right| \le \sqrt {\frac{e}{3}\left( {2N + \frac{9}{2}} \right)} \frac{2}{{3\pi }}6^{\frac{{2N + 3}}{3}} \left| {B_{2N + 2} \left( 0 \right)} \right|\frac{{\sqrt 3 }}{2}\frac{{\Gamma \left( {\frac{{2N + 3}}{3}} \right)}}{{\left| \nu  \right|^{\frac{{2N + 3}}{3}} }},
\]
provided that $\pi <\arg \nu \leq \frac{3\pi}{2}$ and $N\equiv 1 \mod 3$. A very similar argument shows that this bound is also valid when $-\frac{\pi}{2} \leq \arg \nu < 0$. Likewise, for the case $N\equiv 2 \mod 3$, we find
\[
\left| {R_N^{\left( H \right)} \left( \nu  \right)} \right| \le \sqrt {\frac{e}{3}\left( {2N + \frac{5}{2}} \right)} \frac{2}{{3\pi }}6^{\frac{{2N + 1}}{3}} \left| {B_{2N} \left( 0 \right)} \right|\frac{{\sqrt 3 }}{2}\frac{{\Gamma \left( {\frac{{2N + 1}}{3}} \right)}}{{\left| \nu  \right|^{\frac{{2N + 1}}{3}} }}
\]
if $\pi <\arg \nu \leq \frac{3\pi}{2}$ or $-\frac{\pi}{2} \leq \arg \nu < 0$. Finally, if $N\equiv 0\mod 3$, we have
\begin{gather}\label{eq70}
\begin{split}
\left| {R_N^{\left( H \right)} \left( \nu  \right)} \right| \le \; & \sqrt {\frac{e}{3}\left( {2N + \frac{5}{2}} \right)} \frac{2}{{3\pi }}6^{\frac{{2N + 1}}{3}} \left| {B_{2N} \left( 0 \right)} \right|\frac{{\sqrt 3 }}{2}\frac{{\Gamma \left( {\frac{{2N + 1}}{3}} \right)}}{{\left| \nu  \right|^{\frac{{2N + 1}}{3}} }} \\ & + \sqrt {\frac{e}{3}\left( {2N + \frac{9}{2}} \right)} \frac{2}{{3\pi }}6^{\frac{{2N + 3}}{3}} \left| {B_{2N + 2} \left( 0 \right)} \right|\frac{{\sqrt 3 }}{2}\frac{{\Gamma \left( {\frac{{2N + 3}}{3}} \right)}}{{\left| \nu  \right|^{\frac{{2N + 3}}{3}} }},
\end{split}
\end{gather}
when $\pi <\arg \nu \leq \frac{3\pi}{2}$ or $-\frac{\pi}{2} \leq \arg \nu < 0$ and $N\geq 3$. If we replace $\frac{5}{2}$ under the first square root by $\frac{9}{2}$, the estimate becomes true for $N=0$ too. Error bounds of the same form were established in the paper \cite{Nemes}, but with the worse constants $\frac{9}{2}$ and $\frac{13}{2}$ in place of $\frac{5}{2}$ and $\frac{9}{2}$, respectively.

The appearance of the factors proportional to $\sqrt{N}$ in these bounds may give the impression that these estimates are unrealistic for large $N$, but this is not the case. Applying the asymptotic form of the coefficients $B_{2N} \left( 0 \right)$ (see Section \ref{section4}), it can readily be shown that, in particular, \eqref{eq70} implies
\[
R_{3N}^{\left( H \right)} \left( \nu  \right) = \mathcal{O}\left( {\frac{{\left( {2N} \right)^{\frac{1}{6}} }}{{\left( {2\pi } \right)^{2N} }}\left( {\frac{{\Gamma \left( {2N + \frac{1}{3}} \right)}}{{\left| \nu  \right|^{2N + \frac{1}{3}} }} + \frac{{\Gamma \left( {2N + 1} \right)}}{{\left| \nu  \right|^{2N + 1} }}} \right)} \right)
\]
for large $N$. When the asymptotic series truncated optimally, i.e., $2N \approx 2\pi \left| \nu  \right|$, we find with the help of Stirling's formula that the above estimate is equivalent to
\[
R_{3N}^{\left( H \right)} \left( \nu  \right) = \mathcal{O}\left( {\left| \nu  \right|^{ - \frac{1}{3}} e^{ - 2\pi \left| \nu  \right|} } \right).
\]
Noting that $R_{3N}^{\left( H \right)} \left( \nu  \right) = R_{6N}^{\left( H \right)} \left( {\nu ,0} \right)$, this estimate is a special case of the exponentially improved version, given in Theorem \ref{thm3}. Therefore, our bounds are indeed realistic near the Stokes lines.

\subsection{Error bounds for the asymptotic series of $H_\nu ^{\left( 1 \right) \prime} \left( \nu  \right)$ and $H_\nu ^{\left( 2 \right)\prime} \left( \nu  \right)$} Since the remainder in the asymptotic series of $H_\nu ^{\left( 2 \right)\prime} \left( \nu  \right)$ is $ R_N^{\left( H' \right)} \left( {\nu e^{\pi i} } \right)$, it is enough to estimate $R_N^{\left( H '\right)} \left( \nu  \right)$. We will estimate $R_N^{\left( H '\right)} \left( \nu  \right)$ in the large sector $-\pi <\arg \nu <2\pi$. The representation \eqref{eq32} is valid only for $-\frac{\pi}{2} <\arg \nu <\frac{3\pi}{2}$, for the larger sector, we use \eqref{eq33} as the definition of $R_N^{\left( H' \right)} \left( \nu  \right)$. Throughout this subsection, we assume that $N\geq 1$. Algebraic manipulation of \eqref{eq32} shows that
\begin{align}
R_N^{\left( {H'} \right)} \left( \nu  \right) = \; & \frac{{\left( { - 1} \right)^N }}{{\sqrt 3 \pi \nu ^{\frac{{2N + 2}}{3}} }}\int_0^{ + \infty } { \frac{t^{\frac{{2N - 1}}{3}} e^{ - 2\pi t} e^{\frac{{2\left( {2N + 2} \right)\pi i}}{3}}}{{\left( {1 + \left( {t/\nu } \right)^{\frac{2}{3}} e^{\frac{{2\pi i}}{3}} } \right)\left( {1 + \left( {t/\nu } \right)^{\frac{2}{3}} } \right)}}H_{it}^{\left( 1 \right) \prime  } \left( {it} \right)dt} ,
\label{eq34} \\ R_N^{\left( {H'} \right)} \left( \nu  \right) = \; & \frac{{\left( { - 1} \right)^{N + 1} }}{{\sqrt 3 \pi \nu ^{\frac{{2N + 2}}{3}} }}\int_0^{ + \infty } { \frac{t^{\frac{{2N - 1}}{3}} e^{ - 2\pi t} e^{\frac{{2\left( {2N + 2} \right)\pi i}}{3}}}{{\left( {1 + \left( {t/\nu } \right)^{\frac{2}{3}} e^{\frac{{2\pi i}}{3}} } \right)\left( {1 + \left( {t/\nu } \right)^{\frac{2}{3}} } \right)}}H_{it}^{\left( 1 \right) \prime  } \left( {it} \right)dt} \label{eq40} \\ & + \frac{{\left( { - 1} \right)^{N + 1} }}{{\sqrt 3 \pi \nu ^{\frac{{2N + 4}}{3}} }}\int_0^{ + \infty } { \frac{{t^{\frac{{2N + 1}}{3}} e^{ - 2\pi t} e^{\frac{{2\left( {2N + 2} \right)\pi i}}{3}} e^{\frac{\pi }{3}i} }}{{\left( {1 + \left( {t/\nu } \right)^{\frac{2}{3}} e^{\frac{{2\pi i}}{3}} } \right)\left( {1 + \left( {t/\nu } \right)^{\frac{2}{3}} } \right)}}H_{it}^{\left( 1 \right) \prime  } \left( {it} \right)dt} ,
\nonumber \\ R_N^{\left( {H'} \right)} \left( \nu  \right) = \; & \frac{{\left( { - 1} \right)^N }}{{\sqrt 3 \pi \nu ^{\frac{{2N + 4}}{3}} }}\int_0^{ + \infty } { \frac{{t^{\frac{{2N + 1}}{3}} e^{ - 2\pi t} e^{\frac{{2\left( {2N + 2} \right)\pi i}}{3}} e^{\frac{\pi }{3}i} }}{{\left( {1 + \left( {t/\nu } \right)^{\frac{2}{3}} e^{\frac{{2\pi i}}{3}} } \right)\left( {1 + \left( {t/\nu } \right)^{\frac{2}{3}} } \right)}}H_{it}^{\left( 1 \right) \prime  } \left( {it} \right)dt} \label{eq39}
\end{align}
according to whether $N\equiv 0 \mod 3$, $N\equiv 1 \mod 3$ or $N\equiv 2 \mod 3$, respectively. We consider the case $N\equiv 0 \mod 3$. First we derive a bound suitable when $-\frac{\pi}{2}<\arg \nu <\frac{3\pi}{2}$. We introduce the integral representation of $e^{ - 2\pi t}  H_{it}^{\left( 1 \right)\prime} \left( {it} \right)$ and perform the change of variable $is/f\left( u \right) = t$ in \eqref{eq34}, to deduce
\begin{align*}
& R_N^{\left( {H'} \right)} \left( \nu  \right) = \frac{{\left( { - 1} \right)^N }}{{\sqrt 3 \pi \nu ^{\frac{{2N + 2}}{3}} }}\int_0^{ + \infty } { \frac{t^{\frac{{2N - 1}}{3}} e^{\frac{{2\left( {2N + 2} \right)\pi i}}{3}}}{{\left( {1 + \left( {t/\nu } \right)^{\frac{2}{3}} e^{\frac{{2\pi i}}{3}} } \right)\left( {1 + \left( {t/\nu } \right)^{\frac{2}{3}} } \right)}}\frac{1}{{\pi i}}\int_{\mathscr{P}^{\left( { - 1} \right)} \left( {\frac{\pi }{2}} \right)} {e^{itf\left( u \right)} \sinh udu} dt} 
\\ & = \frac{{\left( { - 1} \right)^N }}{{\sqrt 3 \pi \nu ^{\frac{{2N + 2}}{3}} }}\int_0^{ + \infty } {s^{\frac{{2N - 1}}{3}} e^{ - s} e^{\frac{{2\left( {2N + 2} \right)\pi i}}{3}} \frac{1}{{\pi i}}\int_{\mathscr{P}^{\left( { - 1} \right)} \left( {\frac{\pi }{2}} \right)} {\frac{{\left( {i/f\left( u \right)} \right)^{\frac{{2N + 2}}{3}} }}{{\left( {1 + \left( {is/\nu f\left( u \right)} \right)^{\frac{2}{3}} e^{\frac{{2\pi i}}{3}} } \right)\left( {1 + \left( {is/\nu f\left( u \right)} \right)^{\frac{2}{3}} } \right)}}\sinh udu} ds} .
\end{align*}
With the notation of the previous subsection, we can write the contour integral as
\begin{align*}
& \frac{1}{{\pi i}}\int_{\mathscr{P}^{\left( { - 1} \right)} \left( {\frac{\pi }{2}} \right)} {\frac{{\left( {i/f\left( u \right)} \right)^{\frac{{2N + 2}}{3}} }}{{\left( {1 + \left( {is/\nu f\left( u \right)} \right)^{\frac{2}{3}} e^{\frac{{2\pi i}}{3}} } \right)\left( {1 + \left( {is/\nu f\left( u \right)} \right)^{\frac{2}{3}} } \right)}}\sinh udu}
\\ = \; & \frac{2}{\pi }\int_{\mathscr{P}_ + ^{\left( { - 1} \right)} \left( {\frac{\pi }{2}} \right)} {\frac{{\left( {i/f\left( u \right)} \right)^{\frac{{2N + 2}}{3}} }}{{\left( {1 + \left( {is/\nu f\left( u \right)} \right)^{\frac{2}{3}} e^{\frac{{2\pi i}}{3}} } \right)\left( {1 + \left( {is/\nu f\left( u \right)} \right)^{\frac{2}{3}} } \right)}}\sinh x\cos ydy} 
\\ & + \frac{2}{\pi }\int_{\mathscr{P}_ + ^{\left( { - 1} \right)} \left( {\frac{\pi }{2}} \right)} {\frac{{\left( {i/f\left( u \right)} \right)^{\frac{{2N + 2}}{3}} }}{{\left( {1 + \left( {is/\nu f\left( u \right)} \right)^{\frac{2}{3}} e^{\frac{{2\pi i}}{3}} } \right)\left( {1 + \left( {is/\nu f\left( u \right)} \right)^{\frac{2}{3}} } \right)}}\cosh x\sin ydx} ,
\end{align*}
whence
\begin{gather}\label{eq35}
\begin{split}
R_N^{\left( {H'} \right)} \left( \nu  \right) = \; & \frac{{\left( { - 1} \right)^N }}{{\sqrt 3 \pi \nu ^{\frac{{2N + 2}}{3}} }}\int_0^{ + \infty } s^{\frac{{2N - 1}}{3}} e^{ - s} e^{\frac{{2\left( {2N + 2} \right)\pi i}}{3}} \\ & \times \frac{2}{\pi }\int_{\mathscr{P}_ + ^{\left( { - 1} \right)} \left( {\frac{\pi }{2}} \right)} {\frac{{\left( {i/f\left( u \right)} \right)^{\frac{{2N + 2}}{3}} }}{{\left( {1 + \left( {is/\nu f\left( u \right)} \right)^{\frac{2}{3}} e^{\frac{{2\pi i}}{3}} } \right)\left( {1 + \left( {is/\nu f\left( u \right)} \right)^{\frac{2}{3}} } \right)}}\sinh x\cos ydy} ds
\\ & + \frac{{\left( { - 1} \right)^N }}{{\sqrt 3 \pi \nu ^{\frac{{2N + 2}}{3}} }}\int_0^{ + \infty } s^{\frac{{2N - 1}}{3}} e^{ - s} e^{\frac{{2\left( {2N + 2} \right)\pi i}}{3}} \\ & \times \frac{2}{\pi }\int_{\mathscr{P}_ + ^{\left( { - 1} \right)} \left( {\frac{\pi }{2}} \right)} {\frac{{\left( {i/f\left( u \right)} \right)^{\frac{{2N + 2}}{3}} }}{{\left( {1 + \left( {is/\nu f\left( u \right)} \right)^{\frac{2}{3}} e^{\frac{{2\pi i}}{3}} } \right)\left( {1 + \left( {is/\nu f\left( u \right)} \right)^{\frac{2}{3}} } \right)}}\cosh x\sin ydx} ds .
\end{split}
\end{gather}
Similarly, we have the following expression for the coefficients $D_{2N + 1} \left( 0 \right)$:
\begin{align*}
\Gamma \left( {\frac{{2N + 2}}{3}} \right)6^{\frac{{2N + 2}}{3}} D_{2N + 1} \left( 0 \right) = \; &  \left( { - 1} \right)^{N + 1} \int_0^{ + \infty } {t^{\frac{{2N - 1}}{3}} e^{ - 2\pi t} H_{it}^{\left( 1 \right) \prime  } \left( {it} \right)dt} 
\\ = \; &  \left( { - 1} \right)^{N + 1} \int_0^{ + \infty } {t^{\frac{{2N - 1}}{3}} \frac{1}{{\pi i}}\int_{\mathscr{P}^{\left( { - 1} \right)} \left( {\frac{\pi }{2}} \right)} {e^{itf\left( u \right)} \sinh udu} dt}\\
 =\; & \left( { - 1} \right)^{N + 1} \int_0^{ + \infty } {s^{\frac{{2N - 1}}{3}} e^{ - s} \frac{2}{\pi }\int_{\mathscr{P}_ + ^{\left( { - 1} \right)} \left( {\frac{\pi }{2}} \right)} {\left( {\frac{i}{{f\left( u \right)}}} \right)^{\frac{{2N + 2}}{3}} \sinh x\cos ydy} ds} \\ & + \left( { - 1} \right)^{N + 1} \int_0^{ + \infty } {s^{\frac{{2N - 1}}{3}} e^{ - s} \frac{2}{\pi }\int_{\mathscr{P}_ + ^{\left( { - 1} \right)} \left( {\frac{\pi }{2}} \right)} {\left( {\frac{i}{{f\left( u \right)}}} \right)^{\frac{{2N + 2}}{3}} \cosh x\sin ydx} ds} .
\end{align*}
We note that along the path $\mathscr{P}_ + ^{\left( { - 1} \right)} \left( {\frac{\pi }{2}} \right)$, we have $x>0$ and $- 2\pi  \le y <  - \frac{{3\pi }}{2}$ and both of them are increasing, and therefore the integrals in the last expression are all non-negative. In particular, we have $\left( { - 1} \right)^{N + 1} D_{2N + 1} \left( 0 \right) = \left| {D_{2N + 1} \left( 0 \right)} \right|$. Applying the inequality \eqref{eq36} for \eqref{eq35} and using the above representation of $D_{2N + 1} \left( 0 \right)$, we obtain the error bound
\[
\left| {R_N^{\left( {H'} \right)} \left( \nu  \right)} \right| \le \frac{2}{{3\pi }}6^{\frac{{2N + 2}}{3}} \left| {D_{2N + 1} \left( 0 \right)} \right|\frac{{\sqrt 3 }}{2}\frac{{\Gamma \left( {\frac{{2N + 2}}{3}} \right)}}{{\left| \nu  \right|^{\frac{{2N + 2}}{3}} }} \begin{cases} \left|\sec \theta \right| & \; \text{ if } \; { - \frac{\pi }{2} < \theta  < 0 \; \text{ or } \; \pi  < \theta  < \frac{{3\pi }}{2}}, \\ 1 & \; \text{ if } \; {0 \le \theta  \le \pi }, \end{cases}
\]
if $N\equiv 0 \mod 3$. To derive the corresponding bounds for the sectors $-\pi <\arg \nu<0$ and $\pi <\arg \nu <2\pi$, we again use the path rotational argument and Meijer's lemma as in the previous subsection, to find
\[
\left| {R_N^{\left( {H'} \right)} \left( \nu  \right)} \right| \le \frac{{\left| {\sec \left( {\theta  - \varphi^\ast } \right)} \right|}}{{\cos ^{\frac{{2N + 2}}{3}} \varphi^\ast }}\frac{2}{{3\pi }}6^{\frac{{2N + 2}}{3}} \left| {D_{2N + 1} \left( 0 \right)} \right|\frac{{\sqrt 3 }}{2}\frac{{\Gamma \left( {\frac{{2N + 2}}{3}} \right)}}{{\left| \nu  \right|^{\frac{{2N + 2}}{3}} }},
\]
if $N\equiv 0 \mod 3$, where $\varphi^\ast$ is the unique solution of the equation
\[
\sin \left( {\theta  - 2\varphi^\ast} \right) = \frac{2N-1}{2N + 5}\sin \theta ,
\]
that satisfies $ - \frac{{3\pi }}{2} + \theta  < \varphi^\ast < \frac{\pi }{2} $ if $ \frac{{3\pi }}{2} \le \theta  < 2\pi$; $0 < \varphi^\ast <  - \pi  + \theta$ if $\pi  < \theta  < \frac{{3\pi }}{2}$; $ - \frac{\pi}{2} < \varphi^\ast < \frac{\pi}{2} +\theta$ if $-\pi < \theta  \le -\frac{\pi}{2}$; and $\theta < \varphi^\ast < 0$ if $-\frac{\pi}{2} < \theta  < 0$.

With the same argument applied for the integral \eqref{eq39}, the analogous bounds for the case $N\equiv 2 \mod 3$ are found to be
\[
\left| {R_N^{\left( {H'} \right)} \left( \nu  \right)} \right| \le \frac{2}{{3\pi }}6^{\frac{{2N + 4}}{3}} \left| {D_{2N + 3} \left( 0 \right)} \right|\frac{{\sqrt 3 }}{2}\frac{{\Gamma \left( {\frac{{2N + 4}}{3}} \right)}}{{\left| \nu  \right|^{\frac{{2N + 4}}{3}} }} \begin{cases} \left|\sec \theta \right| & \; \text{ if } \; { - \frac{\pi }{2} < \theta  < 0 \; \text{ or } \; \pi  < \theta  < \frac{{3\pi }}{2}}, \\ 1 & \; \text{ if } \; {0 \le \theta  \le \pi }, \end{cases}
\]
and
\[
\left| {R_N^{\left( {H'} \right)} \left( \nu  \right)} \right| \le \frac{{\left| {\sec \left( {\theta  - \varphi ^{\ast\ast} } \right)} \right|}}{{\cos ^{\frac{{2N + 4}}{3}} \varphi^{\ast\ast} }}\frac{2}{{3\pi }}6^{\frac{{2N + 4}}{3}} \left| {D_{2N + 3} \left( 0 \right)} \right|\frac{{\sqrt 3 }}{2}\frac{{\Gamma \left( {\frac{{2N + 4}}{3}} \right)}}{{\left| \nu  \right|^{\frac{{2N + 4}}{3}} }},
\]
where $\varphi^{\ast\ast}$ is the unique solution of the equation
\[
\sin \left( {\theta  - 2\varphi^{\ast\ast}} \right) = \frac{2N+1}{2N + 7}\sin \theta ,
\]
that satisfies $ - \frac{{3\pi }}{2} + \theta  < \varphi^{\ast\ast} < \frac{\pi }{2} $ if $ \frac{{3\pi }}{2} \le \theta  < 2\pi$; $0 < \varphi^{\ast\ast} <  - \pi  + \theta$ if $\pi  < \theta  < \frac{{3\pi }}{2}$; $ - \frac{\pi}{2} < \varphi^{\ast\ast} < \frac{\pi}{2} +\theta$ if $-\pi < \theta  \le -\frac{\pi}{2}$; and $\theta < \varphi^{\ast\ast} < 0$ if $-\frac{\pi}{2} < \theta  < 0$.

Finally, from \eqref{eq40}, one finds that for $N\equiv 1 \mod 3$, we have
\begin{align*}
\left| {R_N^{\left( {H'} \right)} \left( \nu  \right)} \right| \le \; & \left(\frac{2}{{3\pi }}6^{\frac{{2N + 2}}{3}} \left| {D_{2N + 1} \left( 0 \right)} \right|\frac{{\sqrt 3 }}{2}\frac{{\Gamma \left( {\frac{{2N + 2}}{3}} \right)}}{{\left| \nu  \right|^{\frac{{2N + 2}}{3}} }}\right. \\ & \left.+ \frac{2}{{3\pi }}6^{\frac{{2N + 4}}{3}} \left| {D_{2N + 3} \left( 0 \right)} \right|\frac{{\sqrt 3 }}{2}\frac{{\Gamma \left( {\frac{{2N + 4}}{3}} \right)}}{{\left| \nu  \right|^{\frac{{2N + 4}}{3}} }}\right) \begin{cases} \left|\sec \theta \right| & \; \text{ if } \; { - \frac{\pi }{2} < \theta  < 0 \; \text{ or } \; \pi  < \theta  < \frac{{3\pi }}{2}}, \\ 1 & \; \text{ if } \; {0 \le \theta  \le \pi }, \end{cases}
\end{align*}
and
\begin{align*}
\left| {R_N^{\left( {H'} \right)} \left( \nu  \right)} \right| \le \; & \frac{{\left| {\sec \left( {\theta  - \varphi^\ast} \right)} \right|}}{{\cos ^{\frac{{2N + 2}}{3}} \varphi^\ast}}\frac{2}{{3\pi }}6^{\frac{{2N + 2}}{3}} \left| {D_{2N + 1} \left( 0 \right)} \right|\frac{{\sqrt 3 }}{2}\frac{{\Gamma \left( {\frac{{2N + 2}}{3}} \right)}}{{\left| \nu  \right|^{\frac{{2N + 2}}{3}} }} \\ & + \frac{{\left| {\sec \left( {\theta  - \varphi^{\ast\ast} } \right)} \right|}}{{\cos ^{\frac{{2N + 4}}{3}} \varphi^{\ast\ast} }}\frac{2}{{3\pi }}6^{\frac{{2N + 4}}{3}} \left| {D_{2N + 3} \left( 0 \right)} \right|\frac{{\sqrt 3 }}{2}\frac{{\Gamma \left( {\frac{{2N + 4}}{3}} \right)}}{{\left| \nu  \right|^{\frac{{2N + 4}}{3}} }},
\end{align*}
where $\varphi^\ast$ and $\varphi^{\ast \ast}$ are the same as above.

Again, it is possible to derive bounds that are suitable near the Stokes lines $\arg \nu = \frac{3\pi}{2}$ and $\arg \nu = -\frac{\pi}{2}$. Since the argument is the same as in the previous subsection, we just state the final results. If $\pi <\arg \nu \leq \frac{3\pi}{2}$ or $-\frac{\pi}{2} \leq \arg \nu < 0$, then
\[
\left| {R_N^{\left( {H'} \right)} \left( \nu  \right)} \right| \le \sqrt {\frac{e}{3}\left( {2N + \frac{7}{2}} \right)} \frac{2}{{3\pi }}6^{\frac{{2N + 2}}{3}} \left| {D_{2N + 1} \left( 0 \right)} \right|\frac{{\sqrt 3 }}{2}\frac{{\Gamma \left( {\frac{{2N + 2}}{3}} \right)}}{{\left| \nu  \right|^{\frac{{2N + 2}}{3}} }},
\]
when $N\equiv 0 \mod 3$;
\[
\left| {R_N^{\left( {H'} \right)} \left( \nu  \right)} \right| \le \sqrt {\frac{e}{3}\left( {2N + \frac{{11}}{2}} \right)} \frac{2}{{3\pi }}6^{\frac{{2N + 4}}{3}} \left| {D_{2N + 3} \left( 0 \right)} \right|\frac{{\sqrt 3 }}{2}\frac{{\Gamma \left( {\frac{{2N + 4}}{3}} \right)}}{{\left| \nu  \right|^{\frac{{2N + 4}}{3}} }},
\]
if $N\equiv 2 \mod 3$; and
\begin{align*}
\left| {R_N^{\left( {H'} \right)} \left( \nu  \right)} \right| \le \; & \sqrt {\frac{e}{3}\left( {2N + \frac{7}{2}} \right)} \frac{2}{{3\pi }}6^{\frac{{2N + 2}}{3}} \left| {D_{2N + 1} \left( 0 \right)} \right|\frac{{\sqrt 3 }}{2}\frac{{\Gamma \left( {\frac{{2N + 2}}{3}} \right)}}{{\left| \nu  \right|^{\frac{{2N + 2}}{3}} }} \\ & +\sqrt {\frac{e}{3}\left( {2N + \frac{{11}}{2}} \right)} \frac{2}{{3\pi }}6^{\frac{{2N + 4}}{3}} \left| {D_{2N + 3} \left( 0 \right)} \right|\frac{{\sqrt 3 }}{2}\frac{{\Gamma \left( {\frac{{2N + 4}}{3}} \right)}}{{\left| \nu  \right|^{\frac{{2N + 4}}{3}} }},
\end{align*}
when $N\equiv 1 \mod 3$.

\subsection{Error bounds for the asymptotic series of $J_\nu \left( \nu  \right)$ and $Y_\nu \left( \nu  \right)$} First, we shall estimate $R_N^{\left( J\right)} \left( \nu  \right)$ in the sector $\left|\arg \nu \right|<\pi$. The representation \eqref{eq15} is valid only if $\left|\arg \nu \right|<\frac{\pi}{2}$, for the larger sector, we use \eqref{eq41} as the definition of $R_N^{\left( J\right)} \left( \nu  \right)$. Algebraic manipulation of \eqref{eq15} shows that
\begin{align}
R_N^{\left( J \right)} \left( \nu  \right) & = \frac{{\left( { - 1} \right)^N }}{{2\sqrt 3 \pi \nu ^{\frac{{2N + 1}}{3}} }}\int_0^{ + \infty } {t^{\frac{{2N - 2}}{3}} e^{ - 2\pi t} \frac{{1 - \left( {t/\nu } \right)^{\frac{4}{3}} }}{{1 + \left( {t/\nu } \right)^2 }}iH_{it}^{\left( 1 \right)} \left( {it} \right)dt} , \nonumber \\
R_N^{\left( J \right)} \left( \nu  \right) & = \frac{{\left( { - 1} \right)^N }}{{2\sqrt 3 \pi \nu ^{\frac{{2N + 3}}{3}} }}\int_0^{ + \infty } {t^{\frac{{2N}}{3}} e^{ - 2\pi t} \frac{{1 + \left( {t/\nu } \right)^{\frac{2}{3}} }}{{1 + \left( {t/\nu } \right)^2 }}iH_{it}^{\left( 1 \right)} \left( {it} \right)dt} , \nonumber \\
R_N^{\left( J \right)} \left( \nu  \right) & = \frac{{\left( { - 1} \right)^{N + 1} }}{{2\sqrt 3 \pi \nu ^{\frac{{2N + 1}}{3}} }}\int_0^{ + \infty } {t^{\frac{{2N - 2}}{3}} e^{ - 2\pi t} \frac{{1 + \left( {t/\nu } \right)^{\frac{2}{3}} }}{{1 + \left( {t/\nu } \right)^2 }}iH_{it}^{\left( 1 \right)} \left( {it} \right)dt} \nonumber
\end{align}
according to whether $N\equiv 0 \mod 3$, $N\equiv 1 \mod 3$ or $N\equiv 2 \mod 3$, respectively. The estimates which are suitable in the sector $\left|\arg \nu \right|<\frac{\pi}{2}$, were derived in the paper \cite{Nemes}, and they are given by
\begin{align*}
\left| {R_N^{\left( J \right)} \left( \nu  \right)} \right| \le \; & \frac{1}{{3\pi }}6^{\frac{{2N + 1}}{3}} \left| {B_{2N} \left( 0 \right)} \right|\frac{{\sqrt 3 }}{2}\frac{{\Gamma \left( {\frac{{2N + 1}}{3}} \right)}}{{\left| \nu  \right|^{\frac{{2N + 1}}{3}} }} \begin{cases} \left|\csc\left(2\theta\right)\right| & \; \text{ if } \; \frac{\pi}{4} < \left|\theta\right| <\frac{\pi}{2}, \\ 1 & \; \text{ if } \; \left|\theta\right| \leq \frac{\pi}{4} \end{cases} \\ & + \frac{1}{{3\pi }}6^{\frac{{2N + 5}}{3}} \left| {B_{2N + 4} \left( 0 \right)} \right|\frac{{\sqrt 3 }}{2}\frac{{\Gamma \left( {\frac{{2N + 5}}{3}} \right)}}{{\left| \nu  \right|^{\frac{{2N + 5}}{3}} }} \begin{cases} \left|\csc\left(2\theta\right)\right| & \; \text{ if } \; \frac{\pi}{4} < \left|\theta\right| <\frac{\pi}{2}, \\ 0 & \; \text{ if } \; \left|\theta\right| \leq \frac{\pi}{4}, \end{cases}
\end{align*}
for $N\equiv 0 \mod 3$;
\begin{align*}
\left| {R_N^{\left( J \right)} \left( \nu  \right)} \right| \le \; & \left( {\frac{1}{{3\pi }}6^{\frac{{2N + 3}}{3}} \left| {B_{2N + 2} \left( 0 \right)} \right|\frac{{\sqrt 3 }}{2}\frac{{\Gamma \left( {\frac{{2N + 3}}{3}} \right)}}{{\left| \nu  \right|^{\frac{{2N + 3}}{3}} }}}\right. \\ & \left.{+ \frac{1}{{3\pi }}6^{\frac{{2N + 5}}{3}} \left| {B_{2N + 4} \left( 0 \right)} \right|\frac{{\sqrt 3 }}{2}\frac{{\Gamma \left( {\frac{{2N + 5}}{3}} \right)}}{{\left| \nu  \right|^{\frac{{2N + 5}}{3}} }}} \right) \begin{cases} \left|\csc\left(2\theta\right)\right| & \; \text{ if } \; \frac{\pi}{4} < \left|\theta\right| <\frac{\pi}{2}, \\ 1 & \; \text{ if } \; \left|\theta\right| \leq \frac{\pi}{4}, \end{cases}
\end{align*}
if $N\equiv 1 \mod 3$; and
\begin{align*}
\left| {R_N^{\left( J \right)} \left( \nu  \right)} \right| \le \; & \left( {\frac{1}{{3\pi }}6^{\frac{{2N + 1}}{3}} \left| {B_{2N} \left( 0 \right)} \right|\frac{{\sqrt 3 }}{2}\frac{{\Gamma \left( {\frac{{2N + 1}}{3}} \right)}}{{\left| \nu  \right|^{\frac{{2N + 1}}{3}} }} }\right. \\ & \left.{+ \frac{1}{{3\pi }}6^{\frac{{2N + 3}}{3}} \left| {B_{2N + 2} \left( 0 \right)} \right|\frac{{\sqrt 3 }}{2}\frac{{\Gamma \left( {\frac{{2N + 3}}{3}} \right)}}{{\left| \nu  \right|^{\frac{{2N + 3}}{3}} }}} \right) \begin{cases} \left|\csc\left(2\theta\right)\right| & \; \text{ if } \; \frac{\pi}{4} < \left|\theta\right| <\frac{\pi}{2}, \\ 1 & \; \text{ if } \; \left|\theta\right| \leq \frac{\pi}{4}, \end{cases}
\end{align*}
for $N\equiv 2 \mod 3$. The reference \cite{Nemes} also contains sharp bounds for the case $\nu>0$. The derivation of the corresponding estimates for the sectors $\frac{\pi}{4}<\left|\arg \nu\right|<\pi$ is very similar to the arguments in the previous subsections, except instead of \eqref{eq36}, one has to use
\begin{equation}\label{eq66}
\frac{1}{{\left| {1 + re^{ - 2i\vartheta } } \right|}} \le \begin{cases} \left|\csc\left(2\vartheta\right)\right| & \; \text{ if } \; \frac{\pi}{4} < \left|\vartheta\right| <\frac{\pi}{2}, \\ 1 & \; \text{ if } \; \left|\vartheta\right| \leq \frac{\pi}{4}, \end{cases}
\end{equation}
for $r>0$, and another lemma of Meijer's \cite[p. 956]{Meijer}. The final results are
\begin{align*}
\left| {R_N^{\left( J \right)} \left( \nu  \right)} \right| \le \; & \frac{{\left| {\csc \left( {2\left( {\theta  - \varphi^\ast} \right)} \right)} \right|}}{{\cos ^{\frac{{2N + 1}}{3}} \varphi^\ast}}\frac{1}{{3\pi }}6^{\frac{{2N + 1}}{3}} \left| {B_{2N} \left( 0 \right)} \right|\frac{{\sqrt 3 }}{2}\frac{{\Gamma \left( {\frac{{2N + 1}}{3}} \right)}}{{\left| \nu  \right|^{\frac{{2N + 1}}{3}} }} \\ & + \frac{{\left| {\csc \left( {2\left( {\theta  - \varphi^{\ast\ast\ast} } \right)} \right)} \right|}}{{\cos ^{\frac{{2N + 5}}{3}} \varphi^{\ast\ast\ast}}}\frac{1}{{3\pi }}6^{\frac{{2N + 5}}{3}} \left| {B_{2N + 4} \left( 0 \right)} \right|\frac{{\sqrt 3 }}{2}\frac{{\Gamma \left( {\frac{{2N + 5}}{3}} \right)}}{{\left| \nu  \right|^{\frac{{2N + 5}}{3}} }},
\end{align*}
for $N\equiv 0 \mod 3$;
\begin{align*}
\left| {R_N^{\left( J \right)} \left( \nu  \right)} \right| \le \; & \frac{{\left| {\csc \left( {2\left( {\theta  - \varphi^{\ast\ast} } \right)} \right)} \right|}}{{\cos ^{\frac{{2N + 3}}{3}} \varphi^{\ast\ast} }}\frac{1}{{3\pi }}6^{\frac{{2N + 3}}{3}} \left| {B_{2N + 2} \left( 0 \right)} \right|\frac{{\sqrt 3 }}{2}\frac{{\Gamma \left( {\frac{{2N + 3}}{3}} \right)}}{{\left| \nu  \right|^{\frac{{2N + 3}}{3}} }} \\ & + \frac{{\left| {\csc \left( {2\left( {\theta  - \varphi^{\ast\ast\ast} } \right)} \right)} \right|}}{{\cos ^{\frac{{2N + 5}}{3}} \varphi^{\ast\ast\ast}}}\frac{1}{{3\pi }}6^{\frac{{2N + 5}}{3}} \left| {B_{2N + 4} \left( 0 \right)} \right|\frac{{\sqrt 3 }}{2}\frac{{\Gamma \left( {\frac{{2N + 5}}{3}} \right)}}{{\left| \nu  \right|^{\frac{{2N + 5}}{3}} }},
\end{align*}
if $N\equiv 1 \mod 3$; and
\begin{align*}
\left| {R_N^{\left( J \right)} \left( \nu  \right)} \right| \le \; & \frac{{\left| {\csc \left( {2\left( {\theta  - \varphi^\ast} \right)} \right)} \right|}}{{\cos ^{\frac{{2N + 1}}{3}} \varphi^\ast}}\frac{1}{{3\pi }}6^{\frac{{2N + 1}}{3}} \left| {B_{2N} \left( 0 \right)} \right|\frac{{\sqrt 3 }}{2}\frac{{\Gamma \left( {\frac{{2N + 1}}{3}} \right)}}{{\left| \nu  \right|^{\frac{{2N + 1}}{3}} }} \\ & + \frac{{\left| {\csc \left( {2\left( {\theta  - \varphi^{\ast\ast} } \right)} \right)} \right|}}{{\cos ^{\frac{{2N + 3}}{3}} \varphi^{\ast\ast} }}\frac{1}{{3\pi }}6^{\frac{{2N + 3}}{3}} \left| {B_{2N + 2} \left( 0 \right)} \right|\frac{{\sqrt 3 }}{2}\frac{{\Gamma \left( {\frac{{2N + 3}}{3}} \right)}}{{\left| \nu  \right|^{\frac{{2N + 3}}{3}} }},
\end{align*}
for $N\equiv 2 \mod 3$, where $\varphi^\ast$, $\varphi^{\ast \ast}$ and $\varphi^{\ast \ast \ast}$ are the unique solutions of the implicit equations
\[
\left( {2N + 7 } \right)\cos \left( {3\varphi^\ast - 2\theta} \right) = \left( {2N - 5 } \right)\cos \left( {\varphi^\ast- 2\theta } \right),
\]
\[
\left( {2N + 9 } \right)\cos \left( {3\varphi^{\ast\ast} - 2\theta} \right) = \left( {2N - 3 } \right)\cos \left( {\varphi^{\ast\ast}- 2\theta } \right)
\]
and
\[
\left( {2N + 11 } \right)\cos \left( {3\varphi^{\ast\ast\ast} - 2\theta} \right) = \left( {2N -1 } \right)\cos \left( {\varphi^{\ast\ast\ast}- 2\theta } \right),
\]
that satisfy $-\frac{\pi}{2}+\theta < \varphi^\ast,\varphi^{\ast\ast},\varphi^{\ast\ast\ast} <\frac{\pi}{2}$ if $\frac{3\pi}{4} \leq \theta <\pi$; $-\frac{\pi}{2}+\theta < \varphi^\ast,\varphi^{\ast\ast},\varphi^{\ast\ast\ast} <-\frac{\pi}{4}+\theta$ if $\frac{\pi}{2} \leq \theta <\frac{3\pi}{4}$; $0 < \varphi^\ast,\varphi^{\ast\ast},\varphi^{\ast\ast\ast} <-\frac{\pi}{4}+\theta$ if $\frac{\pi}{4} < \theta <\frac{\pi}{2}$; $ -\frac{\pi}{2} < \varphi^\ast,\varphi^{\ast\ast},\varphi^{\ast\ast\ast} < \frac{\pi}{2}+\theta$ if $-\pi <\theta \leq -\frac{3\pi}{4}$; $\frac{\pi}{4}+\theta < \varphi^\ast,\varphi^{\ast\ast},\varphi^{\ast\ast\ast} < \frac{\pi}{2}+\theta$ if $ -\frac{3\pi}{4} < \theta  \leq -\frac{\pi}{2}$; and $\frac{\pi}{4}+\theta < \varphi^\ast,\varphi^{\ast\ast},\varphi^{\ast\ast\ast} < 0$ if $-\frac{\pi}{2} < \theta <-\frac{\pi}{4}$.

Again, it is possible to derive simpler bounds that are useful near the Stokes lines $\arg \nu =\pm \frac{\pi}{2}$. It can be shown that for $\frac{\pi}{4} < \left|\arg \nu\right| \leq \frac{\pi}{2}$ and $N\geq 4$, we have
\begin{align*}
\left| {R_N^{\left( J \right)} \left( \nu  \right)} \right| \le \; & \frac{1}{2}\sqrt {\frac{e}{3}\left( {2N + \frac{{11}}{2}} \right)} \frac{1}{{3\pi }}6^{\frac{{2N + 1}}{3}} \left| {B_{2N} \left( 0 \right)} \right|\frac{{\sqrt 3 }}{2}\frac{{\Gamma \left( {\frac{{2N + 1}}{3}} \right)}}{{\left| \nu  \right|^{\frac{{2N + 1}}{3}} }} \\ & + \frac{1}{2}\sqrt {\frac{e}{3}\left( {2N + \frac{{19}}{2}} \right)} \frac{1}{{3\pi }}6^{\frac{{2N + 5}}{3}} \left| {B_{2N + 4} \left( 0 \right)} \right|\frac{{\sqrt 3 }}{2}\frac{{\Gamma \left( {\frac{{2N + 5}}{3}} \right)}}{{\left| \nu  \right|^{\frac{{2N + 5}}{3}} }},
\end{align*}
if $N\equiv 0 \mod 3$;
\begin{align*}
\left| {R_N^{\left( J \right)} \left( \nu  \right)} \right| \le \; & \frac{1}{2}\sqrt {\frac{e}{3}\left( {2N + \frac{{15}}{2}} \right)} \frac{1}{{3\pi }}6^{\frac{{2N + 3}}{3}} \left| {B_{2N + 2} \left( 0 \right)} \right|\frac{{\sqrt 3 }}{2}\frac{{\Gamma \left( {\frac{{2N + 3}}{3}} \right)}}{{\left| \nu  \right|^{\frac{{2N + 3}}{3}} }} \\ & + \frac{1}{2}\sqrt {\frac{e}{3}\left( {2N + \frac{{19}}{2}} \right)} \frac{1}{{3\pi }}6^{\frac{{2N + 5}}{3}} \left| {B_{2N + 4} \left( 0 \right)} \right|\frac{{\sqrt 3 }}{2}\frac{{\Gamma \left( {\frac{{2N + 5}}{3}} \right)}}{{\left| \nu  \right|^{\frac{{2N + 5}}{3}} }},
\end{align*}
for $N\equiv 1 \mod 3$; and
\begin{align*}
\left| {R_N^{\left( J \right)} \left( \nu  \right)} \right| \le \; & \frac{1}{2}\sqrt {\frac{e}{3}\left( {2N + \frac{{11}}{2}} \right)} \frac{1}{{3\pi }}6^{\frac{{2N + 1}}{3}} \left| {B_{2N} \left( 0 \right)} \right|\frac{{\sqrt 3 }}{2}\frac{{\Gamma \left( {\frac{{2N + 1}}{3}} \right)}}{{\left| \nu  \right|^{\frac{{2N + 1}}{3}} }} \\ & + \frac{1}{2}\sqrt {\frac{e}{3}\left( {2N + \frac{{15}}{2}} \right)} \frac{1}{{3\pi }}6^{\frac{{2N + 3}}{3}} \left| {B_{2N + 2} \left( 0 \right)} \right|\frac{{\sqrt 3 }}{2}\frac{{\Gamma \left( {\frac{{2N + 3}}{3}} \right)}}{{\left| \nu  \right|^{\frac{{2N + 3}}{3}} }},
\end{align*}
if $N\equiv 2 \mod 3$. These bounds are improvements over the similar results in the paper \cite{Nemes}.

Lastly, we consider the error bounds for the expansion of $Y_\nu\left(\nu\right)$. In the sector $\left|\arg \nu\right|<\frac{\pi}{2}$, $R_N^{\left( Y \right)} \left( \nu  \right)$ is given by \eqref{eq16}, for the wider range $\left|\arg \nu\right|<\pi$ we extend it via the expression \eqref{eq42}. The formula \eqref{eq16} can be simplified to
\begin{align*}
R_N^{\left( Y \right)} \left( \nu  \right) & = \frac{{\left( { - 1} \right)^{N + 1} }}{{2\pi \nu ^{\frac{{2N + 1}}{3}} }}\int_0^{ + \infty } {t^{\frac{{2N - 2}}{3}} e^{ - 2\pi t} \frac{{1 + \left( {t/\nu } \right)^{\frac{4}{3}} }}{{1 + \left( {t/\nu } \right)^2 }}iH_{it}^{\left( 1 \right)} \left( {it} \right)dt} ,\\
R_N^{\left( Y \right)} \left( \nu  \right) & = \frac{{\left( { - 1} \right)^N }}{{2\pi \nu ^{\frac{{2N + 3}}{3}} }}\int_0^{ + \infty } {t^{\frac{{2N}}{3}} e^{ - 2\pi t} \frac{{1 - \left( {t/\nu } \right)^{\frac{2}{3}} }}{{1 + \left( {t/\nu } \right)^2 }}iH_{it}^{\left( 1 \right)} \left( {it} \right)dt} ,\\
R_N^{\left( Y \right)} \left( \nu  \right) & = \frac{{\left( { - 1} \right)^{N + 1} }}{{2\pi \nu ^{\frac{{2N + 1}}{3}} }}\int_0^{ + \infty } {t^{\frac{{2N - 2}}{3}} e^{ - 2\pi t} \frac{{1 - \left( {t/\nu } \right)^{\frac{2}{3}} }}{{1 + \left( {t/\nu } \right)^2 }}iH_{it}^{\left( 1 \right)} \left( {it} \right)dt}
\end{align*}
according to whether $N\equiv 0 \mod 3$, $N\equiv 1 \mod 3$ or $N\equiv 2 \mod 3$. The error bounds suitable in the sector $\left|\arg \nu\right|<\frac{\pi}{2}$, were derived in the paper \cite{Nemes}. With the notations of the present paper, these estimates can be written as
\begin{align*}
\left| {R_N^{\left( Y \right)} \left( \nu  \right)} \right| \le \; & \left( \frac{2}{{3\pi }}6^{\frac{{2N + 1}}{3}} \left| {B_{2N} \left( 0 \right)} \right|\frac{3}{4}\frac{{\Gamma \left( {\frac{{2N + 1}}{3}} \right)}}{{\left| \nu  \right|^{\frac{{2N + 1}}{3}} }}\right. \\ & \left. + \frac{2}{{3\pi }}6^{\frac{{2N + 5}}{3}} \left| {B_{2N + 4} \left( 0 \right)} \right|\frac{3}{4}\frac{{\Gamma \left( {\frac{{2N + 5}}{3}} \right)}}{{\left| \nu  \right|^{\frac{{2N + 5}}{3}} }}\right) \begin{cases} \left|\csc\left(2\theta\right)\right| & \; \text{ if } \; \frac{\pi}{4} < \left|\theta\right| <\frac{\pi}{2}, \\ 1 & \; \text{ if } \; \left|\theta\right| \leq \frac{\pi}{4}, \end{cases}
\end{align*}
when $N\equiv 0 \mod 3$;
\[
\left| {R_N^{\left( Y \right)} \left( \nu  \right)} \right| \le \frac{2}{{3\pi }}6^{\frac{{2N + 3}}{3}} \left| {B_{2N + 2} \left( 0 \right)} \right|\frac{3}{4}\frac{{\Gamma \left( {\frac{{2N + 3}}{3}} \right)}}{{\left| \nu  \right|^{\frac{{2N + 3}}{3}} }}  \begin{cases} \left|\csc\left(2\theta\right)\right| & \; \text{ if } \; \frac{\pi}{4} < \left|\theta\right| <\frac{\pi}{2}, \\ 1 & \; \text{ if } \; \left|\theta\right| \leq \frac{\pi}{4}, \end{cases}
\]
if $N\equiv 1 \mod 3$;
\[
\left| {R_N^{\left( Y \right)} \left( \nu  \right)} \right| \le \frac{2}{{3\pi }}6^{\frac{{2N + 1}}{3}} \left| {B_{2N} \left( 0 \right)} \right|\frac{3}{4}\frac{{\Gamma \left( {\frac{{2N + 1}}{3}} \right)}}{{\left| \nu  \right|^{\frac{{2N + 1}}{3}} }} \begin{cases} \left|\csc\left(2\theta\right)\right| & \; \text{ if } \; \frac{\pi}{4} < \left|\theta\right| <\frac{\pi}{2}, \\ 1 & \; \text{ if } \; \left|\theta\right| \leq \frac{\pi}{4}, \end{cases}
\]
when $N\equiv 2 \mod 3$, respectively. The reference \cite{Nemes} also contains sharp bounds for the case $\nu>0$. The corresponding estimates for the sectors $\frac{\pi}{4}<\left|\arg \nu\right|<\pi$ are found to be
\begin{align*}
\left| {R_N^{\left( Y \right)} \left( \nu  \right)} \right| \le  \; & \frac{{\left| {\csc \left( {2\left( {\theta  - \varphi^\ast} \right)} \right)} \right|}}{{\cos ^{\frac{{2N + 1}}{3}} \varphi^\ast}} \frac{2}{{3\pi }}6^{\frac{{2N + 1}}{3}} \left| {B_{2N} \left( 0 \right)} \right|\frac{3}{4}\frac{{\Gamma \left( {\frac{{2N + 1}}{3}} \right)}}{{\left| \nu  \right|^{\frac{{2N + 1}}{3}} }} \\ & + \frac{{\left| {\csc \left( {2\left( {\theta  - \varphi^{\ast\ast\ast} } \right)} \right)} \right|}}{{\cos ^{\frac{{2N + 5}}{3}} \varphi^{\ast\ast\ast}}}\frac{2}{{3\pi }}6^{\frac{{2N + 5}}{3}} \left| {B_{2N + 4} \left( 0 \right)} \right|\frac{3}{4}\frac{{\Gamma \left( {\frac{{2N + 5}}{3}} \right)}}{{\left| \nu  \right|^{\frac{{2N + 5}}{3}} }},
\end{align*}
when $N\equiv 0 \mod 3$;
\[
\left| {R_N^{\left( Y \right)} \left( \nu  \right)} \right| \le \frac{{\left| {\csc \left( {2\left( {\theta  - \varphi^{\ast\ast} } \right)} \right)} \right|}}{{\cos ^{\frac{{2N + 3}}{3}} \varphi^{\ast\ast}}}\frac{2}{{3\pi }}6^{\frac{{2N + 3}}{3}} \left| {B_{2N + 2} \left( 0 \right)} \right|\frac{3}{4}\frac{{\Gamma \left( {\frac{{2N + 3}}{3}} \right)}}{{\left| \nu  \right|^{\frac{{2N + 3}}{3}} }},
\]
if $N\equiv 1 \mod 3$;
\[
\left| {R_N^{\left( Y \right)} \left( \nu  \right)} \right| \le \frac{{\left| {\csc \left( {2\left( {\theta  - \varphi^\ast} \right)} \right)} \right|}}{{\cos ^{\frac{{2N + 1}}{3}} \varphi^\ast}} \frac{2}{{3\pi }}6^{\frac{{2N + 1}}{3}} \left| {B_{2N} \left( 0 \right)} \right|\frac{3}{4}\frac{{\Gamma \left( {\frac{{2N + 1}}{3}} \right)}}{{\left| \nu  \right|^{\frac{{2N + 1}}{3}} }},
\]
when $N\equiv 2 \mod 3$, respectively. Here $\varphi^\ast$, $\varphi^{\ast\ast}$ and $\varphi^{\ast\ast\ast}$ are the same as in the bounds for $R_N^{\left( J \right)} \left( \nu  \right)$. In deriving the second and the third bound, we applied the inequality
\begin{equation}\label{eq69}
\left| {\frac{{1 - r^{\frac{1}{3}}e^{-\frac{2}{3}i\vartheta} }}{{1 + r e^{-2i\vartheta} }}} \right| \le \begin{cases} \left|\csc\left(2\vartheta\right)\right| & \; \text{ if } \; \frac{\pi}{4} < \left|\vartheta\right| <\frac{\pi}{2}, \\ 1 & \; \text{ if } \; \left|\vartheta\right| \leq \frac{\pi}{4}, \end{cases}
\end{equation}
for $r>0$, whose proof is given in \cite[Appendix B]{Nemes}.

Again, it is possible to derive simpler bounds that are effective near the Stokes lines $\arg \nu =\pm \frac{\pi}{2}$. It can be shown that for $\frac{\pi}{4} < \left|\arg \nu\right| \leq \frac{\pi}{2}$ and $N\geq 4$, one has
\begin{align*}
\left| {R_N^{\left( Y \right)} \left( \nu  \right)} \right| \le \; & \frac{1}{2}\sqrt {\frac{e}{3}\left( {2N + \frac{{11}}{2}} \right)}  \frac{2}{{3\pi }}6^{\frac{{2N + 1}}{3}} \left| {B_{2N} \left( 0 \right)} \right|\frac{3}{4}\frac{{\Gamma \left( {\frac{{2N + 1}}{3}} \right)}}{{\left| \nu  \right|^{\frac{{2N + 1}}{3}} }} \\ & + \frac{1}{2}\sqrt {\frac{e}{3}\left( {2N + \frac{{19}}{2}} \right)} \frac{2}{{3\pi }}6^{\frac{{2N + 5}}{3}} \left| {B_{2N + 4} \left( 0 \right)} \right|\frac{3}{4}\frac{{\Gamma \left( {\frac{{2N + 5}}{3}} \right)}}{{\left| \nu  \right|^{\frac{{2N + 5}}{3}} }},
\end{align*}
when $N\equiv 0 \mod 3$;
\[
\left| {R_N^{\left( Y \right)} \left( \nu  \right)} \right| \le \frac{1}{2}\sqrt {\frac{e}{3}\left( {2N + \frac{{15}}{2}} \right)} \frac{2}{{3\pi }}6^{\frac{{2N + 3}}{3}} \left| {B_{2N + 2} \left( 0 \right)} \right|\frac{3}{4}\frac{{\Gamma \left( {\frac{{2N + 3}}{3}} \right)}}{{\left| \nu  \right|^{\frac{{2N + 3}}{3}} }},
\]
if $N\equiv 1 \mod 3$;
\[
\left| {R_N^{\left( Y \right)} \left( \nu  \right)} \right| \le \frac{1}{2}\sqrt {\frac{e}{3}\left( {2N + \frac{{11}}{2}} \right)}  \frac{2}{{3\pi }}6^{\frac{{2N + 1}}{3}} \left| {B_{2N} \left( 0 \right)} \right|\frac{3}{4}\frac{{\Gamma \left( {\frac{{2N + 1}}{3}} \right)}}{{\left| \nu  \right|^{\frac{{2N + 1}}{3}} }},
\]
when $N\equiv 2 \mod 3$, respectively. These bounds are improvements over the similar results in the paper \cite{Nemes}.

\subsection{Error bounds for the asymptotic series of $J'_\nu \left( \nu  \right)$ and $Y'_\nu \left( \nu  \right)$} First, we consider the estimation of $R_N^{\left( J'\right)} \left( \nu  \right)$ in the sector $\left|\arg \nu \right|<\pi$. The representation \eqref{eq65} is valid only if $\left|\arg \nu \right|<\frac{\pi}{2}$, for the larger sector, we use \eqref{eq64} as the definition of $R_N^{\left( J'\right)} \left( \nu  \right)$. Throughout this subsection, we assume that $N\geq 1$. Algebraic manipulation of \eqref{eq65} shows that
\[
R_N^{\left( {J'} \right)} \left( \nu  \right) = \frac{{\left( { - 1} \right)^{N + 1} }}{{2\sqrt 3 \pi \nu ^{\frac{{2N + 2}}{3}} }}\int_0^{ + \infty } {t^{\frac{{2N - 1}}{3}} e^{ - 2\pi t} \frac{{1 + \left( {t/\nu } \right)^{\frac{2}{3}} }}{{1 + \left( {t/\nu } \right)^2 }}H_{it}^{\left( 1 \right)\prime  } \left( {it} \right)dt} ,
\]
\[
R_N^{\left( {J'} \right)} \left( \nu  \right) = \frac{{\left( { - 1} \right)^N }}{{2\sqrt 3 \pi \nu ^{\frac{{2N + 2}}{3}} }}\int_0^{ + \infty } {t^{\frac{{2N - 1}}{3}} e^{ - 2\pi t} \frac{{1 - \left( {t/\nu } \right)^{\frac{4}{3}} }}{{1 + \left( {t/\nu } \right)^2 }}H_{it}^{\left( 1 \right)\prime  } \left( {it} \right)dt} ,
\]
\[
R_N^{\left( {J'} \right)} \left( \nu  \right) = \frac{{\left( { - 1} \right)^N }}{{2\sqrt 3 \pi \nu ^{\frac{{2N + 4}}{3}} }}\int_0^{ + \infty } {t^{\frac{{2N + 1}}{3}} e^{ - 2\pi t} \frac{{1 + \left( {t/\nu } \right)^{\frac{2}{3}} }}{{1 + \left( {t/\nu } \right)^2 }}H_{it}^{\left( 1 \right)\prime  } \left( {it} \right)dt} 
\]
according to whether $N\equiv 0 \mod 3$, $N\equiv 1 \mod 3$ or $N\equiv 2 \mod 3$, respectively. Employing the inequalities
\[
\left| {\frac{{1 - r^{\frac{1}{3}} e^{ - \frac{4}{3}i\vartheta } }}{{1 + re^{ - 2i\vartheta } }}} \right| \le 1 \; \text{ for } \; r>0 \; \text{ and } \; \left| \vartheta  \right| \le \frac{\pi }{4},
\]
and \eqref{eq66}, we obtain the error bounds
\begin{align*}
\left| {R_N^{\left( {J'} \right)} \left( \nu  \right)} \right| \le \; & \left( {\frac{1}{{3\pi }}6^{\frac{{2N + 2}}{3}} \left|D_{2N + 1} \left( 0 \right)\right|\frac{{\sqrt 3 }}{2}\frac{{\Gamma \left( {\frac{{2N + 2}}{3}} \right)}}{{\left| \nu  \right|^{\frac{{2N + 2}}{3}} }}}\right. \\ & \left.{+ \frac{1}{{3\pi }}6^{\frac{{2N + 4}}{3}} \left|D_{2N + 3} \left( 0 \right)\right|\frac{{\sqrt 3 }}{2}\frac{{\Gamma \left( {\frac{{2N + 4}}{3}} \right)}}{{\left| \nu  \right|^{\frac{{2N + 4}}{3}} }}} \right) \begin{cases} \left|\csc\left(2\theta\right)\right| & \; \text{ if } \; \frac{\pi}{4} < \left|\theta\right| <\frac{\pi}{2}, \\ 1 & \; \text{ if } \; \left|\theta\right| \leq \frac{\pi}{4}, \end{cases}
\end{align*}
for $N \equiv 0 \mod 3$;
\begin{align*}
\left| {R_N^{\left( {J'} \right)} \left( \nu  \right)} \right| \le \; & \frac{1}{{3\pi }}6^{\frac{{2N + 2}}{3}} \left|D_{2N + 1} \left( 0 \right)\right|\frac{{\sqrt 3 }}{2}\frac{{\Gamma \left( {\frac{{2N + 2}}{3}} \right)}}{{\left|\nu\right|^{\frac{{2N + 2}}{3}} }} \begin{cases} \left|\csc\left(2\theta\right)\right| & \; \text{ if } \; \frac{\pi}{4} < \left|\theta\right| <\frac{\pi}{2}, \\ 1 & \; \text{ if } \; \left|\theta\right| \leq \frac{\pi}{4}, \end{cases} \\ & + \frac{1}{{3\pi }}6^{\frac{{2N + 6}}{3}} \left|D_{2N + 5} \left( 0 \right)\right|\frac{{\sqrt 3 }}{2}\frac{{\Gamma \left( {\frac{{2N + 6}}{3}} \right)}}{{\left|\nu\right|^{\frac{{2N + 6}}{3}} }} \begin{cases} \left|\csc\left(2\theta\right)\right| & \; \text{ if } \; \frac{\pi}{4} < \left|\theta\right| <\frac{\pi}{2}, \\ 0 & \; \text{ if } \; \left|\theta\right| \leq \frac{\pi}{4}, \end{cases}
\end{align*}
if $N \equiv 1 \mod 3$; and
\begin{align*}
\left| {R_N^{\left( {J'} \right)} \left( \nu  \right)} \right| \le \; & \left( {\frac{1}{{3\pi }}6^{\frac{{2N + 4}}{3}} \left|D_{2N + 3} \left( 0 \right)\right|\frac{{\sqrt 3 }}{2}\frac{{\Gamma \left( {\frac{{2N + 4}}{3}} \right)}}{{\left|\nu\right| ^{\frac{{2N + 4}}{3}} }} }\right. \\ & \left.{+ \frac{1}{{3\pi }}6^{\frac{{2N + 6}}{3}} \left|D_{2N + 5} \left( 0 \right)\right|\frac{{\sqrt 3 }}{2}\frac{{\Gamma \left( {\frac{{2N + 6}}{3}} \right)}}{{\left|\nu\right| ^{\frac{{2N + 6}}{3}} }}} \right) \begin{cases} \left|\csc\left(2\theta\right)\right| & \; \text{ if } \; \frac{\pi}{4} < \left|\theta\right| <\frac{\pi}{2}, \\ 1 & \; \text{ if } \; \left|\theta\right| \leq \frac{\pi}{4}, \end{cases}
\end{align*}
for $N \equiv 2 \mod 3$.

When $\nu$ is real and positive, we can obtain more precise estimates. Indeed as $0<\frac{1}{1+\left(t/\nu\right)^2}<1$ and $H_{it}^{\left( 1 \right)\prime  } \left( {it} \right)>0$ for $t,\nu>0$, the main value theorem of integration shows that
\[
\left( { - 1} \right)^{N + 1} R_N^{\left( {J'} \right)} \left( \nu  \right) = \frac{1}{{3\pi }}6^{\frac{{2N + 2}}{3}} \left| {D_{2N + 1} \left( 0 \right)} \right|\frac{{\sqrt 3 }}{2}\frac{{\Gamma \left( {\frac{{2N + 2}}{3}} \right)}}{{\nu ^{\frac{{2N + 2}}{3}} }}\Theta _1  + \frac{1}{{3\pi }}6^{\frac{{2N + 4}}{3}} \left| {D_{2N + 3} \left( 0 \right)} \right|\frac{{\sqrt 3 }}{2}\frac{{\Gamma \left( {\frac{{2N + 4}}{3}} \right)}}{{\nu ^{\frac{{2N + 4}}{3}} }}\Theta _2 ,
\]
when $N \equiv 0 \mod 3$;
\[
\left( { - 1} \right)^N R_N^{\left( {J'} \right)} \left( \nu  \right) = \frac{1}{{3\pi }}6^{\frac{{2N + 2}}{3}} \left| {D_{2N + 1} \left( 0 \right)} \right|\frac{{\sqrt 3 }}{2}\frac{{\Gamma \left( {\frac{{2N + 2}}{3}} \right)}}{{\nu ^{\frac{{2N + 2}}{3}} }}\Theta _1  - \frac{1}{{3\pi }}6^{\frac{{2N + 6}}{3}} \left| {D_{2N + 5} \left( 0 \right)} \right|\frac{{\sqrt 3 }}{2}\frac{{\Gamma \left( {\frac{{2N + 6}}{3}} \right)}}{{\nu ^{\frac{{2N + 6}}{3}} }}\Theta _3 ,
\]
if $N \equiv 1 \mod 3$;
\[
\left( { - 1} \right)^N R_N^{\left( {J'} \right)} \left( \nu  \right) = \frac{1}{{3\pi }}6^{\frac{{2N + 4}}{3}} \left| {D_{2N + 3} \left( 0 \right)} \right|\frac{{\sqrt 3 }}{2}\frac{{\Gamma \left( {\frac{{2N + 4}}{3}} \right)}}{{\nu ^{\frac{{2N + 4}}{3}} }}\Theta _2  + \frac{1}{{3\pi }}6^{\frac{{2N + 6}}{3}} \left| {D_{2N + 5} \left( 0 \right)} \right|\frac{{\sqrt 3 }}{2}\frac{{\Gamma \left( {\frac{{2N + 6}}{3}} \right)}}{{\nu ^{\frac{{2N + 6}}{3}} }}\Theta _3 ,
\]
when $N \equiv 2 \mod 3$, respectively. Here $0<\Theta_i<1$ ($i=1,2,3$) is an appropriate number depending on $\nu$ and $N$. We remark that Watson \cite[p. 260]{Watson} showed that for $\nu>0$,
\[
J'_\nu  \left( \nu  \right) < \frac{1}{{3\pi }}6^{\frac{2}{3}} D_1 \left( 0 \right)\frac{{\sqrt 3 }}{2}\frac{{\Gamma \left( {\frac{2}{3}} \right)}}{{\nu ^{\frac{2}{3}} }} = \frac{{3^{\frac{1}{6}} \Gamma \left( {\frac{2}{3}} \right)}}{{2^{\frac{1}{3}} \pi \nu ^{\frac{2}{3}} }},
\]
or, in other words, $R_1^{\left( {J'} \right)} \left( \nu  \right) < 0$ for positive $\nu$.

To derive the corresponding estimates for the sectors $\frac{\pi}{4}<\left|\arg \nu\right|<\pi$, we use the inequality \eqref{eq66}, Meijer's lemma \cite[p. 956]{Meijer} and a similar argument to that used in estimating $R_N^{\left( {H'} \right)} \left( \nu  \right)$. The final results are
\begin{align*}
\left| {R_N^{\left( {J'} \right)} \left( \nu  \right)} \right| \le \; & \frac{{\left| {\csc \left( {2\left( {\theta  - \varphi^{\ast} } \right)} \right)} \right|}}{{\cos ^{\frac{{2N + 2}}{3}} \varphi^{\ast} }}\frac{1}{{3\pi }}6^{\frac{{2N + 2}}{3}} \left| {D_{2N + 1} \left( 0 \right)} \right|\frac{{\sqrt 3 }}{2}\frac{{\Gamma \left( {\frac{{2N + 2}}{3}} \right)}}{{\left| \nu  \right|^{\frac{{2N + 2}}{3}} }} \\ &+ \frac{{\left| {\csc \left( {2\left( {\theta  - \varphi^{\ast\ast} } \right)} \right)} \right|}}{{\cos ^{\frac{{2N + 4}}{3}} \varphi^{\ast\ast} }}\frac{1}{{3\pi }}6^{\frac{{2N + 4}}{3}} \left| {D_{2N + 3} \left( 0 \right)} \right|\frac{{\sqrt 3 }}{2}\frac{{\Gamma \left( {\frac{{2N + 4}}{3}} \right)}}{{\left| \nu  \right|^{\frac{{2N + 4}}{3}} }},
\end{align*}
for $N \equiv 0 \mod 3$;
\begin{align*}
\left| {R_N^{\left( {J'} \right)} \left( \nu  \right)} \right| \le \; & \frac{{\left| {\csc \left( {2\left( {\theta  - \varphi^{\ast} } \right)} \right)} \right|}}{{\cos ^{\frac{{2N + 2}}{3}} \varphi^{\ast} }}\frac{1}{{3\pi }}6^{\frac{{2N + 2}}{3}} \left| {D_{2N + 1} \left( 0 \right)} \right|\frac{{\sqrt 3 }}{2}\frac{{\Gamma \left( {\frac{{2N + 2}}{3}} \right)}}{{\left| \nu  \right|^{\frac{{2N + 2}}{3}} }} \\ &+ \frac{{\left| {\csc \left( {2\left( {\theta  - \varphi^{\ast\ast\ast}} \right)} \right)} \right|}}{{\cos ^{\frac{{2N + 6}}{3}} \varphi^{\ast\ast\ast} }}\frac{1}{{3\pi }}6^{\frac{{2N + 6}}{3}} \left| {D_{2N + 5} \left( 0 \right)} \right|\frac{{\sqrt 3 }}{2}\frac{{\Gamma \left( {\frac{{2N + 6}}{3}} \right)}}{{\left| \nu  \right|^{\frac{{2N + 6}}{3}} }},
\end{align*}
if $N \equiv 1 \mod 3$; and
\begin{align*}
\left| {R_N^{\left( {J'} \right)} \left( \nu  \right)} \right| \le \; & \frac{{\left| {\csc \left( {2\left( {\theta  - \varphi^{\ast\ast} } \right)} \right)} \right|}}{{\cos ^{\frac{{2N + 4}}{3}} \varphi^{\ast\ast} }}\frac{1}{{3\pi }}6^{\frac{{2N + 4}}{3}} \left| {D_{2N + 3} \left( 0 \right)} \right|\frac{{\sqrt 3 }}{2}\frac{{\Gamma \left( {\frac{{2N + 4}}{3}} \right)}}{{\left| \nu  \right|^{\frac{{2N + 4}}{3}} }} \\ &+ \frac{{\left| {\csc \left( {2\left( {\theta  - \varphi^{\ast\ast\ast}} \right)} \right)} \right|}}{{\cos ^{\frac{{2N + 6}}{3}} \varphi^{\ast\ast\ast}}}\frac{1}{{3\pi }}6^{\frac{{2N + 6}}{3}} \left| {D_{2N + 5} \left( 0 \right)} \right|\frac{{\sqrt 3 }}{2}\frac{{\Gamma \left( {\frac{{2N + 6}}{3}} \right)}}{{\left| \nu  \right|^{\frac{{2N + 6}}{3}} }},
\end{align*}
for $N \equiv 2 \mod 3$, where $\varphi^\ast$, $\varphi^{\ast \ast}$ and $\varphi^{\ast \ast \ast}$ are the unique solutions of the implicit equations
\[
\left( {2N + 8 } \right)\cos \left( {3\varphi^\ast - 2\theta} \right) = \left( {2N - 6 } \right)\cos \left( {\varphi^\ast- 2\theta } \right),
\]
\[
\left( {2N + 10 } \right)\cos \left( {3\varphi^{\ast\ast} - 2\theta} \right) = \left( {2N - 2 } \right)\cos \left( {\varphi^{\ast\ast}- 2\theta } \right)
\]
and
\[
\left( {2N + 12 } \right)\cos \left( {3\varphi^{\ast\ast\ast} - 2\theta} \right) = 2N \cos \left( {\varphi^{\ast\ast\ast}- 2\theta } \right),
\]
that satisfy $-\frac{\pi}{2}+\theta < \varphi^\ast,\varphi^{\ast\ast},\varphi^{\ast\ast\ast} <\frac{\pi}{2}$ if $\frac{3\pi}{4} \leq \theta <\pi$; $-\frac{\pi}{2}+\theta < \varphi^\ast,\varphi^{\ast\ast},\varphi^{\ast\ast\ast} <-\frac{\pi}{4}+\theta$ if $\frac{\pi}{2} \leq \theta <\frac{3\pi}{4}$; $0 < \varphi^\ast,\varphi^{\ast\ast},\varphi^{\ast\ast\ast} <-\frac{\pi}{4}+\theta$ if $\frac{\pi}{4} < \theta <\frac{\pi}{2}$; $ -\frac{\pi}{2} < \varphi^\ast,\varphi^{\ast\ast},\varphi^{\ast\ast\ast} < \frac{\pi}{2}+\theta$ if $-\pi <\theta \leq -\frac{3\pi}{4}$; $\frac{\pi}{4}+\theta < \varphi^\ast,\varphi^{\ast\ast},\varphi^{\ast\ast\ast} < \frac{\pi}{2}+\theta$ if $ -\frac{3\pi}{4} < \theta  \leq -\frac{\pi}{2}$; and $\frac{\pi}{4}+\theta < \varphi^\ast,\varphi^{\ast\ast},\varphi^{\ast\ast\ast} < 0$ if $-\frac{\pi}{2} < \theta <-\frac{\pi}{4}$.

For the simpler bounds that are useful near the Stokes lines $\arg \nu =\pm \frac{\pi}{2}$, it is readily established that for $\frac{\pi}{4} < \left|\arg \nu\right| \leq \frac{\pi}{2}$ and $N\geq 4$, we have
\begin{align*}
\left| {R_N^{\left( {J'} \right)} \left( \nu  \right)} \right| \le \; & \frac{1}{2}\sqrt {\frac{e}{3}\left( {2N + \frac{{13}}{2}} \right)} \frac{1}{{3\pi }}6^{\frac{{2N + 2}}{3}} \left| {D_{2N + 1} \left( 0 \right)} \right|\frac{{\sqrt 3 }}{2}\frac{{\Gamma \left( {\frac{{2N + 2}}{3}} \right)}}{{\left| \nu  \right|^{\frac{{2N + 2}}{3}} }} \\ & + \frac{1}{2}\sqrt {\frac{e}{3}\left( {2N + \frac{{17}}{2}} \right)} \frac{1}{{3\pi }}6^{\frac{{2N + 4}}{3}} \left| {D_{2N + 3} \left( 0 \right)} \right|\frac{{\sqrt 3 }}{2}\frac{{\Gamma \left( {\frac{{2N + 4}}{3}} \right)}}{{\left| \nu  \right|^{\frac{{2N + 4}}{3}} }},
\end{align*}
when $N \equiv 0 \mod 3$;
\begin{align*}
\left| {R_N^{\left( {J'} \right)} \left( \nu  \right)} \right| \le \; & \frac{1}{2}\sqrt {\frac{e}{3}\left( {2N + \frac{{13}}{2}} \right)} \frac{1}{{3\pi }}6^{\frac{{2N + 2}}{3}} \left| {D_{2N + 1} \left( 0 \right)} \right|\frac{{\sqrt 3 }}{2}\frac{{\Gamma \left( {\frac{{2N + 2}}{3}} \right)}}{{\left| \nu  \right|^{\frac{{2N + 2}}{3}} }} \\ & + \frac{1}{2}\sqrt {\frac{e}{3}\left( {2N + \frac{{21}}{2}} \right)} \frac{1}{{3\pi }}6^{\frac{{2N + 6}}{3}} \left| {D_{2N + 5} \left( 0 \right)} \right|\frac{{\sqrt 3 }}{2}\frac{{\Gamma \left( {\frac{{2N + 6}}{3}} \right)}}{{\left| \nu  \right|^{\frac{{2N + 6}}{3}} }},
\end{align*}
for $N \equiv 1 \mod 3$;
\begin{align*}
\left| {R_N^{\left( {J'} \right)} \left( \nu  \right)} \right| \le \; & \frac{1}{2}\sqrt {\frac{e}{3}\left( {2N + \frac{{17}}{2}} \right)} \frac{1}{{3\pi }}6^{\frac{{2N + 4}}{3}} \left| {D_{2N + 3} \left( 0 \right)} \right|\frac{{\sqrt 3 }}{2}\frac{{\Gamma \left( {\frac{{2N + 4}}{3}} \right)}}{{\left| \nu  \right|^{\frac{{2N + 4}}{3}} }} \\ & + \frac{1}{2}\sqrt {\frac{e}{3}\left( {2N + \frac{{21}}{2}} \right)} \frac{1}{{3\pi }}6^{\frac{{2N + 6}}{3}} \left| {D_{2N + 5} \left( 0 \right)} \right|\frac{{\sqrt 3 }}{2}\frac{{\Gamma \left( {\frac{{2N + 6}}{3}} \right)}}{{\left| \nu  \right|^{\frac{{2N + 6}}{3}} }},
\end{align*}
when $N \equiv 2 \mod 3$.

Finally, we consider the error estimates for the expansion of $Y'_\nu\left(\nu\right)$. In the sector $\left|\arg \nu\right|<\frac{\pi}{2}$, $R_N^{\left( Y' \right)} \left( \nu  \right)$ is given by the formula \eqref{eq67}, for the wider range $\left|\arg \nu\right|<\pi$ we define it by the expression \eqref{eq68}. Simple algebraic manipulation shows that the formula \eqref{eq67} can be written as
\begin{align*}
R_N^{\left( {Y'} \right)} \left( \nu  \right) & = \frac{{\left( { - 1} \right)^{N + 1} }}{{2\pi \nu ^{\frac{{2N + 2}}{3}} }}\int_0^{ + \infty } {t^{\frac{{2N - 1}}{3}} e^{ - 2\pi t} \frac{{1 - \left( {t/\nu } \right)^{\frac{2}{3}} }}{{1 + \left( {t/\nu } \right)^2 }}H_{it}^{\left( 1 \right)\prime  } \left( {it} \right)dt} ,\\
R_N^{\left( {Y'} \right)} \left( \nu  \right) & = \frac{{\left( { - 1} \right)^{N + 1} }}{{2\pi \nu ^{\frac{{2N + 2}}{3}} }}\int_0^{ + \infty } {t^{\frac{{2N - 1}}{3}} e^{ - 2\pi t} \frac{{1 + \left( {t/\nu } \right)^{\frac{4}{3}} }}{{1 + \left( {t/\nu } \right)^2 }}H_{it}^{\left( 1 \right)^\prime  } \left( {it} \right)dt} ,\\
R_N^{\left( {Y'} \right)} \left( \nu  \right) & = \frac{{\left( { - 1} \right)^N }}{{2\pi \nu ^{\frac{{2N + 4}}{3}} }}\int_0^{ + \infty } {t^{\frac{{2N + 1}}{3}} e^{ - 2\pi t} \frac{{1 - \left( {t/\nu } \right)^{\frac{2}{3}} }}{{1 + \left( {t/\nu } \right)^2 }}H_{it}^{\left( 1 \right)^\prime  } \left( {it} \right)dt} 
\end{align*}
according to whether $N\equiv 0 \mod 3$, $N\equiv 1 \mod 3$ or $N\equiv 2 \mod 3$. Applying the inequality \eqref{eq66} in the second and the
inequality \eqref{eq69} in the first and the third case, we deduce
\[
\left| {R_N^{\left( {Y'} \right)} \left( \nu  \right)} \right| \le \frac{2}{{3\pi }}6^{\frac{{2N + 2}}{3}} \left| {D_{2N + 1} \left( 0 \right)} \right|\frac{3}{4}\frac{{\Gamma \left( {\frac{{2N + 2}}{3}} \right)}}{{\left| \nu  \right|^{\frac{{2N + 2}}{3}} }} \begin{cases} \left|\csc\left(2\theta\right)\right| & \; \text{ if } \; \frac{\pi}{4} < \left|\theta\right| <\frac{\pi}{2}, \\ 1 & \; \text{ if } \; \left|\theta\right| \leq \frac{\pi}{4}, \end{cases}
\]
when $N \equiv 0 \mod 3$;
\begin{align*}
\left| {R_N^{\left( {Y'} \right)} \left( \nu  \right)} \right| \le \; & \left( {\frac{2}{{3\pi }}6^{\frac{{2N + 2}}{3}} \left| {D_{2N + 1} \left( 0 \right)} \right|\frac{3}{4}\frac{{\Gamma \left( {\frac{{2N + 2}}{3}} \right)}}{{\left| \nu  \right|^{\frac{{2N + 2}}{3}} }} }\right.\\ & \left.{+ \frac{2}{{3\pi }}6^{\frac{{2N + 6}}{3}} \left| {D_{2N + 5} \left( 0 \right)} \right|\frac{3}{4}\frac{{\Gamma \left( {\frac{{2N + 6}}{3}} \right)}}{{\left| \nu  \right|^{\frac{{2N + 6}}{3}} }}} \right) \begin{cases} \left|\csc\left(2\theta\right)\right| & \; \text{ if } \; \frac{\pi}{4} < \left|\theta\right| <\frac{\pi}{2}, \\ 1 & \; \text{ if } \; \left|\theta\right| \leq \frac{\pi}{4}, \end{cases}
\end{align*}
for $N \equiv 1 \mod 3$;
\[
\left| {R_N^{\left( {Y'} \right)} \left( \nu  \right)} \right| \le \frac{2}{{3\pi }}6^{\frac{{2N + 4}}{3}} \left| {D_{2N + 3} \left( 0 \right)} \right|\frac{3}{4}\frac{{\Gamma \left( {\frac{{2N + 4}}{3}} \right)}}{{\left| \nu  \right|^{\frac{{2N + 4}}{3}} }} \begin{cases} \left|\csc\left(2\theta\right)\right| & \; \text{ if } \; \frac{\pi}{4} < \left|\theta\right| <\frac{\pi}{2}, \\ 1 & \; \text{ if } \; \left|\theta\right| \leq \frac{\pi}{4}, \end{cases}
\]
when $N \equiv 2 \mod 3$, respectively. When $\nu$ is real and positive, we have the following better estimates:
\[
\left( { - 1} \right)^{N + 1} R_N^{\left( {Y'} \right)} \left( \nu  \right) = \frac{2}{{3\pi }}6^{\frac{{2N + 2}}{3}} \left| {D_{2N + 1} \left( 0 \right)} \right|\frac{3}{4}\frac{{\Gamma \left( {\frac{{2N + 2}}{3}} \right)}}{{ \nu^{\frac{{2N + 2}}{3}} }}\Xi _1  - \frac{2}{{3\pi }}6^{\frac{{2N + 4}}{3}} \left| {D_{2N + 3} \left( 0 \right)} \right|\frac{3}{4}\frac{{\Gamma \left( {\frac{{2N + 4}}{3}} \right)}}{{ \nu^{\frac{{2N + 4}}{3}} }}\Xi _2 ,
\]
if $N \equiv 0 \mod 3$;
\[
\left( { - 1} \right)^{N + 1} R_N^{\left( {Y'} \right)} \left( \nu  \right) = \frac{2}{{3\pi }}6^{\frac{{2N + 2}}{3}} \left| {D_{2N + 1} \left( 0 \right)} \right|\frac{3}{4}\frac{{\Gamma \left( {\frac{{2N + 2}}{3}} \right)}}{{\nu^{\frac{{2N + 2}}{3}} }}\Xi _1  + \frac{2}{{3\pi }}6^{\frac{{2N + 6}}{3}} \left| {D_{2N + 5} \left( 0 \right)} \right|\frac{3}{4}\frac{{\Gamma \left( {\frac{{2N + 6}}{3}} \right)}}{{ \nu^{\frac{{2N + 6}}{3}} }}\Xi _3 ,
\]
when $N \equiv 1 \mod 3$;
\[
\left( { - 1} \right)^N R_N^{\left( {Y'} \right)} \left( \nu  \right) = \frac{2}{{3\pi }}6^{\frac{{2N + 4}}{3}} \left| {D_{2N + 3} \left( 0 \right)} \right|\frac{3}{4}\frac{{\Gamma \left( {\frac{{2N + 4}}{3}} \right)}}{{\nu^{\frac{{2N + 4}}{3}} }}\Xi _2  - \frac{2}{{3\pi }}6^{\frac{{2N + 6}}{3}} \left| {D_{2N + 5} \left( 0 \right)} \right|\frac{3}{4}\frac{{\Gamma \left( {\frac{{2N + 6}}{3}} \right)}}{{\nu^{\frac{{2N + 6}}{3}} }}\Xi _3 ,
\]
if $N \equiv 2 \mod 3$. Here $0<\Xi_i<1$ ($i=1,2,3$) is a suitable number depending on $\nu$ and $N$. The second estimate shows that $R_1^{\left( {Y'} \right)} \left( \nu  \right)$ is positive if $\nu$ is positive, whence we obtain the following Watson-type estimate:
\[
Y'_\nu  \left( \nu  \right) > \frac{2}{{3\pi }}6^{\frac{2}{3}} D_1 \left( 0 \right)\frac{3}{4}\frac{{\Gamma \left( {\frac{2}{3}} \right)}}{{\nu ^{\frac{2}{3}} }} = \frac{{3^{\frac{2}{3}} \Gamma \left( {\frac{2}{3}} \right)}}{{2^{\frac{1}{3}} \pi \nu ^{\frac{2}{3}} }},
\]
for $\nu >0$. Therefore, the leading order asymptotics of $Y'_\nu  \left( \nu  \right)$ is always a lower estimate to the actual value if $\nu$ is positive.

The corresponding estimates for the sectors $\frac{\pi}{4}<\left|\arg \nu\right|<\pi$ are found to be
\[
\left| {R_N^{\left( {Y'} \right)} \left( \nu  \right)} \right| \le \frac{{\left| {\csc \left( {2\left( {\theta  - \varphi^{\ast} } \right)} \right)} \right|}}{{\cos ^{\frac{{2N + 2}}{3}} \varphi^{\ast} }}\frac{2}{{3\pi }}6^{\frac{{2N + 2}}{3}} \left| {D_{2N + 1} \left( 0 \right)} \right|\frac{3}{4}\frac{{\Gamma \left( {\frac{{2N + 2}}{3}} \right)}}{{\left| \nu  \right|^{\frac{{2N + 2}}{3}} }},
\]
for $N \equiv 0 \mod 3$;
\begin{align*}
\left| {R_N^{\left( {Y'} \right)} \left( \nu  \right)} \right| \le \; & \frac{{\left| {\csc \left( {2\left( {\theta  - \varphi^{\ast} } \right)} \right)} \right|}}{{\cos ^{\frac{{2N + 2}}{3}} \varphi^{\ast} }}\frac{2}{{3\pi }}6^{\frac{{2N + 2}}{3}} \left| {D_{2N + 1} \left( 0 \right)} \right|\frac{3}{4}\frac{{\Gamma \left( {\frac{{2N + 2}}{3}} \right)}}{{\left| \nu  \right|^{\frac{{2N + 2}}{3}} }} \\ & + \frac{{\left| {\csc \left( {2\left( {\theta  - \varphi^{\ast\ast\ast} } \right)} \right)} \right|}}{{\cos ^{\frac{{2N + 6}}{3}} \varphi^{\ast\ast\ast} }}\frac{2}{{3\pi }}6^{\frac{{2N + 6}}{3}} \left| {D_{2N + 5} \left( 0 \right)} \right|\frac{3}{4}\frac{{\Gamma \left( {\frac{{2N + 6}}{3}} \right)}}{{\left| \nu  \right|^{\frac{{2N + 6}}{3}} }},
\end{align*}
when $N \equiv 1 \mod 3$;
\[
\left| {R_N^{\left( {Y'} \right)} \left( \nu  \right)} \right| \le \frac{{\left| {\csc \left( {2\left( {\theta  - \varphi^{\ast\ast} } \right)} \right)} \right|}}{{\cos ^{\frac{{2N + 4}}{3}} \varphi^{\ast\ast} }}\frac{2}{{3\pi }}6^{\frac{{2N + 4}}{3}} \left| {D_{2N + 3} \left( 0 \right)} \right|\frac{3}{4}\frac{{\Gamma \left( {\frac{{2N + 4}}{3}} \right)}}{{\left| \nu  \right|^{\frac{{2N + 4}}{3}} }},
\]
for $N \equiv 2 \mod 3$. Here $\varphi^\ast$, $\varphi^{\ast\ast}$ and $\varphi^{\ast\ast\ast}$ are the same as in the bounds for $R_N^{\left( J' \right)} \left( \nu  \right)$.

If $\frac{\pi}{4} < \left|\arg \nu\right| \leq \frac{\pi}{2}$ and $N\geq 4$, we have the following simpler bounds useful near the Stokes lines $\arg \nu =\pm \frac{\pi}{2}$:
\[
\left| {R_N^{\left( {Y'} \right)} \left( \nu  \right)} \right| \le \frac{1}{2}\sqrt {\frac{e}{3}\left( {2N + \frac{{13}}{2}} \right)} \frac{2}{{3\pi }}6^{\frac{{2N + 2}}{3}} \left| {D_{2N + 1} \left( 0 \right)} \right|\frac{3}{4}\frac{{\Gamma \left( {\frac{{2N + 2}}{3}} \right)}}{{\left| \nu  \right|^{\frac{{2N + 2}}{3}} }},
\]
\begin{align*}
\left| {R_N^{\left( {Y'} \right)} \left( \nu  \right)} \right| \le \; & \frac{1}{2}\sqrt {\frac{e}{3}\left( {2N + \frac{{13}}{2}} \right)} \frac{2}{{3\pi }}6^{\frac{{2N + 2}}{3}} \left| {D_{2N + 1} \left( 0 \right)} \right|\frac{3}{4}\frac{{\Gamma \left( {\frac{{2N + 2}}{3}} \right)}}{{\left| \nu  \right|^{\frac{{2N + 2}}{3}} }} \\ &+ \frac{1}{2}\sqrt {\frac{e}{3}\left( {2N + \frac{{21}}{2}} \right)} \frac{2}{{3\pi }}6^{\frac{{2N + 6}}{3}} \left| {D_{2N + 5} \left( 0 \right)} \right|\frac{3}{4}\frac{{\Gamma \left( {\frac{{2N + 6}}{3}} \right)}}{{\left| \nu  \right|^{\frac{{2N + 6}}{3}} }},
\end{align*}
\[
\left| {R_N^{\left( {Y'} \right)} \left( \nu  \right)} \right| \le \frac{1}{2}\sqrt {\frac{e}{3}\left( {2N + \frac{{17}}{2}} \right)} \frac{2}{{3\pi }}6^{\frac{{2N + 4}}{3}} \left| {D_{2N + 3} \left( 0 \right)} \right|\frac{3}{4}\frac{{\Gamma \left( {\frac{{2N + 4}}{3}} \right)}}{{\left| \nu  \right|^{\frac{{2N + 4}}{3}} }}
\]
according to whether $N\equiv 0 \mod 3$, $N\equiv 1 \mod 3$ or $N\equiv 2 \mod 3$, respectively.

\section{Asymptotics for the late coefficients}\label{section4}

In this section, we investigate the asymptotic nature of the coefficients $B_n \left( \kappa  \right)$ as $n \to +\infty$. For our purposes, the most appropriate representation of these coefficients is the second integral formula in \eqref{eq8}. From \eqref{eq12}, it follows that for any $t > 0$, $M\geq 0$ and $\left| {\Re \left( \kappa  \right)} \right| < \frac{M + 1}{3}$, it holds that
\[
H_{it \pm \kappa }^{\left( 1 \right)} \left( {it} \right) = \frac{2}{{3\pi i}}\sum\limits_{m = 0}^{M - 1} {6^{\frac{{m + 1}}{3}} i^m B_m \left( { \mp \kappa } \right)\sin \left( {\frac{{\left( {m + 1} \right)\pi }}{3}} \right)\frac{{\Gamma \left( {\frac{{m + 1}}{3}} \right)}}{{t^{\frac{{m + 1}}{3}} }}}  + R_M^{\left( {H} \right)} \left( {it, \mp \kappa } \right).
\]
Substituting into \eqref{eq8} and using the property $B_m\left(-\kappa\right)=\left(-1\right)^m B_m\left(\kappa\right)$, gives the expansion
\begin{gather}\label{eq47}
\begin{split}
B_n \left( \kappa  \right) = \; & \frac{2}{{3\pi \left( {12\pi } \right)^{\frac{n}{3}} }}\left(\sum\limits_{m = 0}^{M - 1} {\left( {12\pi } \right)^{\frac{m}{3}} B_m \left( \kappa  \right)\sin \left( {\frac{{\left( {m + 1} \right)\pi }}{3}} \right)\Gamma \left( {\frac{{m + 1}}{3}} \right)\frac{{\Gamma \left( {\frac{{n - m}}{3}} \right)}}{{\Gamma \left( {\frac{{n + 1}}{3}} \right)}}} \right. \\ & \Biggl.\times \cos \left( {2\pi \kappa  - \frac{\pi }{2}\left( {n + m} \right)} \right)+ A_M \left( {n,\kappa } \right)\Biggr),
\end{split}
\end{gather}
for $0 \le M \le n - 1$ and $\left| {\Re \left( \kappa  \right)} \right| < \frac{M + 1}{3}$, provided that $n\geq 1$. The remainder term $A_M \left( {n,\kappa } \right)$ is given by the integral formula
\[
A_M \left( {n,\kappa } \right) = \frac{{3\pi \left( {2\pi } \right)^{\frac{n}{3}} }}{{6^{\frac{1}{3}} 4\Gamma \left( {\frac{{n + 1}}{3}} \right)}}\int_0^{ + \infty } {t^{\frac{{n - 2}}{3}} e^{ - 2\pi t} i\left( {e^{\left( {2\pi \kappa  - \frac{\pi }{2}n} \right)i} R_M^{\left( {H} \right)} \left( {it, - \kappa } \right) + e^{ - \left( {2\pi \kappa  - \frac{\pi }{2}n} \right)i} R_M^{\left( {H} \right)} \left( {it,\kappa } \right)} \right)dt} .
\]
In Appendix \ref{appendixb}, we show that if $\left|\Re \left( \kappa  \right)\right| < \frac{M + 1}{3}$, then
\[
\left| {R_M^{\left( {H} \right)} \left( {it,\pm\kappa} \right)} \right| \le \frac{{C_M \left( \pm\kappa \right)}}{{t^{\frac{M + 1}{3}} }},
\]
where the constants $C_M \left( \pm\kappa \right)$ depend only on $M$ and $\kappa$. Applying this result, trivial estimation shows that
\[
A_M \left( {n,\kappa } \right) = \mathcal{O}_{\kappa ,M} \left( {\frac{{\Gamma \left( {\frac{{n - M}}{3}} \right)}}{{\Gamma \left( {\frac{{n + 1}}{3}} \right)}}} \right)
\]
for fixed $\kappa$ and $M$ as $n\to +\infty$. For fixed $m\geq 0$ and large $n$, we have
\begin{equation}\label{eq71}
\frac{{\Gamma \left( {\frac{{n - m}}{3}} \right)}}{{\Gamma \left( {\frac{{n + 1}}{3}} \right)}} \sim \left( {\frac{3}{n}} \right)^{\frac{{m + 1}}{3}} ,
\end{equation}
and therefore
\begin{equation}\label{eq61}
B_n \left( \kappa  \right) \sim \frac{2}{{3\pi \left( {12\pi } \right)^{\frac{n}{3}} }}\sum\limits_{m = 0}^{\infty} {\left( {12\pi } \right)^{\frac{m}{3}} B_m \left( \kappa  \right)\sin \left( {\frac{{\left( {m + 1} \right)\pi }}{3}} \right)\Gamma \left( {\frac{{m + 1}}{3}} \right)\frac{{\Gamma \left( {\frac{{n - m}}{3}} \right)}}{{\Gamma \left( {\frac{{n + 1}}{3}} \right)}}\cos \left( {2\pi \kappa  - \frac{\pi }{2}\left( {n + m} \right)} \right)} ,
\end{equation}
as $n\to +\infty$ with $\kappa$ being fixed. Expansions of this type are called inverse factorial series in the literature. It is seen from \eqref{eq71}, that numerically, the character of the expansion \eqref{eq61} is similar to the character of standard Poincar\'e asymptotic power series.

Using the leading order behaviour of the $B_m \left( \kappa  \right)$'s, it is easy to show that the index of the numerically least term in the expansion \eqref{eq47} (or in \eqref{eq61}) is $m \approx \frac{n}{2}$. Therefore, by neglecting the remainder term $A_M \left( {n,\kappa } \right)$ and choosing the optimal $M\approx \frac{n}{2}$, the resulting approximation is the best possible we can get from \eqref{eq47}. Numerical examples illustrating the efficacy of the expansion \eqref{eq47}, truncated optimally, are given in Table \ref{table1}.

\begin{table*}[!ht]
\begin{center}
\begin{tabular}
[c]{ l r @{\,}c@{\,} l}\hline
 & \\ [-1ex]
 values of $n$, $\kappa$ and $M$ & $n=100$, $\kappa=3$, $M=50$ & & \\ [1ex]
 exact numerical value of $\Gamma\left(\frac{n+1}{3}\right)\left|B_n\left(\kappa\right)\right|$ & $0.7745012865285354362490235$ & $\times$ & $10^{-17}$ \\ [1ex]
 approximation \eqref{eq47} to $\Gamma\left(\frac{n+1}{3}\right)\left|B_n\left(\kappa\right)\right|$ & $0.7745012865519805241476135$ & $\times$ & $10^{-17}$  \\ [1ex]
 error & $-0.234450878985899$ & $\times$ & $10^{-27}$\\ [1ex] \hline
 & \\ [-1ex]
 values of $n$, $\kappa$ and $M$ &  $n=100$, $\kappa=2+2i$, $M=50$ & & \\ [1ex]
 exact numerical value of $\Gamma\left(\frac{n+1}{3}\right)\left|B_n\left(\kappa\right)\right|$ & $0.5529397074469944063403455$ & $\times$ & $10^{-12}$ \\ [1ex]
  approximation \eqref{eq47} to $\Gamma\left(\frac{n+1}{3}\right)\left|B_n\left(\kappa\right)\right|$ & $0.5529395059975027990719201$ & $\times$ & $10^{-12}$ \\ [1ex]
 error & $0.2014494916072684254$ & $\times$ & $10^{-18}$\\ [1ex] \hline
 & \\ [-1ex]
 values of $n$, $\kappa$ and $M$ &  $n=200$, $\kappa=5$, $M=100$ & & \\ [1ex]
 exact numerical value of $\Gamma\left(\frac{n+1}{3}\right)\left|B_n\left(\kappa\right)\right|$ & $0.2945913249283174576021141$ & $\times$ & $10^{-12}$ \\ [1ex]
  approximation \eqref{eq47} to $\Gamma\left(\frac{n+1}{3}\right)\left|B_n\left(\kappa\right)\right|$ & $0.2945913249283174576024119$ & $\times$ & $10^{-12}$ \\ [1ex]
 error & $-0.2978$ & $\times$ & $10^{-33}$\\ [1ex] \hline
 & \\ [-1ex]
 values of $n$, $\kappa$ and $M$ &  $n=200$, $\kappa=4+3i$, $M=100$ & & \\ [1ex]
 exact numerical value of $\Gamma\left(\frac{n+1}{3}\right)\left|B_n\left(\kappa\right)\right|$ & $0.1508584308199912914799076$ & $\times$ & $10^{-5}$ \\ [1ex]
 approximation \eqref{eq47} to $\Gamma\left(\frac{n+1}{3}\right)\left|B_n\left(\kappa\right)\right|$ & $0.1508584308199923691209822$ & $\times$ & $10^{-5}$ \\ [1ex]
 error & $-0.10776410746$ & $\times$ & $10^{-19}$\\ [-1ex] 
  & \\\hline
\end{tabular}
\end{center}
\caption{Approximations for $\Gamma\left(\frac{n+1}{3}\right)\left|B_n\left(\kappa\right)\right|$ with various $n$ and $\kappa$, using \eqref{eq47}.}
\label{table1}
\end{table*}

As a consequence of a result due to Cauchy and Meissel, Watson gave the asymptotic approximation
\[
B_{2n} \left( 0  \right) \sim \frac{{\left(-1\right)^n\left( {\frac{2}{3}} \right)^{\frac{2}{3}} }}{{\Gamma \left( {\frac{2}{3}} \right)\left( {n + \frac{1}{3}} \right)^{\frac{1}{3}} \left( {12\pi } \right)^{\frac{2n}{3}} }}
\]
as $n\to +\infty$, and remarked that there seems to be no very simple approximate formula for the general case. From the first two terms of \eqref{eq47}, one can prove the following simple asymptotic formula
\[
B_n \left( \kappa  \right) \sim \frac{{\left( {\frac{2}{3}} \right)^{\frac{2}{3}} }}{{\Gamma \left( {\frac{2}{3}} \right)\left( {\frac{n}{2} + \frac{1}{3}} \right)^{\frac{1}{3}} \left( {12\pi } \right)^{\frac{n}{3}} }}\cos \left( {2\pi \kappa  - \frac{\pi }{2}n} \right) + \frac{{\left( {\frac{2}{3}} \right)^{\frac{1}{3}} }}{{\Gamma \left( {\frac{1}{3}} \right)\left( {\frac{n}{2} + \frac{1}{3}} \right)^{\frac{2}{3}} \left( {12\pi } \right)^{\frac{{n - 1}}{3}} }}\kappa \sin \left( {2\pi \kappa  - \frac{\pi }{2}n} \right)
\]
as $n\to +\infty$, which reduces to Watson's approximation when $n=2n$ and $\kappa=0$.

From the formula $2D_n \left( \kappa  \right) = B_n \left( {\kappa  + 1} \right) - B_n \left( {\kappa  - 1} \right)$, it follows that $D_n \left( \kappa  \right)$ has the asymptotic series given in the right-hand side of \eqref{eq61} with $B_m \left( \kappa  \right)$ being replaced by $D_m \left( \kappa  \right)$. As a special case, we find
\begin{align*}
 - \frac{{6^{\frac{{2n + 1}}{3}} }}{4}\frac{{\Gamma \left( {\frac{{2n + 2}}{3}} \right)}}{{\Gamma \left( {\frac{1}{3}} \right)}}D_{2n + 1} \left( 0 \right) \sim \; & \left( {\frac{{16}}{{\pi ^2 }}} \right)^{\frac{1}{3}} \beta \sqrt 3 \frac{{\left( { - 1} \right)^n }}{{\left( {2\pi } \right)^{\frac{{2n + 1}}{3}} }}\left( {\Gamma \left( {\frac{{2n + 1}}{3} - \frac{1}{3}} \right) + \frac{{\left( {2\pi } \right)^{\frac{2}{3}} }}{{30\beta }}\Gamma \left( {\frac{{2n + 1}}{3} - 1} \right) }\right. \\ & \left.{ - \frac{{46\pi ^2 }}{{1575}}\Gamma \left( {\frac{{2n + 1}}{3} - \frac{7}{3}} \right) +  \cdots } \right),
\end{align*}
which is Dingle's \cite[p. 172]{Dingle} expansion for the late coefficients of the asymptotic series of $J'_\nu\left(\nu\right)$ and $Y'_\nu\left(\nu\right)$. The constant $\beta$ is frequently used by Dingle in his discussion of the asymptotic expansions of integrals with second order saddle points. Its value is given by
\[
\beta  = \frac{{6^{\frac{1}{3}} \Gamma \left( {\frac{5}{3}} \right)}}{{4\Gamma \left( {\frac{1}{3}} \right)}} = \frac{{\Gamma \left( {\frac{5}{6}} \right)}}{{2 \cdot 3^{\frac{2}{3}} \sqrt \pi}} \approx 0.1530827453.
\]
We remark that Dingle proved his expansion by formal, non-rigorous methods.

More accurate approximations could be derived for the coefficients $B_n \left( \kappa  \right)$ by estimating the remainder $A_M \left( {n,\kappa } \right)$, rather than bounding it, but we do not discuss the details here.

\section{Exponentially improved asymptotic expansions}\label{section5} We shall find it convenient to express our exponentially improved expansions in terms of the (scaled) Terminant function, which is defined in terms of the Incomplete gamma function as
\[
\widehat T_p \left( w \right) = \frac{{e^{\pi ip} \Gamma \left( p \right)}}{{2\pi i}}\Gamma \left( {1 - p,w} \right) = \frac{e^{\pi ip} w^{1 - p} e^{ - w} }{2\pi i}\int_0^{ + \infty } {\frac{{t^{p - 1} e^{ - t} }}{w + t}dt} \; \text{ for } \; p>0 \; \text{ and } \; \left| \arg w \right| < \pi ,
\]
and by analytic continuation elsewhere. Olver \cite[equations (4.5) and (4.6)]{Olver4} showed that when $p \sim \left|w\right|$ and $w \to \infty$, we have
\begin{equation}\label{eq55}
\widehat T_p \left( w \right) = \begin{cases} \mathcal{O}\left( {e^{ - w - \left| w \right|} } \right) & \; \text{ if } \; \left| {\arg w} \right| \le \pi, \\ \mathcal{O}\left(1\right) & \; \text{ if } \; - 3\pi  < \arg w \le  - \pi. \end{cases}
\end{equation}
Concerning the smooth transition of the Stokes discontinuities, we will use the more precise asymptotic formula
\begin{equation}\label{eq57}
e^{ - 2\pi ip} \widehat T_p \left( w \right) =  - \frac{1}{2} + \frac{1}{2}\mathop{\text{erf}} \left( { - \overline {c\left( { - \varphi } \right)} \sqrt {\frac{1}{2}\left| w \right|} } \right) + \mathcal{O}\left( {\frac{{e^{ - \frac{1}{2}\left| w \right|\overline {c^2 \left( { - \varphi } \right)} } }}{{\left| w \right|^{\frac{1}{2}} }}} \right)
\end{equation}
for $- 3\pi  + \delta  \le \arg w \le \pi  - \delta$, $0 < \delta \le 2\pi$. Here $\varphi = \arg w$ and $\mathop{\text{erf}}$ denotes the Error function. The quantity $c\left( \varphi  \right)$ is defined implicitly by the equation
\[
\frac{1}{2}c^2 \left( \varphi  \right) = 1 + i\left( {\varphi  - \pi } \right) - e^{i\left( {\varphi  - \pi } \right)},
\]
and corresponds to the branch of $c\left( \varphi  \right)$ which has the following expansion in the neighbourhood of $\varphi = \pi$:
\begin{equation}\label{eq58}
c\left( \varphi  \right) = \left( {\varphi  - \pi } \right) + \frac{i}{6}\left( {\varphi  - \pi } \right)^2  - \frac{1}{{36}}\left( {\varphi  - \pi } \right)^3  - \frac{i}{{270}}\left( {\varphi  - \pi } \right)^4  +  \cdots .
\end{equation}
For complete asymptotic expansions, see Olver \cite{Olver5}. We remark that Olver uses the different notation $F_p \left( w \right) = ie^{ - \pi ip} \widehat T_p \left( w \right)$ for the Terminant function and the other branch of the function $c\left( \varphi  \right)$. For further properties of the Terminant function, see, for example, Paris and Kaminski \cite[Chapter 6]{Paris3}.

\subsection{Proof of Theorem \ref{thm3}} First, we suppose that $\left|\arg z\right|<\frac{\pi}{2}$. We write \eqref{eq12} with $N=0$ in the form
\begin{align*}
H_\nu ^{\left( 1 \right)} \left( z \right) = & \frac{{e^{ - \frac{\pi }{3}i} }}{{\sqrt 3 \pi }}\frac{{ie^{2\pi i\kappa } }}{{2z^{\frac{1}{3}} }}\int_0^{ + \infty } {\frac{{t^{ - \frac{2}{3}} e^{ - 2\pi t} }}{{1 + it/z}}H_{it + \kappa }^{\left( 1 \right)} \left( {it} \right)dt}  + \frac{{e^{ - \frac{\pi }{3}i} }}{{\sqrt 3 \pi }}\frac{{ie^{ - 2\pi i\kappa } }}{{2z^{\frac{1}{3}} }}\int_0^{ + \infty } {\frac{{t^{ - \frac{2}{3}} e^{ - 2\pi t} }}{{1 - it/z}}H_{it - \kappa }^{\left( 1 \right)} \left( {it} \right)dt} 
\\ & + \frac{{e^{\frac{\pi }{3}i} }}{{\sqrt 3 \pi }}\frac{{e^{2\pi i\kappa } }}{{2z^{\frac{2}{3}} }}\int_0^{ + \infty } {\frac{{t^{ - \frac{1}{3}} e^{ - 2\pi t} }}{{1 + it/z}}H_{it + \kappa }^{\left( 1 \right)} \left( {it} \right)dt}  - \frac{{e^{\frac{\pi }{3}i} }}{{\sqrt 3 \pi }}\frac{{e^{ - 2\pi i\kappa } }}{{2z^{\frac{2}{3}} }}\int_0^{ + \infty } {\frac{{t^{ - \frac{1}{3}} e^{ - 2\pi t} }}{{1 - it/z}}H_{it - \kappa }^{\left( 1 \right)} \left( {it} \right)dt} .
\end{align*}
Let $N$ and $M$ be arbitrary positive integers and suppose that $\left| {\Re \left( \kappa  \right)} \right| < \min \left( {N - \frac{2}{3},M - \frac{1}{3}} \right)$. Using the expression \eqref{eq48} and the second representation in \eqref{eq8}, we find that
\begin{multline*}
H_\nu ^{\left( 1 \right)} \left( z \right) = \frac{{e^{ - \frac{\pi }{3}i} }}{{\sqrt 3 \pi z^{\frac{1}{3}} }}\sum\limits_{n = 0}^{N - 1} {\left( { - 1} \right)^n 6^{n + \frac{1}{3}} B_{3n} \left( \kappa  \right)\frac{{\Gamma \left( {n + \frac{1}{3}} \right)}}{{z^n }}} \\ + \frac{{e^{\frac{\pi }{3}i} }}{{\sqrt 3 \pi z^{\frac{2}{3}} }}\sum\limits_{m = 0}^{M - 1} {\left( { - 1} \right)^m 6^{m + \frac{2}{3}} B_{3m + 1} \left( \kappa  \right)\frac{{\Gamma \left( {m + \frac{2}{3}} \right)}}{{z^m }}}  + R_{N,M}^{\left( {H} \right)} \left( {z,\kappa } \right),
\end{multline*}
where
\begin{gather}\label{eq49}
\begin{split}
R_{N,M}^{\left( {H} \right)} \left( {z,\kappa } \right) = \; & \left( { - i} \right)^N \frac{{e^{ - \frac{\pi }{3}i} }}{{\sqrt 3 \pi }}\frac{{ie^{2\pi i\kappa } }}{{2z^{N + \frac{1}{3}} }}\int_0^{ + \infty } {\frac{{t^{N - \frac{2}{3}} e^{ - 2\pi t} }}{{1 + it/z}}H_{it + \kappa }^{\left( 1 \right)} \left( {it} \right)dt}  \\ & + i^N \frac{{e^{ - \frac{\pi }{3}i} }}{{\sqrt 3 \pi }}\frac{{ie^{ - 2\pi i\kappa } }}{{2z^{N + \frac{1}{3}} }}\int_0^{ + \infty } {\frac{{t^{N - \frac{2}{3}} e^{ - 2\pi t} }}{{1 - it/z}}H_{it - \kappa }^{\left( 1 \right)} \left( {it} \right)dt} 
\\ & + \left( { - i} \right)^M \frac{{e^{\frac{\pi }{3}i} }}{{\sqrt 3 \pi }}\frac{{e^{2\pi i\kappa } }}{{2z^{M + \frac{2}{3}} }}\int_0^{ + \infty } {\frac{{t^{M - \frac{1}{3}} e^{ - 2\pi t} }}{{1 + it/z}}H_{it + \kappa }^{\left( 1 \right)} \left( {it} \right)dt} \\ & - i^M \frac{{e^{\frac{\pi }{3}i} }}{{\sqrt 3 \pi }}\frac{{e^{ - 2\pi i\kappa } }}{{2z^{M + \frac{2}{3}} }}\int_0^{ + \infty } {\frac{{t^{M - \frac{1}{3}} e^{ - 2\pi t} }}{{1 - it/z}}H_{it - \kappa }^{\left( 1 \right)} \left( {it} \right)dt} .
\end{split}
\end{gather}
Note that $R_{N,N}^{\left( H \right)} \left( {z,\kappa } \right) = R_{3N}^{\left( H \right)} \left( {z,\kappa } \right)$. Assume that $K$ and $L$ are integers such that $0\leq K< 3N$ and $0\leq L< 3M+1$. We apply \eqref{eq12} again to expand the functions $H_{it\pm \kappa }^{\left( 1 \right)} \left( {it} \right)$ under the integral in \eqref{eq49}, to obtain
\begin{align*}
R_{N,M}^{\left( {H} \right)} \left( {z,\kappa } \right) = \; & ie^{ - \frac{\pi }{3}i} \frac{{e^{ - 2\pi i\nu } }}{{\sqrt 3 }}\frac{2}{{3\pi }}\sum\limits_{k = 0}^{K - 1} {6^{\frac{{k + 1}}{3}} B_k \left( \kappa  \right)\sin \left( {\frac{{\left( {k + 1} \right)\pi }}{3}} \right)\frac{{\Gamma \left( {\frac{{k + 1}}{3}} \right)}}{{z^{\frac{{k + 1}}{3}} }}\widehat T_{N - \frac{k}{3}} \left( { - 2\pi iz} \right)} \\
& - i\frac{{e^{2\pi i\nu } }}{{\sqrt 3 }}\frac{2}{{3\pi }}\sum\limits_{k = 0}^{K - 1} {6^{\frac{{k + 1}}{3}} B_k \left( \kappa  \right)e^{\frac{{2\left( {k + 1} \right)\pi i}}{3}} \sin \left( {\frac{{\left( {k + 1} \right)\pi }}{3}} \right)\frac{{\Gamma \left( {\frac{{k + 1}}{3}} \right)}}{{z^{\frac{{k + 1}}{3}} }}\widehat T_{N - \frac{k}{3}} \left( {2\pi iz} \right)} \\
& - ie^{\frac{\pi }{3}i} \frac{{e^{ - 2\pi i\nu } }}{{\sqrt 3 }}\frac{2}{{3\pi }}\sum\limits_{\ell  = 0}^{L - 1} {6^{\frac{{\ell  + 1}}{3}} B_\ell  \left( \kappa  \right)\sin \left( {\frac{{\left( {\ell  + 1} \right)\pi }}{3}} \right)\frac{{\Gamma \left( {\frac{{\ell  + 1}}{3}} \right)}}{{z^{\frac{{\ell  + 1}}{3}} }}\widehat T_{M - \frac{{\ell  - 1}}{3}} \left( { - 2\pi iz} \right)} \\
& + i\frac{{e^{2\pi i\nu } }}{{\sqrt 3 }}\frac{2}{{3\pi }}\sum\limits_{\ell  = 0}^{L - 1} {6^{\frac{{\ell  + 1}}{3}} B_\ell  \left( \kappa  \right)e^{\frac{{2\left( {\ell  + 1} \right)\pi i}}{3}} \sin \left( {\frac{{\left( {\ell  + 1} \right)\pi }}{3}} \right)\frac{{\Gamma \left( {\frac{{\ell  + 1}}{3}} \right)}}{{z^{\frac{{\ell  + 1}}{3}} }}\widehat T_{M - \frac{{\ell  - 1}}{3}} \left( {2\pi iz} \right)} \\
& + R_{N,M,K,L}^{\left( {H} \right)} \left( {z,\kappa } \right),
\end{align*}
with
\begin{gather}\label{eq51}
\begin{split}
R_{N,M,K,L}^{\left( {H} \right)} \left( {z,\kappa } \right) = \; & \left( { - i} \right)^N \frac{{e^{ - \frac{\pi }{3}i} }}{{\sqrt 3 \pi }}\frac{{ie^{2\pi i\kappa } }}{{2z^{N + \frac{1}{3}} }}\int_0^{ + \infty } {\frac{{t^{N - \frac{2}{3}} e^{ - 2\pi t} }}{{1 + it/z}}R_K^{\left( {H} \right)} \left( {it, - \kappa } \right)dt} \\ & + i^N \frac{{e^{ - \frac{\pi }{3}i} }}{{\sqrt 3 \pi }}\frac{{ie^{ - 2\pi i\kappa } }}{{2z^{N + \frac{1}{3}} }}\int_0^{ + \infty } {\frac{{t^{N - \frac{2}{3}} e^{ - 2\pi t} }}{{1 - it/z}}R_K^{\left( {H} \right)} \left( {it,\kappa } \right)dt} 
\\ & + \left( { - i} \right)^M \frac{{e^{\frac{\pi }{3}i} }}{{\sqrt 3 \pi }}\frac{{e^{2\pi i\kappa } }}{{2z^{M + \frac{2}{3}} }}\int_0^{ + \infty } {\frac{{t^{M - \frac{1}{3}} e^{ - 2\pi t} }}{{1 + it/z}}R_L^{\left( {H} \right)} \left( {it, - \kappa } \right)dt} \\ & - i^M \frac{{e^{\frac{\pi }{3}i} }}{{\sqrt 3 \pi }}\frac{{e^{ - 2\pi i\kappa } }}{{2z^{M + \frac{2}{3}} }}\int_0^{ + \infty } {\frac{{t^{M - \frac{1}{3}} e^{ - 2\pi t} }}{{1 - it/z}}R_L^{\left( {H} \right)} \left( {it,\kappa } \right)dt} .
\end{split}
\end{gather}
Here we have made use of the definition of the Terminant function, and the properties $\nu = z-\kappa$ and $B_n\left(-\kappa\right)=\left(-1\right)^n B_n\left(\kappa\right)$. We shall estimate the remainder $R_{N,M,K,L}^{\left( {H} \right)} \left( {z,\kappa } \right)$ in the closed sector $\left|\arg z\right|\leq\frac{\pi}{2}$, in which $R_{N,M,K,L}^{\left( {H} \right)} \left( {z,\kappa } \right)$ is well-defined by analytic continuation. In Appendix \ref{appendixb}, we show that there is a constant $C_N \left( \kappa  \right)>0$ depending only on $N$ and $\kappa$, such that if $0 \leq \arg z \leq \pi$ and $\left|\Re \left( \kappa  \right)\right| < \frac{N + 1}{3}$, then
\begin{equation}\label{eq52}
\left| {R_N^{\left( {H} \right)} \left( {z,\kappa} \right)} \right| \le \frac{{C_N \left( \kappa \right)}}{{\left| z \right|^{\frac{{N + 1}}{3}} }}.
\end{equation}
Suppose that $\left| {\Re \left( \kappa  \right)} \right| < \min \left( {\frac{{K + 1}}{3},\frac{{L + 1}}{3}} \right)$. Employing the bound \eqref{eq52} and the inequality \eqref{eq50}, one finds that
\begin{gather}\label{eq53}
\begin{split}
\left| {R_{N,M,K,L}^{\left( {H} \right)} \left( {z,\kappa } \right)} \right| \le \; & \frac{1}{{\sqrt 3 \pi }}\frac{1}{{2\left| z \right|^{N + \frac{1}{3}} }}\frac{{\Gamma \left( {N - \frac{K}{3}} \right)}}{{\left( {2\pi } \right)^{N - \frac{K}{3}} }}\left( {C_K \left( { - \kappa } \right)\left| {e^{2\pi i\kappa } } \right| + C_K \left( \kappa  \right)\left| {e^{ - 2\pi i\kappa } } \right|} \right)\left| {\sec \theta } \right|
\\ & + \frac{1}{{\sqrt 3 \pi }}\frac{1}{{2\left| z \right|^{M + \frac{2}{3}} }}\frac{{\Gamma \left( {M - \frac{{L - 1}}{3}} \right)}}{{\left( {2\pi } \right)^{M - \frac{{L - 1}}{3}} }}\left( {C_L \left( { - \kappa } \right)\left| {e^{2\pi i\kappa } } \right| + C_L \left( \kappa  \right)\left| {e^{ - 2\pi i\kappa } } \right|} \right)\left| {\sec \theta } \right|,
\end{split}
\end{gather}
when $\left|\arg z\right|<\frac{\pi}{2}$ and $\left|\arg z\right|$ is not too close to $\frac{\pi}{2}$.

Now we consider the bound which is suitable near $\left|\arg z\right|=\frac{\pi}{2}$. Let $0<\varphi,\varphi'< \frac{\pi}{2}$ be two acute angles that may depend on $N$, $M$, $K$ and $L$. Suppose that $-\frac{\pi}{2}\leq \theta < -\max\left(\varphi,\varphi'\right)$, where $\theta = \arg z$. We rotate the path of integration in the first integral in \eqref{eq51} by $-\varphi$ and in the third integral by $-\varphi'$, and apply the estimates \eqref{eq52} and \eqref{eq50}, to deduce
\begin{align*}
\left| {R_{N,M,K,L}^{\left( {H} \right)} \left( {z,\kappa } \right)} \right| \le \; & \frac{{\sec \left( {\theta  + \varphi } \right)}}{{\cos ^{N - \frac{K}{3}} \varphi }}\frac{1}{{\sqrt 3 \pi }}\frac{1}{{2\left| z \right|^{N + \frac{1}{3}} }}\frac{{\Gamma \left( {N - \frac{K}{3}} \right)}}{{\left( {2\pi } \right)^{N - \frac{K}{3}} }}C_K \left( { - \kappa } \right)\left| {e^{2\pi i\kappa } } \right|
\\ & + \frac{1}{{\sqrt 3 \pi }}\frac{1}{{2\left| z \right|^{N + \frac{1}{3}} }}\frac{{\Gamma \left( {N - \frac{K}{3}} \right)}}{{\left( {2\pi } \right)^{N - \frac{K}{3}} }}C_K \left( \kappa  \right)\left| {e^{ - 2\pi i\kappa } } \right|
\\ & + \frac{{\sec \left( {\theta  + \varphi' } \right)}}{{\cos ^{M - \frac{{L - 1}}{3}} \varphi' }}\frac{1}{{\sqrt 3 \pi }}\frac{1}{{2\left| z \right|^{M + \frac{2}{3}} }}\frac{{\Gamma \left( {M - \frac{{L - 1}}{3}} \right)}}{{\left( {2\pi } \right)^{M - \frac{{L - 1}}{3}} }}C_L \left( { - \kappa } \right)\left| {e^{2\pi i\kappa } } \right|
\\ & + \frac{1}{{\sqrt 3 \pi }}\frac{1}{{2\left| z \right|^{M + \frac{2}{3}} }}\frac{{\Gamma \left( {M - \frac{{L - 1}}{3}} \right)}}{{\left( {2\pi } \right)^{M - \frac{{L - 1}}{3}} }}C_L \left( \kappa  \right)\left| {e^{ - 2\pi i\kappa } } \right|.
\end{align*}
To minimise the factors involving the trigonometric functions, we choose
\[
\varphi  = \arctan \left( {\left( {N - \frac{K}{3}} \right)^{ - \frac{1}{2}} } \right) \; \text{ and } \; \varphi'  = \arctan \left( {\left( {M - \frac{L-1}{3}} \right)^{ - \frac{1}{2}} } \right),
\]
and note that
\[
\frac{{\sec \left( {\theta  + \varphi } \right)}}{{\cos ^{N - \frac{K}{3}} \varphi }} \le \frac{{\sec \left( { - \frac{\pi }{2} + \varphi } \right)}}{{\cos ^{N - \frac{K}{3}} \varphi }} = \left( {1 + \frac{1}{{N - \frac{K}{3}}}} \right)^{\frac{1}{2}\left( {N - \frac{K}{3} + 1} \right)} \sqrt {N - \frac{K}{3}}  \le 2\sqrt {N - \frac{K}{3}} 
\]
and
\[
\frac{{\sec \left( {\theta  + \varphi '} \right)}}{{\cos ^{N - \frac{K}{3}} \varphi '}} \le \frac{{\sec \left( { - \frac{\pi }{2} + \varphi '} \right)}}{{\cos ^{N - \frac{K}{3}} \varphi '}} = \left( {1 + \frac{1}{{M - \frac{{L - 1}}{3}}}} \right)^{\frac{1}{2}\left( {M - \frac{{L - 1}}{3} + 1} \right)} \sqrt {M - \frac{{L - 1}}{3}}  \le 2\sqrt {M - \frac{{L - 1}}{3}} ,
\]
for any $-\frac{\pi}{2}\leq \theta < -\max\left(\varphi,\varphi'\right)$. Therefore, we have
\begin{align*}
\left| {R_{N,M,K,L}^{\left( {H} \right)} \left( {z,\kappa } \right)} \right| \le \; & \sqrt {N - \frac{K}{3}} \frac{1}{{\sqrt 3 \pi }}\frac{1}{{\left| z \right|^{N + \frac{1}{3}} }}\frac{{\Gamma \left( {N - \frac{K}{3}} \right)}}{{\left( {2\pi } \right)^{N - \frac{K}{3}} }}C_K \left( { - \kappa } \right)\left| {e^{2\pi i\kappa } } \right| \\
& + \frac{1}{{\sqrt 3 \pi }}\frac{1}{{2\left| z \right|^{N + \frac{1}{3}} }}\frac{{\Gamma \left( {N - \frac{K}{3}} \right)}}{{\left( {2\pi } \right)^{N - \frac{K}{3}} }}C_K \left( \kappa  \right)\left| {e^{ - 2\pi i\kappa } } \right|
\\ & + \sqrt {M - \frac{{L - 1}}{3}} \frac{1}{{\sqrt 3 \pi }}\frac{1}{{\left| z \right|^{M + \frac{2}{3}} }}\frac{{\Gamma \left( {M - \frac{{L - 1}}{3}} \right)}}{{\left( {2\pi } \right)^{M - \frac{{L - 1}}{3}} }}C_L \left( { - \kappa } \right)\left| {e^{2\pi i\kappa } } \right|
\\ & + \frac{1}{{\sqrt 3 \pi }}\frac{1}{{2\left| z \right|^{M + \frac{2}{3}} }}\frac{{\Gamma \left( {M - \frac{{L - 1}}{3}} \right)}}{{\left( {2\pi } \right)^{M - \frac{{L - 1}}{3}} }}C_L \left( \kappa  \right)\left| {e^{ - 2\pi i\kappa } } \right|
\end{align*}
when $\arg z$ is close to $-\frac{\pi}{2}$. Comparing this estimate with \eqref{eq53}, we conclude that
\[
R_{N,M,K,L}^{\left( {H} \right)} \left( {z,\kappa } \right) = \mathcal{O}_{\kappa ,K} \left( {\sqrt N \frac{1}{{\left| z \right|^{N + \frac{1}{3}} }}\frac{{\Gamma \left( {N - \frac{K}{3}} \right)}}{{\left( {2\pi } \right)^{N - \frac{K}{3}} }}} \right) + \mathcal{O}_{\kappa ,L} \left( {\sqrt M \frac{1}{{\left| z \right|^{M + \frac{2}{3}} }}\frac{{\Gamma \left( {M - \frac{{L - 1}}{3}} \right)}}{{\left( {2\pi } \right)^{M - \frac{{L - 1}}{3}} }}} \right)
\]
for $-\frac{\pi}{2} \leq \arg z \leq 0$. By the reflection principle, $\left| {R_{N,M,K,L}^{\left( {H} \right)} \left( {\bar z,\kappa } \right)} \right| = \left| {\overline {R_{N,M,K,L}^{\left( {H} \right)} \left( {z,\bar \kappa } \right)} } \right| = \left| {R_{N,M,K,L}^{\left( {H} \right)} \left( {z,\bar \kappa } \right)} \right|$, which shows that the above order estimate is true in the larger sector $\left|\arg z\right| \leq \frac{\pi}{2}$. Suppose that $N = 2\pi \left| z \right| + \rho$, $M = 2\pi \left| z \right| + \sigma$ where $\rho$ and $\sigma$ are bounded. An application of Stirling's formula shows that
\begin{equation}\label{eq54}
R_{N,M,K,L}^{\left( {H} \right)} \left( {z,\kappa } \right) = \mathcal{O}_{\kappa ,K,\rho } \left( {\frac{{e^{ - 2\pi \left| z \right|} }}{{\left| z \right|^{\frac{{K + 1}}{3}} }}} \right) + \mathcal{O}_{\kappa ,L,\sigma } \left( {\frac{{e^{ - 2\pi \left| z \right|} }}{{\left| z \right|^{\frac{{L + 1}}{3}} }}} \right)
\end{equation}
as $z\to \infty$ in the closed sector $\left|\arg z\right| \leq \frac{\pi}{2}$.

Next, we consider the sector $\frac{\pi}{2} \leq \arg z \leq \frac{3\pi}{2}$. Rotating the path of integration in the second and the fourth integrals in \eqref{eq51} and applying the residue theorem gives
\begin{gather}\label{eq72}
\begin{split}
& R_{N,M,K,L}^{\left( {H} \right)} \left( {z,\kappa } \right) = i\frac{{e^{2\pi i\nu } }}{{\sqrt 3 }}R_K^{\left( {H} \right)} \left( {z,\kappa } \right) - i\frac{{e^{2\pi i\nu } }}{{\sqrt 3 }}R_L^{\left( {H} \right)} \left( {z,\kappa } \right)
\\ & + \left( { - i} \right)^N \frac{{e^{ - \frac{\pi }{3}i} }}{{\sqrt 3 \pi }}\frac{{ie^{2\pi i\kappa } }}{{2z^{N + \frac{1}{3}} }}\int_0^{ + \infty } {\frac{{t^{N - \frac{2}{3}} e^{ - 2\pi t} }}{{1 + it/z}}R_K^{\left( {H} \right)} \left( {it, - \kappa } \right)dt}  + i^N \frac{{e^{ - \frac{\pi }{3}i} }}{{\sqrt 3 \pi }}\frac{{ie^{ - 2\pi i\kappa } }}{{2z^{N + \frac{1}{3}} }}\int_0^{ + \infty } {\frac{{t^{N - \frac{2}{3}} e^{ - 2\pi t} }}{{1 - it/z}}R_K^{\left( {H} \right)} \left( {it,\kappa } \right)dt} 
\\ & + \left( { - i} \right)^M \frac{{e^{\frac{\pi }{3}i} }}{{\sqrt 3 \pi }}\frac{{e^{2\pi i\kappa } }}{{2z^{M + \frac{2}{3}} }}\int_0^{ + \infty } {\frac{{t^{M - \frac{1}{3}} e^{ - 2\pi t} }}{{1 + it/z}}R_L^{\left( {H} \right)} \left( {it, - \kappa } \right)dt}  - i^M \frac{{e^{\frac{\pi }{3}i} }}{{\sqrt 3 \pi }}\frac{{e^{ - 2\pi i\kappa } }}{{2z^{M + \frac{2}{3}} }}\int_0^{ + \infty } {\frac{{t^{M - \frac{1}{3}} e^{ - 2\pi t} }}{{1 - it/z}}R_L^{\left( {H} \right)} \left( {it,\kappa } \right)dt} ,
\end{split}
\end{gather}
for $\frac{\pi}{2} < \arg z < \frac{3\pi}{2}$. The four integrals can be estimated using the techniques from the previous part (for example by writing $z = \left(ze^{-\pi i}\right)e^{\pi i}$), and we find that they satisfy the order estimate given in the right-hand side of \eqref{eq54}. It follows that when $K = L$, the bound \eqref{eq54} remains valid in the wider sector $-\frac{\pi}{2} \leq \arg z \leq \frac{3\pi}{2}$. Otherwise, we have
\begin{align*}
\left| {R_{N,M,K,L}^{\left( H \right)} \left( {z,\kappa } \right)} \right| \le \; & \frac{{e^{ - 2\pi \Im \left( z \right)} e^{2\pi \Im \left( \kappa  \right)} }}{{\sqrt 3 }}\left| {R_K^{\left( H \right)} \left( {z,\kappa } \right)} \right| + \frac{{e^{ - 2\pi \Im \left( z \right)} e^{2\pi \Im \left( \kappa  \right)} }}{{\sqrt 3 }}\left| {R_L^{\left( H \right)} \left( {z,\kappa } \right)} \right|
\\ & + \mathcal{O}_{\kappa ,K,\rho } \left( {\frac{{e^{ - 2\pi \left| z \right|} }}{{\left| z \right|^{\frac{{K + 1}}{3}} }}} \right) + \mathcal{O}_{\kappa ,L,\sigma } \left( {\frac{{e^{ - 2\pi \left| z \right|} }}{{\left| z \right|^{\frac{{L + 1}}{3}} }}} \right).
\end{align*}
In Appendix \ref{appendixb}, we show that $R_K^{\left( H \right)} \left( {z,\kappa } \right) = \mathcal{O}_{\kappa ,K,\delta} \left( {\left|z\right|^{ - \frac{{K + 1}}{3}} } \right)$ and $R_L^{\left( H \right)} \left( {z,\kappa } \right) = \mathcal{O}_{\kappa ,L,\delta} \left( {\left|z\right|^{ - \frac{{L + 1}}{3}} } \right)$ as $z\to \infty$ in $ - \pi  + \delta  \leq \arg z \leq 2\pi  - \delta$ ($0<\delta\leq \pi$), whence
\[
R_{N,M,K,L}^{\left( H \right)} \left( {z,\kappa } \right) = \mathcal{O}_{\kappa ,K,\rho } \left( {\frac{{e^{ - 2\pi \Im \left( z \right)} }}{{\left| z \right|^{\frac{{K + 1}}{3}} }}} \right) + \mathcal{O}_{\kappa ,L,\sigma } \left( {\frac{{e^{ - 2\pi \Im \left( z \right)} }}{{\left| z \right|^{\frac{{L + 1}}{3}} }}} \right)
\]
as $z\to \infty$ in the sector $\frac{\pi}{2} \leq \arg z \leq \frac{3\pi}{2}$, using continuity.

Consider now the case $-\frac{3\pi}{2} \leq \arg z \leq -\frac{\pi}{2}$. Rotating the path of integration in the first and the third integrals in \eqref{eq51} and applying the residue theorem yields
\begin{gather}\label{eq62}
\begin{split}
& R_{N,M,K,L}^{\left( H \right)} \left( {z,\kappa } \right) =  - ie^{\frac{\pi }{3}i} \frac{{e^{ - 2\pi i\nu } }}{{\sqrt 3 }}R_K^{\left( H \right)} \left( {ze^{\pi i} , - \kappa } \right) + ie^{ - \frac{\pi }{3}i} \frac{{e^{ - 2\pi i\nu } }}{{\sqrt 3 }}R_L^{\left( H \right)} \left( {ze^{\pi i} , - \kappa } \right)
\\ & + \left( { - i} \right)^N \frac{{e^{ - \frac{\pi }{3}i} }}{{\sqrt 3 \pi }}\frac{{ie^{2\pi i\kappa } }}{{2z^{N + \frac{1}{3}} }}\int_0^{ + \infty } {\frac{{t^{N - \frac{2}{3}} e^{ - 2\pi t} }}{{1 + it/z}}R_K^{\left( H \right)} \left( {it, - \kappa } \right)dt}  + i^N \frac{{e^{ - \frac{\pi }{3}i} }}{{\sqrt 3 \pi }}\frac{{ie^{ - 2\pi i\kappa } }}{{2z^{N + \frac{1}{3}} }}\int_0^{ + \infty } {\frac{{t^{N - \frac{2}{3}} e^{ - 2\pi t} }}{{1 - it/z}}R_K^{\left( H \right)} \left( {it,\kappa } \right)dt} 
\\ & + \left( { - i} \right)^M \frac{{e^{\frac{\pi }{3}i} }}{{\sqrt 3 \pi }}\frac{{e^{2\pi i\kappa } }}{{2z^{M + \frac{2}{3}} }}\int_0^{ + \infty } {\frac{{t^{M - \frac{1}{3}} e^{ - 2\pi t} }}{{1 + it/z}}R_L^{\left( H \right)} \left( {it, - \kappa } \right)dt}  - i^M \frac{{e^{\frac{\pi }{3}i} }}{{\sqrt 3 \pi }}\frac{{e^{ - 2\pi i\kappa } }}{{2z^{M + \frac{2}{3}} }}\int_0^{ + \infty } {\frac{{t^{M - \frac{1}{3}} e^{ - 2\pi t} }}{{1 - it/z}}R_L^{\left( H \right)} \left( {it,\kappa } \right)dt} 
\end{split}
\end{gather}
for $-\frac{3\pi}{2} < \arg z < -\frac{\pi}{2}$. It is easy to see (for example by writing $z = \left(ze^{\pi i}\right)e^{-\pi i}$) that the sum of the four integrals has the order of magnitude given in the right-hand side of \eqref{eq54}. Since
\begin{multline*}
\left| { - ie^{\frac{\pi }{3}i} \frac{{e^{ - 2\pi i\nu } }}{{\sqrt 3 }}R_K^{\left( H \right)} \left( {ze^{\pi i} , - \kappa } \right) + ie^{ - \frac{\pi }{3}i} \frac{{e^{ - 2\pi i\nu } }}{{\sqrt 3 }}R_L^{\left( H \right)} \left( {ze^{\pi i} , - \kappa } \right)} \right| \\ \le \frac{{e^{2\pi \Im \left( z \right)} e^{ - 2\pi \Im \left( \kappa  \right)} }}{{\sqrt 3 }}\left| {R_K^{\left( H \right)} \left( {ze^{\pi i} , - \kappa } \right)} \right| + \frac{{e^{2\pi \Im \left( z \right)} e^{ - 2\pi \Im \left( \kappa  \right)} }}{{\sqrt 3 }}\left| {R_L^{\left( H \right)} \left( {ze^{\pi i} , - \kappa } \right)} \right|,
\end{multline*}
and $R_K^{\left( H \right)} \left( {ze^{\pi i},-\kappa } \right) = \mathcal{O}_{\kappa ,K,\delta} \left( {\left|z\right|^{ - \frac{{K + 1}}{3}} } \right)$ and $R_L^{\left( H \right)} \left( {ze^{\pi i},-\kappa } \right) = \mathcal{O}_{\kappa ,L,\delta} \left( {\left|z\right|^{ - \frac{{L + 1}}{3}} } \right)$ as $z\to \infty$ in $ - 2\pi  + \delta \leq \arg z \leq \pi  - \delta$ ($0<\delta\leq \pi$), we find
\begin{equation}\label{eq73}
R_{N,M,K,L}^{\left( H \right)} \left( {z,\kappa } \right) = \mathcal{O}_{\kappa ,K,\rho } \left( {\frac{{e^{ 2\pi \Im \left( z \right)} }}{{\left| z \right|^{\frac{{K + 1}}{3}} }}} \right) + \mathcal{O}_{\kappa ,L,\sigma } \left( {\frac{{e^{ 2\pi \Im \left( z \right)} }}{{\left| z \right|^{\frac{{L + 1}}{3}} }}} \right)
\end{equation}
as $z\to \infty$ in the sector $-\frac{3\pi}{2} \leq \arg z \leq -\frac{\pi}{2}$, using continuity.

Let $0 < \delta  \le \frac{{3\pi }}{2}$ be a fixed real number. We consider the sector $\frac{{3\pi }}{2} \le \arg z \le 3\pi  - \delta <3\pi$. Rotating the path of integration in the first and the third integrals in \eqref{eq72} and applying the residue theorem yields
\begin{align*}
& R_{N,M,K,L}^{\left( H \right)} \left( {z,\kappa } \right) = ie^{ - \frac{\pi }{3}i} \frac{{e^{ - 2\pi i\nu } }}{{\sqrt 3 }}R_K^{\left( H \right)} \left( {ze^{ - \pi i} , - \kappa } \right) - ie^{\frac{\pi }{3}i} \frac{{e^{ - 2\pi i\nu } }}{{\sqrt 3 }}R_L^{\left( H \right)} \left( {ze^{ - \pi i} , - \kappa } \right)
\\ & +i\frac{{e^{2\pi i\nu } }}{{\sqrt 3 }}R_K^{\left( H \right)} \left( {z,\kappa } \right) - i\frac{{e^{2\pi i\nu } }}{{\sqrt 3 }}R_L^{\left( H \right)} \left( {z,\kappa } \right)
\\ &  + \left( { - i} \right)^N \frac{{e^{ - \frac{\pi }{3}i} }}{{\sqrt 3 \pi }}\frac{{ie^{2\pi i\kappa } }}{{2z^{N + \frac{1}{3}} }}\int_0^{ + \infty } {\frac{{t^{N - \frac{2}{3}} e^{ - 2\pi t} }}{{1 + it/z}}R_K^{\left( H \right)} \left( {it, - \kappa } \right)dt}  + i^N \frac{{e^{ - \frac{\pi }{3}i} }}{{\sqrt 3 \pi }}\frac{{ie^{ - 2\pi i\kappa } }}{{2z^{N + \frac{1}{3}} }}\int_0^{ + \infty } {\frac{{t^{N - \frac{2}{3}} e^{ - 2\pi t} }}{{1 - it/z}}R_K^{\left( H \right)} \left( {it,\kappa } \right)dt} 
\\ & + \left( { - i} \right)^M \frac{{e^{\frac{\pi }{3}i} }}{{\sqrt 3 \pi }}\frac{{e^{2\pi i\kappa } }}{{2z^{M + \frac{2}{3}} }}\int_0^{ + \infty } {\frac{{t^{M - \frac{1}{3}} e^{ - 2\pi t} }}{{1 + it/z}}R_L^{\left( H \right)} \left( {it, - \kappa } \right)dt}  - i^M \frac{{e^{\frac{\pi }{3}i} }}{{\sqrt 3 \pi }}\frac{{e^{ - 2\pi i\kappa } }}{{2z^{M + \frac{2}{3}} }}\int_0^{ + \infty } {\frac{{t^{M - \frac{1}{3}} e^{ - 2\pi t} }}{{1 - it/z}}R_L^{\left( H \right)} \left( {it,\kappa } \right)dt} 
\end{align*}
for $\frac{3\pi}{2} < \arg z < \frac{5\pi}{2}$. The sum of the four integrals can be represented in terms of $R_{N,M,K,L}^{\left( H \right)}$ and by analytic continuation, we establish that
\begin{gather}\label{eq74}
\begin{split}
R_{N,M,K,L}^{\left( H \right)} \left( {z,\kappa } \right) =\; & ie^{ - \frac{\pi }{3}i} \frac{{e^{ - 2\pi i\nu } }}{{\sqrt 3 }}R_K^{\left( H \right)} \left( {ze^{ - \pi i} , - \kappa } \right) - ie^{\frac{\pi }{3}i} \frac{{e^{ - 2\pi i\nu } }}{{\sqrt 3 }}R_L^{\left( H \right)} \left( {ze^{ - \pi i} , - \kappa } \right)
\\ & +i\frac{{e^{2\pi i\nu } }}{{\sqrt 3 }}R_K^{\left( H \right)} \left( {z,\kappa } \right) - i\frac{{e^{2\pi i\nu } }}{{\sqrt 3 }}R_L^{\left( H \right)} \left( {z,\kappa } \right) - R_{N,M,K,L}^{\left( H \right)} \left( {ze^{ - 3\pi i} , - \kappa } \right)
\end{split}
\end{gather}
for $\frac{{3\pi }}{2} \le \arg z \le 3\pi  - \delta <3\pi$. If $K = L$, then we have
\[
\left| {R_{N,M,K,L}^{\left( H \right)} \left( {z,\kappa } \right)} \right| \le \frac{{e^{2\pi \Im \left( \nu  \right)} }}{{\sqrt 3 }}\left| {R_K^{\left( H \right)} \left( {ze^{ - \pi i} , - \kappa } \right)} \right| + \frac{{e^{2\pi \Im \left( \nu  \right)} }}{{\sqrt 3 }}\left| {R_L^{\left( H \right)} \left( {ze^{ - \pi i} , - \kappa } \right)} \right| + \left| {R_{N,M,K,L}^{\left( H \right)} \left( {ze^{ - 3\pi i} , - \kappa } \right)} \right|.
\]
Using the fact that $R_K^{\left( H \right)} \left( {ze^{-\pi i},-\kappa } \right) = \mathcal{O}_{\kappa ,K,\delta} \left( {\left|z\right|^{ - \frac{{K + 1}}{3}} } \right)$ and $R_L^{\left( H \right)} \left( {ze^{-\pi i},-\kappa } \right) = \mathcal{O}_{\kappa ,L,\delta} \left( {\left|z\right|^{ - \frac{{L + 1}}{3}} } \right)$ as $z\to \infty$ in $ \delta \leq \arg z \leq 3\pi  - \delta$, together with the estimates \eqref{eq54} and \eqref{eq73}, we deduce that
\[
R_{N,M,K,L}^{\left( H \right)} \left( {z,\kappa } \right) = \mathcal{O}_{\kappa ,K,\rho ,\delta } \left( {\frac{{e^{2\pi \Im \left( z \right)} }}{{\left| z \right|^{\frac{{K + 1}}{3}} }}} \right) + \mathcal{O}_{\kappa ,L,\sigma ,\delta } \left( {\frac{{e^{2\pi \Im \left( z \right)} }}{{\left| z \right|^{\frac{{L + 1}}{3}} }}} \right)
\]
when $z\to \infty$ in the sector $\frac{{3\pi }}{2} \le \arg z \le 3\pi  - \delta$. Otherwise, from the connection formula \eqref{eq75}, we have
\[
R_K^{\left( H \right)} \left( {z,\kappa } \right) =  - R_K^{\left( H \right)} \left( {ze^{ - 2\pi i} ,\kappa } \right) + R_K^{\left( H \right)} \left( {ze^{ - \pi i} , - \kappa } \right) - e^{ - 2\pi i\nu } H_\nu ^{\left( 2 \right)} \left( {ze^{ - 2\pi i} } \right)
\]
and similarly for $L$, whence \eqref{eq74} can be written as
\begin{align*}
R_{N,M,K,L}^{\left( H \right)} \left( {z,\kappa } \right) = \; & \frac{i}{{\sqrt 3 }}\left( {e^{ - \frac{\pi }{3}i} e^{ - 2\pi i\nu }  + e^{2\pi i\nu } } \right)R_K^{\left( H \right)} \left( {ze^{ - \pi i} , - \kappa } \right) - \frac{i}{{\sqrt 3 }}\left( {e^{\frac{\pi }{3}i} e^{ - 2\pi i\nu }  + e^{2\pi i\nu } } \right)R_L^{\left( H \right)} \left( {ze^{ - \pi i} , - \kappa } \right)
\\ & - i\frac{{e^{2\pi i\nu } }}{{\sqrt 3 }}R_K^{\left( H \right)} \left( {ze^{ - 2\pi i} ,\kappa } \right) + i\frac{{e^{2\pi i\nu } }}{{\sqrt 3 }}R_L^{\left( H \right)} \left( {ze^{ - 2\pi i} ,\kappa } \right) - R_{N,M,K,L}^{\left( H \right)} \left( {ze^{ - 3\pi i} , - \kappa } \right)
\end{align*}
for $\frac{{3\pi }}{2} \le \arg z \le 3\pi  - \delta <3\pi$. It is readily shown that this expression implies the estimate
\[
R_{N,M,K,L}^{\left( H \right)} \left( {z,\kappa } \right) = \mathcal{O}_{\kappa ,K,\rho ,\delta } \left( {\frac{{\cosh \left( {2\pi \Im \left( z \right)} \right)}}{{\left| z \right|^{\frac{{K + 1}}{3}} }}} \right) + \mathcal{O}_{\kappa ,L,\sigma ,\delta } \left( {\frac{{\cosh \left( {2\pi \Im \left( z \right)} \right)}}{{\left| z \right|^{\frac{{L + 1}}{3}} }}} \right)
\]
as $z\to \infty$ in the sector $\frac{{3\pi }}{2} \le \arg z \le 3\pi  - \delta$. The proof of the estimates for the sector $-2\pi < - 2\pi  + \delta  \le \arg z \le  - \frac{{3\pi }}{2}$ is similar.

To omit the assumption $\left| {\Re \left( \kappa  \right)} \right| < \min \left( {\frac{{K + 1}}{3},\frac{{L + 1}}{3}} \right)$, we proceed as follows. Let $K$ and $L$ be arbitrary non-negative integers and $\kappa$ be an arbitrary fixed complex number. Let $K'$ and $L'$ be non-negative integers such that $\left| {\Re \left( \kappa  \right)} \right| < \min \left( {\frac{{K' + 1}}{3},\frac{{L' + 1}}{3}} \right)$. We have
\begin{align*}
R_{N,M,K,L}^{\left( H \right)} \left( {z,\kappa } \right) = \; & ie^{ - \frac{\pi }{3}i} \frac{{e^{ - 2\pi i\nu } }}{{\sqrt 3 }}\frac{2}{{3\pi }}\sum\limits_{k = K}^{K' - 1} {6^{\frac{{k + 1}}{3}} B_k \left( \kappa  \right)\sin \left( {\frac{{\left( {k + 1} \right)\pi }}{3}} \right)\frac{{\Gamma \left( {\frac{{k + 1}}{3}} \right)}}{{z^{\frac{{k + 1}}{3}} }}\widehat T_{N - \frac{k}{3}} \left( { - 2\pi iz} \right)} 
\\ & - i\frac{{e^{2\pi i\nu } }}{{\sqrt 3 }}\frac{2}{{3\pi }}\sum\limits_{k = K}^{K' - 1} {6^{\frac{{k + 1}}{3}} B_k \left( \kappa  \right)e^{\frac{{2\left( {k + 1} \right)\pi i}}{3}} \sin \left( {\frac{{\left( {k + 1} \right)\pi }}{3}} \right)\frac{{\Gamma \left( {\frac{{k + 1}}{3}} \right)}}{{z^{\frac{{k + 1}}{3}} }}\widehat T_{N - \frac{k}{3}} \left( {2\pi iz} \right)} 
\\ & - ie^{\frac{\pi }{3}i} \frac{{e^{ - 2\pi i\nu } }}{{\sqrt 3 }}\frac{2}{{3\pi }}\sum\limits_{\ell  = L}^{L' - 1} {6^{\frac{{\ell  + 1}}{3}} B_\ell  \left( \kappa  \right)\sin \left( {\frac{{\left( {\ell  + 1} \right)\pi }}{3}} \right)\frac{{\Gamma \left( {\frac{{\ell  + 1}}{3}} \right)}}{{z^{\frac{{\ell  + 1}}{3}} }}\widehat T_{M - \frac{{\ell  - 1}}{3}} \left( { - 2\pi iz} \right)} 
\\ & + i\frac{{e^{2\pi i\nu } }}{{\sqrt 3 }}\frac{2}{{3\pi }}\sum\limits_{\ell  = L}^{L' - 1} {6^{\frac{{\ell  + 1}}{3}} B_\ell  \left( \kappa  \right)e^{\frac{{2\left( {\ell  + 1} \right)\pi i}}{3}} \sin \left( {\frac{{\left( {\ell  + 1} \right)\pi }}{3}} \right)\frac{{\Gamma \left( {\frac{{\ell  + 1}}{3}} \right)}}{{z^{\frac{{\ell  + 1}}{3}} }}\widehat T_{M - \frac{{\ell  - 1}}{3}} \left( {2\pi iz} \right)} 
\\ & + R_{N,M,K',L'}^{\left( H \right)} \left( {z,\kappa } \right).
\end{align*}
Employing the previously obtained bounds for $R_{N,M,K',L'}^{\left( H \right)} \left( {z,\kappa } \right)$ and Olver's estimation \eqref{eq55} together with the connection formula for the Terminant function \cite[p. 260]{Paris3}, shows that $R_{N,M,K,L}^{\left( H \right)} \left( {z,\kappa } \right)$ indeed satisfies the order estimates given in Theorem \ref{thm3}.

The corresponding expansion for $H_\nu ^{\left( 1 \right) \prime} \left( z \right)$ follows from the relation $2H_\nu ^{\left( 1 \right) \prime} \left( z \right) = H_{\nu  - 1}^{\left( 1 \right)} \left( z \right) - H_{\nu  + 1}^{\left( 1 \right)} \left( z \right)$.

\subsection{Stokes phenomenon and Berry's transition}

We study the Stokes phenomenon related to the asymptotic expansion of $H_{\nu}^{\left( 1 \right)} \left( z \right)$ occurring when $\arg z$ passes through the value $-\frac{\pi}{2}$. The Stokes line $\arg z =\frac{3\pi}{2}$ and the asymptotic expansion of $H_{\nu}^{\left( 1 \right)\prime} \left( z \right)$ can be treated similarly. From \eqref{eq62} we have
\begin{align*}
& R_{3N}^{\left( H \right)} \left( {z,\kappa } \right) = R_{N,N,0,0}^{\left( H \right)} \left( {z,\kappa } \right) = ie^{\frac{\pi }{3}i} \frac{{e^{ - 2\pi i\nu } }}{{\sqrt 3 }}H_\nu ^{\left( 2 \right)} \left( z \right) - ie^{ - \frac{\pi }{3}i} \frac{{e^{ - 2\pi i\nu } }}{{\sqrt 3 }}H_\nu ^{\left( 2 \right)} \left( z \right) \\ & + \left( { - i} \right)^N \frac{{e^{ - \frac{\pi }{3}i} }}{{\sqrt 3 \pi }}\frac{{ie^{2\pi i\kappa } }}{{2z^{N + \frac{1}{3}} }}\int_0^{ + \infty } {\frac{{t^{N - \frac{2}{3}} e^{ - 2\pi t} }}{{1 + it/z}}H_{it + \kappa }^{\left( 1 \right)} \left( {it} \right)dt}  + i^N \frac{{e^{ - \frac{\pi }{3}i} }}{{\sqrt 3 \pi }}\frac{{ie^{ - 2\pi i\kappa } }}{{2z^{N + \frac{1}{3}} }}\int_0^{ + \infty } {\frac{{t^{N - \frac{2}{3}} e^{ - 2\pi t} }}{{1 - it/z}}H_{it - \kappa }^{\left( 1 \right)} \left( {it} \right)dt} 
\\ & + \left( { - i} \right)^N \frac{{e^{\frac{\pi }{3}i} }}{{\sqrt 3 \pi }}\frac{{e^{2\pi i\kappa } }}{{2z^{N + \frac{2}{3}} }}\int_0^{ + \infty } {\frac{{t^{N - \frac{1}{3}} e^{ - 2\pi t} }}{{1 + it/z}}H_{it + \kappa }^{\left( 1 \right)} \left( {it} \right)dt}  - i^N \frac{{e^{\frac{\pi }{3}i} }}{{\sqrt 3 \pi }}\frac{{e^{ - 2\pi i\kappa } }}{{2z^{N + \frac{2}{3}} }}\int_0^{ + \infty } {\frac{{t^{N - \frac{1}{3}} e^{ - 2\pi t} }}{{1 - it/z}}H_{it - \kappa }^{\left( 1 \right)} \left( {it} \right)dt} 
\end{align*}
when $-\frac{3\pi}{2} < \arg z < -\frac{\pi}{2}$. We can expand the integrals into an asymptotic series in inverse powers of $z$, using the geometric series and the second representation in \eqref{eq8}. For the Hankel function $H_\nu ^{\left( 2 \right)} \left( z \right)$, we use \eqref{eq18}. The final result is the compound asymptotic series
\begin{align*}
R_{3N}^{\left( H \right)} \left( {z,\kappa } \right) \sim &  - ie^{\frac{\pi }{3}i} \frac{{e^{ - 2\pi i\nu } }}{{\sqrt 3 }}\frac{2}{{3\pi }}\sum\limits_{k = 0}^\infty  {6^{\frac{{k + 1}}{3}} B_k \left( \kappa  \right)e^{ - \frac{{2\left( {k + 1} \right)\pi i}}{3}} \sin \left( {\frac{{\left( {k + 1} \right)\pi }}{3}} \right)\frac{{\Gamma \left( {\frac{{k + 1}}{3}} \right)}}{{z^{\frac{{k + 1}}{3}} }}} 
\\ & + ie^{ - \frac{\pi }{3}i} \frac{{e^{ - 2\pi i\nu } }}{{\sqrt 3 }}\frac{2}{{3\pi }}\sum\limits_{\ell  = 0}^\infty  {6^{\frac{{\ell  + 1}}{3}} B_\ell  \left( \kappa  \right)e^{ - \frac{{2\left( {\ell  + 1} \right)\pi i}}{3}} \sin \left( {\frac{{\left( {\ell  + 1} \right)\pi }}{3}} \right)\frac{{\Gamma \left( {\frac{{\ell  + 1}}{3}} \right)}}{{z^{\frac{{\ell  + 1}}{3}} }}} 
\\ & - \frac{2}{{3\pi }}\sum\limits_{n = 3N}^\infty  {6^{\frac{{n + 1}}{3}} B_n \left( \kappa  \right)e^{\frac{{2\left( {n + 1} \right)\pi i}}{3}} \sin \left( {\frac{{\left( {n + 1} \right)\pi }}{3}} \right)\frac{{\Gamma \left( {\frac{{n + 1}}{3}} \right)}}{{z^{\frac{{n + 1}}{3}} }}} 
\end{align*}
as $z\to\infty$ in the sector $-\frac{3\pi}{2} < \arg z < -\frac{\pi}{2}$. Hence, as $\arg z$ decreases through the value $-\frac{\pi}{2}$, the two additional series
\begin{gather}\label{eq63}
\begin{split}
 &  - ie^{\frac{\pi }{3}i} \frac{{e^{ - 2\pi i\nu } }}{{\sqrt 3 }}\frac{2}{{3\pi }}\sum\limits_{k = 0}^\infty  {6^{\frac{{k + 1}}{3}} B_k \left( \kappa  \right)e^{ - \frac{{2\left( {k + 1} \right)\pi i}}{3}} \sin \left( {\frac{{\left( {k + 1} \right)\pi }}{3}} \right)\frac{{\Gamma \left( {\frac{{k + 1}}{3}} \right)}}{{z^{\frac{{k + 1}}{3}} }}} 
\\ & + ie^{ - \frac{\pi }{3}i} \frac{{e^{ - 2\pi i\nu } }}{{\sqrt 3 }}\frac{2}{{3\pi }}\sum\limits_{\ell  = 0}^\infty  {6^{\frac{{\ell  + 1}}{3}} B_\ell  \left( \kappa  \right)e^{ - \frac{{2\left( {\ell  + 1} \right)\pi i}}{3}} \sin \left( {\frac{{\left( {\ell  + 1} \right)\pi }}{3}} \right)\frac{{\Gamma \left( {\frac{{\ell  + 1}}{3}} \right)}}{{z^{\frac{{\ell  + 1}}{3}} }}}
\end{split}
\end{gather}
appear in the asymptotic expansion of $H_{\nu}^{\left( 1 \right)} \left( z \right)$ beside the original series \eqref{eq17}. We have encountered a Stokes phenomenon with Stokes line $\arg z = -\frac{\pi}{2}$.

In the important papers \cite{Berry3, Berry2}, Berry provided a new interpretation of the Stokes phenomenon; he found that assuming optimal truncation, the transition between compound asymptotic expansions is of Error function type, thus yielding a smooth, although very rapid, transition as a Stokes line is crossed.

Using the exponentially improved expansion given in Theorem \ref{thm3}, we show that the Nicholson--Debye expansion for $H_{\nu}^{\left( 1 \right)} \left( z \right)$ exhibits the Berry transition between the asymptotic series across the Stokes line $\arg z = -\frac{\pi}{2}$. More precisely, we shall find that the first few terms of the two series in \eqref{eq63} ``emerge" in a rapid and smooth way as $\arg z$ decreases through $-\frac{\pi}{2}$.

Let us assume that $M,N \approx 2\pi \left|z\right|$. Under these conditions, Olver's estimation \eqref{eq55} gives that
\[
e^{2\pi i\nu } \widehat T_{N - \frac{k}{3}} \left( {2\pi iz} \right),e^{2\pi i\nu } \widehat T_{M - \frac{{\ell  - 1}}{3}} \left( {2\pi iz} \right) = \mathcal{O}\left( {e^{ - 2\pi \left| z \right|} } \right)
\]
as $z\to \infty$ and $-\frac{3\pi}{2} < \arg z < \frac{\pi}{2}$. Therefore, from Theorem \ref{thm3}, we infer that for large $z$, $-\pi <\arg z<\frac{\pi}{2}$, we have
\begin{align*}
& H_\nu ^{\left( 1 \right)} \left( z \right) \approx  \frac{{e^{ - \frac{\pi }{3}i} }}{{\sqrt 3 \pi z^{\frac{1}{3}} }}\sum\limits_{n = 0}^{N - 1} {\left( { - 1} \right)^n 6^{n + \frac{1}{3}} B_{3n} \left( \kappa  \right)\frac{{\Gamma \left( {n + \frac{1}{3}} \right)}}{{z^n }}}  + \frac{{e^{\frac{\pi }{3}i} }}{{\sqrt 3 \pi z^{\frac{2}{3}} }}\sum\limits_{m = 0}^{M - 1} {\left( { - 1} \right)^m 6^{m + \frac{2}{3}} B_{3m + 1} \left( \kappa  \right)\frac{{\Gamma \left( {m + \frac{2}{3}} \right)}}{{z^m }}}
\\ & - ie^{\frac{\pi }{3}i} \frac{{e^{ - 2\pi i\nu } }}{{\sqrt 3 }}\frac{2}{{3\pi }}\sum\limits_{k = 0} {6^{\frac{{k + 1}}{3}} B_k \left( \kappa  \right)e^{ - \frac{{2\left( {k + 1} \right)\pi i}}{3}} \sin \left( {\frac{{\left( {k + 1} \right)\pi }}{3}} \right)\frac{{\Gamma \left( {\frac{{k + 1}}{3}} \right)}}{{z^{\frac{{k + 1}}{3}} }}\left( { - e^{2\pi i\frac{k}{3}} \widehat T_{N - \frac{k}{3}} \left( { - 2\pi iz} \right)} \right)} 
\\ & + ie^{ - \frac{\pi }{3}i} \frac{{e^{ - 2\pi i\nu } }}{{\sqrt 3 }}\frac{2}{{3\pi }}\sum\limits_{\ell  = 0} {6^{\frac{{\ell  + 1}}{3}} B_\ell  \left( \kappa  \right)e^{ - \frac{{2\left( {\ell  + 1} \right)\pi i}}{3}} \sin \left( {\frac{{\left( {\ell  + 1} \right)\pi }}{3}} \right)\frac{{\Gamma \left( {\frac{{\ell  + 1}}{3}} \right)}}{{z^{\frac{{\ell  + 1}}{3}} }}\left( { - e^{2\pi i\frac{{\ell  - 1}}{3}} \widehat T_{M - \frac{{\ell  - 1}}{3}} \left( { - 2\pi iz} \right)} \right)} ,
\end{align*}
where $\sum\nolimits_{k = 0}$ and $\sum\nolimits_{\ell = 0}$ mean that the sums are restricted to the first few terms of the series.

Under the above assumptions on $N$ and $M$, from \eqref{eq57} and \eqref{eq58}, the normalised Terminant functions have the asymptotic behaviour
\[
 - e^{2\pi i\frac{k}{3}} \widehat T_{N - \frac{k}{3}} \left( { - 2\pi iz} \right), - e^{2\pi i\frac{{\ell  - 1}}{3}} \widehat T_{M - \frac{{\ell  - 1}}{3}} \left( { - 2\pi iz} \right) \sim \frac{1}{2} - \frac{1}{2}\mathop{\text{erf}}\left( {\left( {\theta  + \frac{\pi }{2}} \right)\sqrt {\pi \left| z \right|} } \right)
\]
provided that $\arg z = \theta$ is close to $-\frac{\pi}{2}$, $z$ is large and $k,\ell$ are small in comparison with $N$ and $M$. Therefore, when $\theta  >  - \frac{\pi }{2}$, the normalised Terminant functions are exponentially small; for $\theta  =  -\frac{\pi}{2}$, they are asymptotically $\frac{1}{2}$ up to an exponentially small error; and when $\theta  <  - \frac{\pi}{2}$, the normalised Terminant functions are asymptotic to $1$ with an exponentially small error. Thus, the transition across the Stokes line $\arg z = -\frac{\pi}{2}$ is effected rapidly and smoothly.

\section{Discussion}\label{section6}

In this paper, we have discussed in detail the asymptotic expansions of the Hankel functions $H_\nu ^{\left( 1 \right)} \left( {z } \right)$, $H_\nu ^{\left( 2 \right)} \left( {z } \right)$, the Bessel functions $J_\nu  \left( {z} \right)$, $Y_\nu  \left( {z} \right)$ and their derivatives of nearly equal order and argument. We have obtained explicit formulas for the remainder terms of these expansions, and have given asymptotic expansions for their late coefficients. In particular, we have solved a problem of Watson regarding the approximation of the late coefficients. We have derived exponentially improved asymptotic expansions, and have demonstrated the smooth transition of the Stokes discontinuities.

For the special case, when the oder equals to the argument, we have given realistic error bounds for all the asymptotic series. In this special case the error analysis is simple and the bounds are directly related to the absolute value of the first non-vanishing omitted term(s) of the asymptotic series. It seems that in general, it is not possible to derive such simple bounds, because the functions $i H_{it \pm \kappa }^{\left( 1 \right)}\left( {it} \right)$ and $H_{it \pm \kappa }^{\left( 1 \right)\prime}\left( {it} \right)$ ($t>0$) are no longer positive. Nevertheless, we show here how to derive computable bounds for $R_N^{\left( {H } \right)} \left( {z ,\kappa } \right)$ when $-\frac{\pi}{2}<\arg z <\frac{3\pi}{2}$. For simplicity, we consider the special case when $N \equiv 1 \mod 6$. In this case, the representation \eqref{eq9} can be rearranged into the more useful form
\begin{align*}
R_N^{\left( {H } \right)} \left( {z ,\kappa } \right) = \; & \frac{{i^{N + 1} }}{{2\sqrt 3 \pi z^{\frac{{N + 1}}{3}} }}\int_0^{ + \infty } {\frac{{t^{\frac{{N - 2}}{3}} e^{ - 2\pi t} e^{\frac{{2\left( {N + 1} \right)\pi i}}{3}} }}{{\left( {1 + \left( {t/z} \right)^{\frac{2}{3}} e^{\frac{{2\pi i}}{3}} } \right)\left( {1 + \left( {t/z} \right)^{\frac{2}{3}} } \right)}}\left( {e^{2\pi i\kappa } H_{it + \kappa }^{\left( 1 \right)} \left( {it} \right) - e^{ - 2\pi i\kappa } H_{it - \kappa }^{\left( 1 \right)} \left( {it} \right)} \right)dt} 
\\ & + \frac{{i^N }}{{2\sqrt 3 \pi z^{\frac{{N + 4}}{3}} }}\int_0^{ + \infty } {\frac{{t^{\frac{{N + 1}}{3}} e^{ - 2\pi t} e^{\frac{{2\left( {N + 1} \right)\pi i}}{3}} }}{{\left( {1 + \left( {t/z} \right)^{\frac{2}{3}} e^{\frac{{2\pi i}}{3}} } \right)\left( {1 + \left( {t/z} \right)^{\frac{2}{3}} } \right)}}\left( {e^{2\pi i\kappa } H_{it + \kappa }^{\left( 1 \right)} \left( {it} \right) + e^{ - 2\pi i\kappa } H_{it - \kappa }^{\left( 1 \right)} \left( {it} \right)} \right)dt} ,
\end{align*}
provided that $\left|\Re \left( \kappa  \right)\right| < \frac{N + 1}{3}$. Such rearrangements are possible for the other cases too, but the error term may consist of four different integrals. For this form we can use trivial estimation and the inequality \eqref{eq36}, to obtain the bound
\begin{align*}
\left| {R_N^{\left( H \right)} \left( {z,\kappa } \right)} \right| \le \; & \left( {\frac{1}{{2\sqrt 3 \pi \left| z \right|^{\frac{{N + 1}}{3}} }}\int_0^{ + \infty } {t^{\frac{{N - 2}}{3}} e^{ - 2\pi t} \left| {e^{2\pi i\kappa } H_{it + \kappa }^{\left( 1 \right)} \left( {it} \right) - e^{ - 2\pi i\kappa } H_{it - \kappa }^{\left( 1 \right)} \left( {it} \right)} \right|dt} }\right.\\ & \left.{+ \frac{1}{{2\sqrt 3 \pi \left| z \right|^{\frac{{N + 4}}{3}} }}\int_0^{ + \infty } {t^{\frac{{N + 1}}{3}} e^{ - 2\pi t} \left| {e^{2\pi i\kappa } H_{it + \kappa }^{\left( 1 \right)} \left( {it} \right) + e^{ - 2\pi i\kappa } H_{it - \kappa }^{\left( 1 \right)} \left( {it} \right)} \right|dt} } \right) \\ & \times \begin{cases} \left|\sec \theta \right| & \; \text{ if } \; { - \frac{\pi }{2} < \theta  < 0 \; \text{ or } \; \pi  < \theta  < \frac{{3\pi }}{2}}, \\ 1 & \; \text{ if } \; {0 \le \theta  \le \pi }. \end{cases}
\end{align*}
The convergence of the integrals follows from the argument in Appendix \ref{appendixb}. Numerical evaluation of these integrals then leads to explicit error bounds. It is seen from the second representation in \eqref{eq8} that the terms in this bound are closely related to the absolute value of the first two non-vanishing omitted terms of the asymptotic series \eqref{eq17}, though they are not identical to them.

\appendix

\section{Computation of the coefficients $B_n\left(\kappa\right)$}\label{appendixa}

In 1952, Lauwerier \cite{Lauwerier} showed that the coefficients in asymptotic expansions of Laplace-type integrals can be calculated by means of linear recurrence relations. Simple application of his method provides the formula
\[
B_n \left( \kappa  \right) = \frac{{6^{ - \frac{{n + 1}}{3}} }}{{\Gamma \left( {\frac{{n + 1}}{3}} \right)}}\int_0^{ + \infty } {t^{\frac{{n - 2}}{3}} e^{ - \frac{t}{6}} P_n \left( {t,\kappa } \right)dt} ,
\]
where the two-variable polynomials $P_0 \left( x,\kappa \right), P_1 \left( x,\kappa \right), P_2 \left( x,\kappa \right),\ldots$ are given by the recurrence relation
\[
P_n \left( {x,\kappa } \right) = \frac{{\kappa ^n }}{{n!}} - \sum\limits_{k = 1}^{\left\lfloor {\frac{n}{2}} \right\rfloor } {\frac{1}{{\left( {2k + 3} \right)!}}\int_0^x {P_{n - 2k} \left( {t,\kappa } \right)dt} } \; \text{ for } \; n\geq 2,
\]
with $P_0 \left( {x,\kappa } \right) = 1$ and $P_1 \left( {x,\kappa } \right) = \kappa $.

Expanding the higher derivative in \eqref{eq8} using Leibniz's rule and noting that $t^3/\left(\sinh t -t\right)$ is an even function of $t$, we deduce
\begin{equation}\label{eq45}
B_n \left( \kappa  \right) = \sum\limits_{k = 0}^{\left\lfloor {\frac{n}{2}} \right\rfloor } {\frac{{\kappa ^{n - 2k} }}{{\left( {n - 2k} \right)!}}\frac{1}{{\left( {2k} \right)!}}\left[ {\frac{{d^{2k} }}{{dt^{2k} }}\left( {\frac{1}{6}\frac{{t^3 }}{{\sinh t - t}}} \right)^{\frac{{n + 1}}{3}} } \right]_{t = 0} } ,
\end{equation}
which is Airey's \cite{Airey2} representation for the coefficients $B_n\left(\kappa\right)$.

There has been a recent interest in finding explicit formulas for the coefficients in asymptotic expansions of Laplace-type integrals (see \cite{Lopez}, \cite{Nemes2}, \cite{Wojdylo1} and \cite{Wojdylo2}). There are two general formulas for these coefficients, one containing Potential polynomials and one containing Bell polynomials. We derive them here for the special case of the coefficients $B_n\left(\kappa\right)$. Let
\[
\sinh t - t = \sum\limits_{j= 0}^\infty  {a_j t^{j + 3} } ,
\]
so that
\[
a_{2j}  = \frac{1}{{\left( {2j + 3} \right)!}},\; a_{2j + 1}  = 0 \; \text{ for } \; j \ge 0.
\]
Let $0 \leq i \leq j$ be integers and $\rho$ be a complex number. We define the Potential polynomials
\[
\mathsf{A}_{\rho ,j}  = \mathsf{A}_{\rho ,j} \left( {\frac{{a_1 }}{{a_0 }},\frac{{a_2 }}{{a_0 }}, \ldots ,\frac{{a_k }}{{a_0 }}} \right)
\]
and the Bell polynomials
\[
\mathsf{B}_{j,i}  = \mathsf{B}_{j,i} \left( {a_1 ,a_2 , \ldots ,a_{j - i + 1} } \right)
\]
via the expansions
\begin{equation}\label{eq43}
\left( {1 + \sum\limits_{j = 1}^\infty  {\frac{{a_j }}{{a_0 }}t^j } } \right)^\rho  = \sum\limits_{j = 0}^\infty  {\mathsf{A}_{\rho ,j} t^j } \; \text{ and } \; \mathsf{A}_{\rho ,j}  = \sum\limits_{i = 0}^j {\binom{\rho}{i}\frac{1}{a_0^i}\mathsf{B}_{j,i}} .
\end{equation}
Naturally, these polynomials can be defined for arbitrary power series with $a_0\neq 0$. It is possible to express the Potential polynomials with complex parameter in terms of Potential polynomials with integer parameter using the following formula of Comtet \cite[p. 142]{Comtet}
\begin{equation}\label{eq44}
\mathsf{A}_{\rho ,j}  = \frac{\Gamma \left( { - \rho  + j + 1} \right)}{j!\Gamma \left( { - \rho } \right)}\sum\limits_{i = 0}^j {\frac{\left( { - 1} \right)^i}{- \rho  + i}\binom{j}{i}\mathsf{A}_{i,j} } .
\end{equation}
With these notations we can write \eqref{eq45} as
\[
B_n \left( \kappa  \right) = \sum\limits_{k = 0}^{\left\lfloor {\frac{n}{2}} \right\rfloor } {\frac{{\kappa ^{n - 2k} }}{{\left( {n - 2k} \right)!}}\mathsf{A}_{ - \frac{{n + 1}}{3},2k} } .
\]
Using \eqref{eq43} and \eqref{eq44} we find
\[
B_n \left( \kappa  \right) = \sum\limits_{k = 0}^{\left\lfloor {\frac{n}{2}} \right\rfloor } {\frac{{\kappa ^{n - 2k} }}{{\left( {n - 2k} \right)!}}\sum\limits_{j = 0}^{2k} {\left( { - 1} \right)^j 6^j \frac{{\Gamma \left( {\frac{{n + 1}}{3} + j} \right)}}{{j!\Gamma \left( {\frac{{n + 1}}{3}} \right)}}\mathsf{B}_{2k,j}}  } 
\]
and
\begin{equation}\label{eq46}
B_n \left( \kappa  \right) = \sum\limits_{k = 0}^{\left\lfloor {\frac{n}{2}} \right\rfloor } {\frac{{\kappa ^{n - 2k} }}{{\left( {n - 2k} \right)!}}\frac{{3\Gamma \left( {\frac{{n + 1}}{3} + 2k + 1} \right)}}{{\left( {2k} \right)!\Gamma \left( {\frac{{n + 1}}{3}} \right)}}\sum\limits_{j = 0}^{2k} {\frac{{\left( { - 1} \right)^j }}{{n + 3j + 1}} \binom{2k}{j} \mathsf{A}_{j,2k} } } .
\end{equation}
The quantities $\mathsf{B}_{k,j}$ and $\mathsf{A}_{j,k}$ appearing in these formulas may be computed from the recurrence relations
\[
\mathsf{B}_{k,j}  = \sum\limits_{i= 1}^{k - j + 1} {a_i \mathsf{B}_{k - i,j - 1} } \; \text{ and } \; \mathsf{A}_{j,k}  = \sum\limits_{i = 0}^{k} {\frac{a_i}{a_0}\mathsf{A}_{j - 1,k - i} }
\]
with $\mathsf{B}_{0,0} = \mathsf{A}_{0,0} = 1$, $\mathsf{B}_{i,0} = \mathsf{A}_{0,i} = 0$ $\left(i \ge 1\right)$, $\mathsf{B}_{i,1} = a_0 \mathsf{A}_{1,i} = a_i$ (see Nemes \cite{Nemes2}). In the paper \cite{Nemes}, it was shown that the Potential polynomials $\mathsf{A}_{j,2k}$ in \eqref{eq46} can be written as
\[
\mathsf{A}_{j,2k}  = \sum\limits_{i = 0}^j {\left( { - 1} \right)^{j + i} \binom{j}{i}\frac{{2^{2k + 2j} 6^j }}{{\left( {2k + 2j} \right)!}}B_{2k + 2j}^{\left( { - i} \right)} \left( { - \frac{i}{2}} \right)} ,
\]
where $B_n^{\left( \ell  \right)} \left(x\right)$ stands for the generalised Bernoulli polynomials, which are defined by the exponential generating function
\[
\left( \frac{t}{e^t  - 1} \right)^\ell  e^{x t}  = \sum\limits_{n = 0}^\infty  {B_n^{\left( \ell  \right)} \left(x\right)\frac{t^n}{n!}} \; \text{ for } \; \left|t\right| < 2\pi.
\]
For basic properties of these polynomials, see Milne-Thomson \cite{Milne-Thomson} or N\"{o}rlund \cite{Norlund}.

Some other properties of the polynomials $B_n \left( \kappa  \right)$ can be found in the paper of Sch\"{o}be \cite{Schobe}.

\section{Auxiliary estimates}\label{appendixb} In this appendix, we prove some estimates for the remainder $R_N^{\left( {H} \right)} \left( {z,\kappa } \right)$. First, we prove the estimate \eqref{eq52}. Suppose that $0 \leq \arg z \leq \pi$ and $\left|\Re \left( \kappa  \right)\right| < \frac{N + 1}{3}$ with $N\geq 0$. Trivial estimation of \eqref{eq9} yields
\begin{align*}
\left| {R_N^{\left( {H} \right)} \left( {z,\kappa } \right)} \right| \le \; & \frac{{\left| {e^{2\pi i\kappa } } \right|}}{{6\pi \left| z \right|^{\frac{{N + 1}}{3}} }}\int_0^{ + \infty } {t^{\frac{{N - 2}}{3}} e^{ - 2\pi t} \left| {\frac{{e^{\frac{{\left( {N + 1} \right)\pi i}}{3}} }}{{1 + i\left( {t/z} \right)^{\frac{1}{3}} e^{\frac{\pi }{3}i} }} - \frac{{e^{\left( {N + 1} \right)\pi i} }}{{1 - i\left( {t/z} \right)^{\frac{1}{3}} }}} \right|\left| {H_{it + \kappa }^{\left( 1 \right)} \left( {it} \right)} \right|dt} 
\\ & + \frac{{\left| {e^{ - 2\pi i\kappa } } \right|}}{{6\pi \left| z \right|^{\frac{{N + 1}}{3}} }}\int_0^{ + \infty } {t^{\frac{{N - 2}}{3}} e^{ - 2\pi t} \left| {\frac{{e^{\frac{{\left( {N + 1} \right)\pi i}}{3}} }}{{1 - i\left( {t/z} \right)^{\frac{1}{3}} e^{\frac{\pi }{3}i} }} - \frac{{e^{\left( {N + 1} \right)\pi i} }}{{1 + i\left( {t/z} \right)^{\frac{1}{3}} }}} \right|\left| {H_{it - \kappa }^{\left( 1 \right)} \left( {it} \right)} \right|dt} .
\end{align*}
Employing the elementary inequality
\begin{equation}\label{eq50}
\frac{1}{{\left| {1 - re^{i\vartheta } } \right|}} \le \begin{cases} \left|\csc \vartheta \right| & \; \text{ if } \; 0 < \left|\vartheta \text{ mod } 2\pi\right| <\frac{\pi}{2}, \\ 1 & \; \text{ if } \; \frac{\pi}{2} \leq \left|\vartheta \text{ mod } 2\pi\right| \leq \pi, \end{cases}
\end{equation}
for $r>0$, we obtain
\[
\left| {\frac{{e^{\frac{{\left( {N + 1} \right)\pi i}}{3}} }}{{1 + i\left( {t/z} \right)^{\frac{1}{3}} e^{\frac{\pi }{3}i} }} - \frac{{e^{\left( {N + 1} \right)\pi i} }}{{1 - i\left( {t/z} \right)^{\frac{1}{3}} }}} \right| \le \frac{1}{{\left| {1 + i\left( {t/z} \right)^{\frac{1}{3}} e^{\frac{\pi }{3}i} } \right|}} + \frac{1}{{\left| {1 - i\left( {t/z} \right)^{\frac{1}{3}} } \right|}} \le 2 + 2 = 4
\]
and
\[
\left| {\frac{{e^{\frac{{\left( {N + 1} \right)\pi i}}{3}} }}{{1 - i\left( {t/z} \right)^{\frac{1}{3}} e^{\frac{\pi }{3}i} }} - \frac{{e^{\left( {N + 1} \right)\pi i} }}{{1 + i\left( {t/z} \right)^{\frac{1}{3}} }}} \right| \le \frac{1}{{\left| {1 - i\left( {t/z} \right)^{\frac{1}{3}} e^{\frac{\pi }{3}i} } \right|}} + \frac{1}{{\left| {1 + i\left( {t/z} \right)^{\frac{1}{3}} } \right|}} \le 1 + 1 = 2.
\]
Hence, we have the estimate
\[
\left| {R_N^{\left( {H} \right)} \left( {z,\kappa } \right)} \right| \le \frac{{2\left| {e^{2\pi i\kappa } } \right|}}{{3\pi \left| z \right|^{\frac{{N + 1}}{3}} }}\int_0^{ + \infty } {t^{\frac{{N - 2}}{3}} e^{ - 2\pi t} \left| {H_{it + \kappa }^{\left( 1 \right)} \left( {it} \right)} \right|dt}  + \frac{{\left| {e^{ - 2\pi i\kappa } } \right|}}{{3\pi \left| z \right|^{\frac{{N + 1}}{3}} }}\int_0^{ + \infty } {t^{\frac{{N - 2}}{3}} e^{ - 2\pi t} \left| {H_{it - \kappa }^{\left( 1 \right)} \left( {it} \right)} \right|dt} .
\]
It remains to show that the integrals on the right-hand side are convergent. The integrands are continuous functions of $t>0$. The asymptotic expansion \eqref{eq17} shows that for large positive $t$ and fixed $\kappa$, we have
\[
H_{it \pm \kappa }^{\left( 1 \right)} \left( {it} \right)  \sim \frac{2}{{3\pi }}6^{\frac{1}{3}} \frac{{\sqrt 3 }}{2}\frac{{\Gamma \left( {\frac{1}{3}} \right)}}{{t^{\frac{1}{3}} }}. 
\]
We also have
\[
\left| { H_{it \pm \kappa }^{\left( 1 \right)} \left( {it} \right)} \right| = \begin{cases} \mathcal{O}_\kappa  \left( {t^{ - \left| {\Re \left( \kappa  \right)} \right|} } \right) & \; \text{ if } \; \kappa \neq 0, \\ \mathcal{O}\left( {\log t} \right) & \; \text{ if } \; \kappa = 0, \end{cases}
\]
as $t\to 0$ along the positive real axis (see, e.g., \cite[\S 10.7(i)]{NIST}). Therefore the integrals are convergent and we conclude that there is a constant $C_N \left( \kappa  \right)>0$ depending only on $N$ and $\kappa$, such that if $0 \leq \arg z \leq \pi$ and $\left|\Re \left( \kappa  \right)\right| < \frac{N + 1}{3}$, then
\[
\left| {R_N^{\left( {H} \right)} \left( {z,\kappa} \right)} \right| \le \frac{{C_N \left( \kappa \right)}}{{\left| z \right|^{\frac{{N + 1}}{3}} }}.
\]

Finally, we extend the region of validity of the expansions \eqref{eq17} and \eqref{eq18}. It was given by Watson that the $N$th remainders in the asymptotic expansions \eqref{eq17} and \eqref{eq18} satisfy
\[
R_N^{\left( H \right)} \left( {z,\kappa } \right) = \mathcal{O}_{\kappa ,N,\delta} \left( {\frac{1}{{\left| z \right|^{\frac{{N + 1}}{3}} }}} \right) \; \text{ and } \; - R_N^{\left( H \right)} \left( {ze^{\pi i} , - \kappa } \right) = \mathcal{O}_{\kappa ,N,\delta} \left( {\frac{1}{{\left| z \right|^{\frac{{N + 1}}{3}} }}} \right)
\]
as $z\to \infty$ in $\left| {\arg z} \right| \le \pi - \delta  < \pi$, with $0<\delta \leq \pi$ being fixed. However, from the connection between the two remainders it is seen that the first estimate is valid in the larger sector $ - \pi  + \delta \leq \arg z \leq 2\pi  - \delta$, and the second estimate is valid in the larger sector $ - 2\pi  + \delta \leq \arg z \leq \pi  - \delta$.

\end{document}